\documentclass{m2an}
\usepackage{amssymb}
\usepackage{stmaryrd}
\usepackage{mathrsfs}
\usepackage{mathtools}
\usepackage{enumerate}
\usepackage{hyperref}
\usepackage{esint}
\usepackage{subfigure}
\usepackage{booktabs}
\usepackage{psfrag}

\newcommand{\R}{\mathbb{R}}
\newcommand{\ysf}{\mathsf{y}}
\newcommand{\LAP}{{\Delta}}
\newcommand{\usf}{\mathsf{u}}
\newcommand{\asf}{\mathsf{a}}
\newcommand{\bsf}{\mathsf{b}}
\newcommand{\T}{\mathscr{T}}
\newcommand{\Xcal}{\mathcal{X}}
\newcommand{\Ycal}{\mathcal{Y}}
\newcommand{\dist}{ {\textup{\textsf{d}}}_{x_0} }
\newcommand{\distz}{ {\textup{\textsf{d}}}_{z} }
\newcommand{\vsf}{\mathsf{v}}
\newcommand{\psf}{\mathsf{p}}
\newcommand{\wsf}{\mathsf{w}}
\newcommand{\Tr}{\mathbb{T}}
\newcommand{\V}{\mathbb{V}}
\newcommand{\bu}{\bar{\mathbf{u}}}
\newcommand{\Sides}{\mathscr{S}}
\newcommand{\Ne}{\mathcal{N}}
\newcommand{\E}{\mathscr{E}}
\newcommand{\calG}{{\mathcal{G}}}
\newcommand{\GRAD}{\nabla}
\newcommand{\diff}{\, \mbox{\rm d}}
\newcommand{\qsf}{\mathsf{q}}
\newcommand{\rsf}{\mathsf{r}}
\newcommand{\Res}{\mathscr{R}}
\newcommand{\Orara}{\mathscr{O}}

\DeclareMathOperator*{\diam}{diam}
\DeclareMathOperator*{\supp}{supp}

\usepackage[usenames,dvipsnames]{color}

\begin{document}

\title{An a posteriori error analysis for an optimal control problem with point sources}\thanks{AA is partially supported by the USM grant 116.12.1. EO is partially supported by CONICYT through FONDECYT project 3160201. RR is supported by BASAL PFB03 CMM project, Universidad de Chile. AJS is partially supported by NSF grant DMS-1418784.}
\author{Alejandro Allendes}\address{Departamento de Matem\'atica, Universidad T\'ecnica Federico Santa Mar\'ia, Valpara\'iso, Chile.
\texttt{alejandro.allendes@usm.cl}}
\author{Enrique Ot\'arola}\address{Departamento de Matem\'atica, Universidad T\'ecnica Federico Santa Mar\'ia, Valpara\'iso, Chile.
\texttt{enrique.otarola@usm.cl}}
\author{Richard Rankin}\address{Departamento de Matem\'atica, Universidad T\'ecnica Federico Santa Mar\'ia, Valpara\'iso, Chile.
\texttt{richard.rankin@usm.cl}}
\author{Abner J.~Salgado}\address{Department of Mathematics, University of Tennessee, Knoxville, TN 37996, USA. \texttt{asalgad1@utk.edu}}
\date{\today}
\begin{abstract}
We propose and analyze a reliable and efficient a posteriori error estimator for a control--constrained linear--quadratic optimal control problem involving Dirac measures; the control variable corresponds to the amplitude of forces modeled as point sources. The proposed a posteriori error estimator is defined as the sum of two contributions, which are associated with the state and adjoint equations. The estimator associated with the state equation is based on Muckenhoupt weighted Sobolev spaces, while the one associated with the adjoint is in the maximum norm and allows for unbounded right hand sides. The analysis is valid for two and three-dimensional domains. On the basis of the devised a posteriori error estimator, we design a simple adaptive strategy that yields optimal rates of convergence for the numerical examples that we perform.
\end{abstract}
%
%
\maketitle

\section{Introduction}
\label{sec:introduccion}

In this work we shall be interested in the design and analysis of a reliable and efficient a posteriori error estimator for a control--constrained linear--quadratic optimal control problem involving Dirac measures or point sources. To make matters precise, for $n \in \{2,3\}$, we let $\Omega \subset \R^{n}$ be an open and bounded polytopal domain with Lipschitz boundary $\partial \Omega$ and $D$ be a finite ordered subset of $\Omega$ with cardinality $l = \# D$.  Given a desired state $\ysf_d \in L^2(\Omega)$ and a regularization parameter $\lambda > 0$, we define the cost functional 
\begin{equation}
\label{eq:defofJ}
J(\ysf,\mathbf{u}) = \frac12 \| \ysf - \ysf_d \|_{L^2(\Omega)}^2 + \frac\lambda2 \| \mathbf{u} \|_{\R^{l}}^2.
\end{equation}
With these ingredients at hand, we define the \emph{optimal control problem with point sources} as follows: find $\min J(\ysf,\mathbf{u})$ subject to the linear and elliptic state equation 
\begin{equation}
\label{eq:defofPDE}
-\LAP \ysf = \sum_{z \in D} \usf_z \delta_z \text{ in } \Omega, \qquad \ysf = 0 \text{ on } \partial\Omega,
\end{equation}
where $\delta_z$ corresponds to the Dirac delta supported at the point $z \in D$ and 
\begin{equation}
\label{eq:defofu}
\mathbf{u} = \{\usf_z\}_{z \in D} \in \R^{l}, \qquad \asf_z \leq \usf_z \leq \bsf_z  \quad \forall z \in D.  
\end{equation}
Here $\mathbf{u}$ denotes the control variable. The control bounds $\mathbf{a} = \{\asf_z\}_{z \in D}$ and $\mathbf{b} =\{\bsf_z\}_{z \in D}$ both belong to $ \R^{l}$ and satisfy that $\asf_z < \bsf_z$ for all $z \in D$.

Since the state equation \eqref{eq:defofPDE} contains a linear combination of $l$ Dirac measures as a forcing term and $n>1$ the state $\ysf$ does not belong to $H^1(\Omega)$. Consequently, the error analysis involved in the finite element approximation of problem \eqref{eq:defofPDE} is not standard. We refer the reader to \cite{Casas:85,NOS2,Scott:73} for sub-optimal error analyses on quasi--uniform meshes and \cite{ABSV:11,MR3239767} for quasi-optimal results based on graded meshes. 

The mathematical difficulties presented in the study of \eqref{eq:defofPDE} are also present in the analysis of the control problem \eqref{eq:defofJ}--\eqref{eq:defofu}. Based on the function space setting inherited by Muckenhoupt weighted Sobolev spaces \cite{NOS2}, reference \cite{AOS2} provides an a priori error analysis for the control problem \eqref{eq:defofJ}--\eqref{eq:defofu} that relies on the convexity of $\Omega$. The authors propose a fully discrete scheme, on quasi--uniform meshes, that discretizes the state and the corresponding adjoint state using piecewise linear functions and obtain the following error estimates: Let $\epsilon > 0$ and $\ysf_d \in L^q(\Omega)$ for every $q \in (2,\infty)$. If $n=2$, the authors obtain a rate of convergence $\mathcal{O}(h_{\T}^{2-\epsilon})$, in the $\mathbb{R}^{l}$-norm, for the error in the optimal control. If $n=3$, the derived rate is $\mathcal{O}(h_{\T}^{1-\epsilon})$; see \cite[Theorem 5.1]{AOS2}. The fact that these error estimates are not optimal in terms of approximation and the need for both the convexity of $\Omega$ and the higher integrability of the desired state $\ysf_d$ motivate the study of adaptive finite element methods (AFEMs) for the optimal control problem with point sources \eqref{eq:defofJ}--\eqref{eq:defofu}.

AFEMs are iterative methods that improve the quality of the finite element approximation to a partial differential equation (PDE) while striving to keep an optimal distribution of computational resources measured in terms of degrees of freedom. These methods are mainly based on loops of the form
\begin{equation}
\label{eq:loop}
 \textsc{SOLVE} \rightarrow \textsc{ESTIMATE} \rightarrow \textsc{MARK} \rightarrow \textsc{REFINE}. 
\end{equation}
An essential ingredient of an AFEM, which governs the step $\textsc{ESTIMATE}$ in \eqref{eq:loop}, is an a posteriori error estimator. This is a computable quantity that depends on the discrete solution and data and provides information about the local quality of the approximate solution. The a posteriori error analysis for linear second-order elliptic boundary value problems has attained a mature understanding. We refer the reader to \cite{AObook,MR1770058,NSV:09,NV,Verfurth} for an up to date discussion including also the design of AFEMs, their convergence and optimal complexity.

In contrast to the well-established theory for linear elliptic PDEs, the a posteriori error analysis for finite element approximations of a constrained optimal control problem has not yet been fully understood. The main source of difficulty is its inherent nonlinear feature, which appears due to the control constraints. To the best of our knowledge, the first work that provided an advance in this matter is \cite{LiuYan}. These results were later improved in \cite{HHIK} by providing efficiency estimates involving oscillation terms. Recently, these ideas were unified in \cite{KRS}. Unfortunately, the analysis presented in \cite{KRS} relies fundamentally on a particular structure for the underlying problem and the relations among the natural spaces for the state, adjoint state and control; these requirements are not satisfied by the control problem \eqref{eq:defofJ}--\eqref{eq:defofu}. For an up to date survey on a posteriori error analysis for optimal control problems we refer the reader to \cite{AOSR,KRS,MR3485971}.

The main objective of this work is to propose and analyze a reliable and efficient a posteriori error estimator for the optimal control problem with point sources. The proposed error estimator is built on the basis of a suitable error estimator on Muckenhoupt weighted Sobolev spaces that is associated with the state equation and 
a pointwise error estimator that is associated with the adjoint equation. Assuming only that $\Omega$ is a Lipschitz polytope, we prove the global reliability and local efficiency of our proposed error estimator. The analysis is delicate since it involves the interaction of $L^{\infty}(\Omega)$, $\mathbb{R}^l$ and weighted Sobolev spaces, combined with having to deal with the first--order necessary and sufficient optimality condition that characterizes the optimal control $\bar {\mathbf{u}}$. It is important to comment that this work exploits the ideas developed in \cite{AOSR} for the a posteriori error analysis of the so--called \emph{pointwise tracking optimal control problem}. Although the mathematical techniques are similar, the a posteriori error analysis of our control problem does not follow directly from \cite{AOSR}; it requires its own analysis. This is mainly due to the following reasons:
\begin{enumerate}[$\bullet$]
 \item The optimal control variable $\bar{\mathbf{u}}$ belongs to $\mathbb{R}^{l}$, while the one of the problem studied in \cite{AOSR} belongs to $L^2(\Omega)$. This in a sense simplifies the analysis. For instance, as opposed to \cite{AOSR,KRS}, we can obtain local efficiency estimates that do not require convexity of $\Omega$. Nevertheless, it comes with its own set of complications. In particular, the \textit{low regularity} of the state equation.
  
 \item The adjoint problem is a Poisson equation with a forcing term $\ysf - \ysf_d$, \emph{which does not belong to} $L^{\infty}(\Omega)$. Consequently, we must consider a pointwise error indicator that accounts for unbounded right hand sides. Since we were not able to locate one in the literature, in section~\ref{sec:a_posteriori_infty}, we propose such an error indicator and provide its analysis on the basis of \cite{Camacho01072015,demlow2014maximum}. Notice that, thanks to the structure of the control problem of \cite{AOSR}, such an estimator was not needed there.
\end{enumerate}

The outline of this paper is as follows. In section~\ref{sec:notation}, we introduce the notation and functional framework we shall work with. Section \ref{sec:control_problem} contains the description of our control problem and reviews the a priori error analysis developed in \cite{AOS2}. In section \ref{sec:a_posteriori_infty}, we propose and analyze a pointwise a posteriori error estimator for the Laplacian that allows for unbounded right hand sides. Combining this estimator and another one based on Muckenhoupt weighted Sobolev spaces, in section~\ref{sec:a_posteriori} we devise an a posteriori error estimator for our optimal control problem. We show in sections \ref{subsub:reliable} and \ref{subsub:efficient}, its reliability and efficiency, respectively. We conclude, in section~\ref{sec:numex}, with a series of numerical examples that illustrate and go beyond our theory.
\section{Notation and preliminaries}
\label{sec:notation}

Let us fix notation and the setting in which we will operate. Throughout this work, $n \in \{2,3\}$ and $\Omega$ is an open and bounded polytopal domain of $\R^n$ with Lipschitz boundary $\partial \Omega$. If $\Xcal$ and $\Ycal$ are normed vector spaces, we write $\Xcal \hookrightarrow \Ycal$ to denote that $\Xcal$ is continuously embedded in $\Ycal$. We denote by $\Xcal'$ and $\|\cdot\|_{\Xcal}$ the dual and the norm of $\Xcal$, respectively.

For $E \subset \Omega$ of finite Hausdorff $i$-dimension, $i \in \{1,2,3\}$, we denote its measure by $|E|$. The mean value of a function $f$ over a set $E$ is
\[
 \fint_E f  = \frac{1}{|E|}\int_{E} f .
\]

The relation $a \lesssim b$ indicates that $a \leq C b$, with a constant $C$ which is independent of $a$, $b$ and the size of the elements in the mesh. The value of $C$ might change at each occurrence.

\subsection{Weighted Sobolev spaces}
\label{sub:weighted_Sobolev}

We start this section with a notion which will be fundamental for further discussions, that of a weight. A weight is a locally integrable, nonnegative function defined on $\mathbb{R}^n$. If $\omega$ is a weight, we say that $\omega$ belongs to the so--called Muckenhoupt class $A_2$, or that it is an $A_2$-weight,  if there is a constant $C_{\omega}$ such that
\begin{equation}
\label{A_pclass}
C_{\omega} = \sup_{B} \left( \frac{1}{|B|} \int_{B} \omega \right) \left( \frac{1}{|B|} \int_{B} \omega^{-1} \right)  < \infty,
\end{equation}
where the supremum is taken over all balls $B$ in $\R^n$ \cite{Javier,FKS:82,Muckenhoupt,Turesson}. 

We present an important example of a Muckenhoupt weight. Let $x_0$ be an interior point of $\Omega$ and denote by $\dist(x)$ the Euclidean distance to $x_0$. Define $\dist^\alpha(x) = \dist(x)^{\alpha}$. We then have that $\dist^\alpha \in A_2$ if and only if $\alpha \in (-n,n)$. We refer the reader to \cite{Javier,NOS2,Turesson} for more examples of $A_2$-weights and their most important properties. Since it will be necessary for our analysis, here we mention the following \emph{reverse H\"older inequality} for $A_2$-weights. Its proof can be found in \cite[Theorem 7.4]{Javier} or \cite[Lemma 1.2.12]{Turesson}. 

\begin{prpstn}[reverse H\"older inequality]
\label{prop:revHolder}
Let $\omega \in A_2$, then there is $\epsilon > 0$ such that for every ball $B \subset \R^n$ we have
\[
  \left( \fint_B \omega^{1+\epsilon} \right)^{1/(1+\epsilon)} \lesssim \fint_B \omega,
\]
where the hidden constant depends only on the dimension $n$.
\end{prpstn}

To analyze problem \eqref{eq:defofPDE}, we consider Lebesgue and Sobolev spaces with Muckenhoupt weights. If $\Omega$ is an open and bounded domain of $\mathbb{R}^n$ and $\omega \in A_2$, we define 
\begin{equation}
\label{eq:L2wspace}
 L^2(\omega, \Omega) := \left \{ v \in L^1_{\textrm{loc}}(\Omega): \| v \|_{L^2(\omega, \Omega)} := \left( \int_{\Omega} |v|^2 \omega \right)^{\frac{1}{2}} < \infty \right \},
\end{equation}
and 
\begin{equation}
\label{eq:H1wspace}
 H^1(\omega, \Omega) := \left \{ v \in L^2(\omega,\Omega): | \nabla v |\in L^2(\omega, \Omega) \right \},
\end{equation}
with norm 
\begin{equation}
\label{eq:norm}
\| v \|_{H^1(\omega, \Omega)}:= \left( \| v \|^2_{L^2(\omega, \Omega)} + \| \nabla v \|^2_{L^2(\omega, \Omega)} \right)^{\tfrac{1}{2}}.
\end{equation}
We also define $H_0^1(\omega, \Omega)$ as the closure of $C_0^{\infty}(\Omega)$ in $H^1(\omega, \Omega)$. In view of the fact that $\omega$ is an $A_2$-weight, \cite[Theorem 1.3]{FKS:82} guarantees that a weighted Poincar\'e inequality holds and thus that $ \| \nabla v \|_{L^2(\omega, \Omega)} $ is an equivalent norm to $\| v \|_{H^1(\omega, \Omega)}$.

The notion of a Muckenhoupt weight has important consequences and we conclude this section by mentioning some of them. If $\omega \in A_2$, then we have that $H^1(\omega,\Omega)$ is Hilbert and $H^1(\omega,\Omega) \cap C^{\infty}(\Omega)$ is dense in $H^1(\omega,\Omega)$ (cf.~\cite[Proposition 2.1.2, Corollary 2.1.6]{Turesson} and \cite[Theorem~1]{GU}).

\subsection{The Poisson problem in Lipschitz polytopes}
\label{sub:reg_of_y}

In this section we collect some standard results concerning the regularity of the solution to the Poisson problem
\begin{equation}
\label{eq:Fish}
  - \Delta u = f \ \text{in } \Omega, \quad u=0 \ \text{on } \partial\Omega,
\end{equation}
where $\Omega$ is a bounded and Lipschitz, but not necessarily convex, polytope. We begin with the following result \cite{Dauge:92,Grisvard,JK:95,JK:81,MR2641539,Savare:98}.

\begin{prpstn}[higher integrability]
\label{prop:uisW1p}
Let $u \in H^1_0(\Omega)$ denote the unique solution of \eqref{eq:Fish} with $f \in L^2(\Omega)$. There is $q>n$ such that $u \in W^{1,q}(\Omega)$. Moreover,
\[
  \| u \|_{W^{1,q}(\Omega)} \lesssim \| f \|_{L^2(\Omega)},
\]
where the hidden constant is independent of $u$ and $f$. This, in particular, implies that for $\kappa = 1-n/q>0$ we have $u \in C^{0,\kappa}(\bar\Omega)$ with a similar estimate.
\end{prpstn}

We now present a local regularity result, whose proof can be found, for instance, in  \cite[Theorem 9.11]{GT} or \cite[Theorem 12.2.2]{MR3012036}.

\begin{prpstn}[local regularity]
\label{prop:loclreg}
Let $u \in H^1_0(\Omega)$ denote the unique solution of \eqref{eq:Fish} with $f \in L^r(\Omega)$ and $r \in [2,\infty)$. If $U \Subset \Omega$ then $u \in W^{2,r}(U)$ and the following estimate holds
\[
  \| u \|_{W^{2,r}(U)} \lesssim \| u \|_{L^r(\Omega)} + \| f \|_{L^r(\Omega)},
\]
where the hidden constant depends on $\textup{dist}(U,\partial\Omega)$ but is independent of $u$ and $f$.
\end{prpstn}

We remark that, since $\Omega$ is bounded, the estimates of Propositions~\ref{prop:uisW1p} and \ref{prop:loclreg} allow us to obtain that, for every $U \Subset \Omega$,
\begin{equation}
\label{eq:W2p}
  \| u \|_{W^{2,r}(U)} \lesssim  \| f \|_{L^r(\Omega)},  
\end{equation}
where the hidden constant depends on $|\Omega|$ and $\textup{dist}(U,\partial\Omega)$ but is independent of $u$ and $f$.


\section{The optimal control problem with point sources}
\label{sec:control_problem}

In this section we precisely describe and analyze the \emph{optimal control problem with point sources} introduced in section \ref{sec:introduccion}.  We begin by assuming that we are given an ordered set $D \subset \Omega$ with finite cardinality $l$. We define 
\[
d_{D} = \left\{\begin{array}{ll}
\mathrm{dist}(D,\partial \Omega), & \mbox{if }l=1,
\\
\min \left \{ \mathrm{dist}(D,\partial \Omega), \min\{|z-z'|: z,z' \in D, \ z\neq z' \} \right \}, & \mbox{otherwise}.
\end{array}\right.
\]
Since $D$ is finite and $D \subset \Omega$ we have that $d_{D}>0$.
It is then suitable, for our analysis, to define the weight $\rho$ as follows: If $l = 1$, then
\begin{equation}
\label{eq:defofoweight1}
\rho(x) = \distz^{\alpha}(x), 
\end{equation}
otherwise
\begin{equation}
\label{eq:defofoweight}
  \rho(x) = \begin{dcases}
                \distz^{\alpha}(x), & \exists z \in D: \distz(x) < d_{D}/2, \\
                1, & \distz(x) \geq d_{D}/2 ~~ \forall z \in D.
              \end{dcases}
\end{equation}
Here $\alpha \in (n-2,n)$ and $\distz(x) = |x-z|$ denotes the Euclidean distance to $z$. Since $\alpha \in (n-2,n) \subset (-n,n)$, the weight $\rho$ belongs to the Muckenhoupt class $A_2$ \cite{MR3215609}. Consequently, $H^1(\rho,\Omega)$, defined by \eqref{eq:H1wspace}, is a Hilbert space endowed with the norm \eqref{eq:norm}. We state the following embedding result \cite[Lemma 2]{AOSR}.
\begin{lmm}[$H^1_0(\rho,\Omega) \hookrightarrow L^2(\Omega)$]
\label{lem:restalphagen}
If $\alpha \in (n-2,2)$ then $H^1_0(\rho,\Omega) \hookrightarrow L^2(\Omega)$ and we have the following weighted Poincar\'e inequality
\[
  \| v \|_{L^2(\Omega)} \lesssim \| \nabla v \|_{L^2(\rho,\Omega)} \quad \forall v \in H^1_0(\rho,\Omega),
\]
where the hidden constant depends only on $\Omega$ and $d_D$.
\end{lmm}

We define the set of admissible controls by
\begin{equation}
\label{eq:Uz}
\mathcal{U}_{\textrm{ad}} = \left\{ \mathbf{u} \in \R^l: \asf_z \leq 	\usf_z \le \bsf_z  \quad \forall z  \in D \right\},
\end{equation}
where the control bounds $\mathbf{a}$ and $\mathbf{b}$ both belong to $\R^{l}$ and satisfy that $\asf_z < \bsf_z$ for all $z \in D$. The set $\mathcal{U}_{\textrm{ad}}$ is a nonempty, closed, and convex subset of $\R^l$. 

We recall that the cost functional $J$ is defined by \eqref{eq:defofJ} and thus define the optimal control problem with point sources as follows: Find $\min J (\ysf,\mathbf{u})$ subject to the following weak formulation of the state equation \eqref{eq:defofPDE}:
\begin{equation}
\label{eq:state_equation}
\ysf \in H^1_0(\rho, \Omega): \quad   (\nabla \ysf, \nabla \vsf)_{L^2(\Omega)} = \sum_{z \in D } \usf_z \langle \delta_z, \vsf \rangle \quad \forall \vsf \in H^1_0(\rho^{-1},\Omega),
\end{equation}
and the control constraints $\mathbf{u} \in \mathcal{U}_{\textrm{ad}}$. Here $\langle\cdot,\cdot\rangle$ denotes the duality pairing between $H^1_0(\rho^{-1},\Omega)$ and its dual $H^1_0(\rho^{-1},\Omega)'$; the results of \cite[Theorem 2.3]{AGM}, \cite{DAngelo:SINUM} and \cite[Lemma 7.1.3]{KMR} guarantee that $\delta_z \in H^1_0(\rho^{-1},\Omega)'$ for $\alpha \in (n-2,n)$ and that \eqref{eq:state_equation} is well--posed. On the other hand, the continuous embedding of Lemma \ref{lem:restalphagen}, \ie~ $H_0^1(\rho,\Omega) \hookrightarrow L^2(\Omega)$ for $\alpha \in (n-2,2)$, and the fact that $\ysf_d \in L^2(\Omega)$ imply that $\ysf - \ysf_d \in L^2(\Omega)$. We have thus concluded that $J$ is well-defined on $H_0^1(\rho,\Omega) \times \mathbb{R}^l$ for $\alpha \in (n-2,2)$.

To analyze the optimal control problem with point sources we introduce the control--to--state map $S: \mathbb{R}^l \rightarrow H_0^1(\rho,\Omega)$ as follows:  given a control $\mathbf{u}$, the map $S$ associates to it a unique state $\ysf$ that solves problem \eqref{eq:state_equation}. Since $\alpha \in (n-2,2)$, $S$ is well-defined, and with this operator at hand, we define the reduced cost functional
\begin{equation}
\label{eq:f}
f(\mathbf{u}):=\frac{1}{2} \|  S\mathbf{u} - \ysf_{d}  \|^2_{L^2(\Omega)}+ \frac{\lambda}{2} \| \mathbf{u} \|^2_{\mathbb{R}^l}.
\end{equation}
We immediately conclude that $f$ is weakly lower semicontinuous and strictly convex ($\lambda > 0$). This, combined with the fact that $\mathcal{U}_{\textrm{ad}}$ is compact, allows us to obtain the existence and uniqueness of an optimal control $\bar{ \mathbf{u}} \in \mathcal{U}_{\textrm{ad}} $ and an optimal state $\bar \ysf   = S\bar{\mathbf{u}} \in H_0^1(\rho,\Omega)$ that satisfy \eqref{eq:state_equation}; see \cite[Theorem 1.2]{LionsControl} and \cite[Theorem 2.14]{Tbook}. In addition, we have that the optimal control $\bar{ \mathbf{u}}$ satisfies \cite[Lemma 2.1]{Tbook}:
\begin{equation}
\label{eq:first_order}
f'(\bar{ \mathbf{u}})  ( \mathbf{u}-\bar{ \mathbf{u}} ) \geq 0 \qquad \forall  \mathbf{u} \in \mathcal{U}_{\textrm{ad}}.
\end{equation}
This variational inequality is necessary and sufficient for optimality. To explore it, we define the adjoint variable $\psf$ as the unique solution to
\begin{equation}
\label{eq:adjp}
\psf \in H^1_0(\Omega): \quad  (\nabla \wsf,\nabla \psf)_{L^2(\Omega)} = ( \ysf - \ysf_d, \wsf)_{L^2(\Omega)} \quad \forall \wsf \in H^1_0(\Omega).
\end{equation}
Since $\ysf - \ysf_d \in L^2(\Omega)$, the well-posedness of \eqref{eq:adjp} is immediate.  In 
addition, the results of Proposition \ref{prop:uisW1p} guarantee that  $\psf$ is H\"older continuous. We now derive a local regularity result for the adjoint state $\psf$. To do this, we define
\begin{equation}
\label{eq:r_star}
  r^* = \begin{dcases}
         \infty & n =2, \\
          3 & n = 3.
        \end{dcases}
\end{equation} 
Standard arguments reveal the regularity properties of the solution to problem \eqref{eq:state_equation}: $\ysf \in L^{r}(\Omega)$ for every $r < r^{*}$.

\begin{prpstn}[local regularity]
\label{prop:loclreg2}
Let $\psf \in H^1_0(\Omega)$ denote the unique solution of \eqref{eq:adjp}. If $\alpha \in (n-2,2)$, $\ysf_d \in L^r(\Omega)$, for every $r <r^*$, and $U \Subset \Omega$, then $\psf \in W^{2,r}(U)$ for every $r$ such that $2<r<r^*$.
\end{prpstn}
\begin{proof}
Standard regularity results for $\ysf$ and the assumption that $\ysf_d$ satisfies allow us to immediately conclude that $\ysf - \ysf_d \in L^r(\Omega)$ for every $r$ such that $2<r<r^*$. We thus conclude by applying the results of Proposition \ref{prop:loclreg}.
\end{proof}

With this result at hand, we obtain a weighted integrability result for $\psf$.

\begin{prpstn}[weighted integrability]
\label{pro:weighted_integra}
Let $\psf \in H^1_0(\Omega)$ denote the unique solution of \eqref{eq:adjp}. If $\alpha \in (n-2,2)$ and $\ysf_d \in L^r(\Omega)$ for every $r < r^*$, then $\psf \in H_0^1(\rho^{-1},\Omega)$. 
\end{prpstn}
\begin{proof} 
For each $z \in D$, let $B(z)$ denote the ball with center $z$ and radius $d_D/2$. Set $G = \Omega \setminus \cup_{z \in D} B(z)$ and compute
\[
  \int_\Omega \rho^{-1} | \nabla \psf |^2 = \sum_{z \in D} \int_{B(z)} \rho^{-1} | \nabla \psf|^2 + \int_G \rho^{-1} | \nabla \psf |^2.
\]
By definition of $G$, there is a constant $g>0$, that depends only on $d_D$ and $\alpha$, such that $\rho(x) \geq g$ for every $x \in G$, thus
\[
  \int_G \rho^{-1} | \nabla \psf |^2 \leq g^{-1} \int_G |\nabla \psf |^2 \lesssim \| \ysf - \ysf_d \|_{L^2(\Omega)}^2.
\]
We now bound the integral near the support of the Dirac measures. First, since $\#D$ is finite, it suffices to consider a single ball. Next, owing to $B(z) \Subset \Omega$ we invoke Proposition~\ref{prop:loclreg2} and conclude that $\psf \in W^{2,r}(B(z))$ for every $r$ such that $2<r<r^*$. If $n=2$, in view of the fact that $W^{1,r}(B(z) ) \hookrightarrow L^{\infty}(B(z))$ for $r>2$, we have that
\[
  \| \nabla \psf \|_{L^{\infty}(B(z))} 	\lesssim   \| \psf \|_{W^{2,r}(B(z))} \lesssim  \| \ysf - \ysf_d \|_{L^{r}(\Omega)},
\]
and thus that
\[
  \int_{B(z)} \rho^{-1} | \nabla \psf |^2 \lesssim \| \ysf - \ysf_d \|_{L^{r}(\Omega)}^2.
\]

In three dimensions ($n=3$) we do not have that $\psf$ is Lipschitz and, thus, we must employ a different argument. Namely,  if $\epsilon>0$ we have, by H\"older's inequality,
\begin{equation}
\label{eq:savesin3d}
  \int_{B(z)} \rho^{-1} | \nabla \psf |^2 \leq \left(\int_{B(z)} \rho^{-(1+\epsilon)} \right)^{1/(1+\epsilon)} \left( \int_{B(z)} |\nabla \psf|^{2(1+\epsilon)/\epsilon} \right)^{\epsilon/(1+\epsilon)}.
\end{equation}
Now, invoking the reverse H\"older inequality of Proposition~\ref{prop:revHolder} we have
\[
  \left(\int_{B(z)} \rho^{-(1+\epsilon)} \right)^{1/(1+\epsilon)} \lesssim |B(z)|^{1/(1+\epsilon)} \fint_{B(z)} \rho^{-1},
\]
which, since $\#D$ is finite, is uniformly bounded. Finally, we recall that, by Proposition~\ref{prop:loclreg2}, $\nabla \psf \in W^{1,r}(B(z))$ for every $r<3$ which, in turn, implies that $\nabla \psf \in L^q(B(z))$ for all $q \leq 3r/(3-r)$. Choosing, in \eqref{eq:savesin3d}, the value of $\epsilon$ given by Proposition~\ref{prop:revHolder} gives a uniform bound on the first factor. In addition, once $\epsilon$ is fixed, so is $2(1+\epsilon)/\epsilon$. Therefore, since
\[
  \lim_{r\uparrow 3} \frac{3r}{3-r} = \infty \Longrightarrow \exists r_0 \in (1,3) : \frac{3r_0}{3-r_0} \geq 2(1+\epsilon)/\epsilon,
\]
which allows us to conclude that the second factor on the right hand side of \eqref{eq:savesin3d} is also bounded.

This shows that $\psf \in H^1_0(\rho^{-1},\Omega)$.
\end{proof}

With these ingredients at hand, we proceed to show optimality conditions for our problem.
\begin{thrm}[optimality conditions]
\label{thm:optim_conds}
Let $\alpha \in (n-2,2)$. The pair $(\bar{\ysf},\bar{\mathbf{u}}) \in H_0^1(\rho,\Omega) \times\mathbb{R}^l$ is optimal for the optimal control problem with point sources if and only if $\bar{\mathbf{u}} \in \mathcal{U}_{\textrm{ad}}$, $\bar{\ysf} = S \bar{\mathbf{u}}$, and the optimal control $\bar{\mathbf{u}}$ satisfies
\begin{equation}
\label{eq:VIp}
\sum_{z \in D}  ( \bar \psf(z) + \lambda \bar{ \usf }_{z} )  (  \usf_{z}- \bar{ \usf }_{z})  \geq  0 \qquad \forall \mathbf{u}\in \mathcal{U}_{\textrm{ad}},
\end{equation}
where the optimal adjoint state $\bar \psf \in H^1_0(\Omega)$ solves \eqref{eq:adjp} with $\ysf=\bar\ysf$.
\end{thrm}
\begin{proof}
A basic computation reveals that the the first--order optimality condition \eqref{eq:first_order} reads, for every $ \mathbf{u} \in \R^l$, as follows:
\[
0 \leq f'( \bar{ \mathbf{u}} )(\mathbf{u} - \bar {\mathbf{u}} ) = (S \bar{ \mathbf{u}}  - \ysf_d, S(\mathbf{u}-\bar{ \mathbf{u}}) )_{L^2(\Omega)} + \lambda (\bar{ \mathbf{u}}  , \mathbf{u}-\bar{ \mathbf{u}} )_{\R^l}. 
\]
Since the second term on the right hand side of the previous expression is already present in the desired variational inequality \eqref{eq:VIp}, we investigate the first term. To accomplish this task, we first set $\ysf = S  \mathbf{u}$ and $\bar \ysf = S  \bar {\mathbf{u}}$ and notice that $\ysf - \bar \ysf$ solves
\begin{equation}
\label{eq:state_equation_oc}
 (\nabla( \ysf - \bar \ysf), \nabla \vsf)_{L^2(\Omega)} = \sum_{z \in D } ( \usf_z - \bar \usf_z) \langle \delta_z, \vsf \rangle \quad \forall \vsf \in H^1_0(\rho^{-1},\Omega).
\end{equation}
In view of the results of Proposition \ref{pro:weighted_integra}, we are allowed to set $\vsf = \bar \psf$ in \eqref{eq:state_equation_oc}. This yields
\begin{equation}
\label{eq:aux_op_1}
  (\nabla( \ysf - \bar \ysf), \nabla \bar \psf)_{L^2(\Omega)} = \sum_{z \in D } ( \usf_z - \bar \usf_z) \bar \psf(z);
\end{equation}
notice that we have also used the H\"older continuity of $\psf$, given in Proposition \ref{prop:uisW1p}, to guarantee that $\langle \delta_z, \bar \psf  \rangle = \bar \psf(z)$.

Now, we would like to set $\wsf = \ysf - \bar \ysf$ in \eqref{eq:adjp} to conclude that
\begin{equation}
\label{eq:aux_op_2}
 (\nabla(\ysf - \bar \ysf),\nabla \bar \psf)_{L^2(\Omega)} = (\bar \ysf - \ysf_d, \ysf - \bar \ysf)_{L^2(\Omega)}.
\end{equation}
This, on the basis of \eqref{eq:aux_op_1}, would allows us to obtain \eqref{eq:VIp}. Unfortunately $\ysf-\bar\ysf \notin H_0^1(\Omega)$ and thus we need to justify \eqref{eq:aux_op_2} with a different argument. Let $w_n \in C_0^{\infty}(\Omega)$ be such that $w_n \rightarrow \ysf - \bar \ysf$ in $H_0^1(\rho,\Omega)$. Setting $\wsf = w_n$ in \eqref{eq:adjp} yields
\[
  (\nabla w_n,\nabla \bar \psf)_{L^2(\Omega)} = ( \bar \ysf - \ysf_d,w_n)_{L^2(\Omega)}.
\]
The continuity of the variational form for the Dirichlet Laplace operator on $H_0^1(\rho,\Omega) \times H_0^1(\rho^{-1},\Omega)$ \cite[Proposition 2.1]{AGM} and the regularity results of Proposition \ref{pro:weighted_integra} imply that $(\nabla w_n,\nabla \bar \psf)_{L^2(\Omega)}$ converges to $(\nabla (\ysf -\bar \ysf ),\nabla \bar \psf)_{L^2(\Omega)}$ as $n \rightarrow \infty$. Finally, we invoke Lemma \ref{lem:restalphagen} to obtain the convergence of $( \bar \ysf - \ysf_d,w_n)_{L^2(\Omega)}$ to $( \bar \ysf - \ysf_d,\ysf - \bar \ysf)_{L^2(\Omega)}$ as $n \rightarrow \infty$ and conclude the proof.
\end{proof}

To summarize, the pair $(\bar{\ysf},\bar{ \mathbf{u}}) \in H^1_0(\rho, \Omega)\times \mathcal{U}_{\textrm{ad}}$ 
is optimal for the optimal control problem with point sources \eqref{eq:defofJ}--\eqref{eq:defofu} if and only if the triple $(\bar{\ysf},\bar{ \mathbf{u}},\bar{\psf}) \in H^1_0(\rho, \Omega)\times \mathcal{U}_{\textrm{ad}} \times H_0^1(\Omega)$ satisfies the following optimality system:
\begin{equation}\label{eq:optimal_system}
\left\{
\begin{array}{cl}
(\nabla \bar{\ysf}, \nabla \vsf)_{L^2(\Omega)}  =  \displaystyle \sum_{z \in D } \bar{\usf}_z \langle \delta_z, \vsf \rangle 
& \forall \vsf \in H^1_0(\rho^{-1},\Omega),\\
(\nabla \wsf, \nabla \bar{\psf})_{L^2(\Omega)}  =  (\bar{\ysf}-\ysf_{d},\wsf)_{L^2(\Omega)}
& \forall \wsf \in H^1_0(\Omega),
\vspace{0.3cm} \\
\displaystyle \sum_{z \in D}  ( \bar \psf(z) + \lambda \bar{ \usf }_{z} )  (  \usf_{z}- \bar{ \usf }_{z})  \geq  0
&
\forall \mathbf{u}\in \mathcal{U}_{\textrm{ad}}.
\end{array}
\right.
\end{equation}

In the spirit of \cite[section 2.8]{Tbook} and \cite[chapter 2]{LionsControl}, we present the following projection formula for $\bar{\mathbf{u}}$ that is equivalent to \eqref{eq:VIp}. This formula, for $z \in D$, reads: 
\begin{equation}
\label{eq:Pi}
\bar{\usf}_{z} = \max \left\{ \asf_{z}, \min\left\{ \bsf_{z}, -\frac1\lambda \bar{\psf} (z)\right\} \right\}.
\end{equation}

\subsection{Finite element discretization}
\label{sub:fem}

We recall the finite element approximation of the control problem with point sources detailed in \cite{AOS2}. In doing so, we consider $\T = \{T\}$ to be a conforming partition of $\Omega$ into simplices $T$ with size $h_T = \diam(T)$ and define $h_{\T} = \max_{T \in \T} h_T$.  We denote by $\Tr$ the collection of conforming and shape regular meshes that are refinements of an initial mesh $\T_0$. 

Given a mesh $\T \in \Tr$, we define the finite element space of continuous piecewise polynomials of degree one as
\begin{equation}
\V(\T) = \left\{v_{\T} \in C^0( \bar \Omega): {v_{\T}}_{|T} \in \mathbb{P}_1(T) \ \forall T \in \T, \ v_{\T} = 0 \textrm{ on } \partial \Omega \right\}.
\label{eq:defFESpace}
\end{equation}

With these ingredients at hand, we present a finite element discretization for our optimal control problem. The optimal state and adjoint state are discretized on the basis of $\V(\T)$. We remark that no discretization is needed for the optimal control variable: the admissible set $\mathcal{U}_{\textrm{ad}}$ is a subset of a finite dimensional space. Then, the discrete counterpart of \eqref{eq:defofJ}--\eqref{eq:defofu} reads: Find
$
\min J (\ysf_{\T},\mathbf{u}_\T)
$
subject to the discrete state equation
\begin{equation}
\label{eq:discrete_state_equation}
\ysf_{\T} \in \V(\T): \quad (\nabla \ysf_{\T}, \nabla \vsf_{\T} )_{L^2(\Omega)} = \sum_{z \in D } \usf_{\T,z} \langle \delta_z, \vsf_{\T} \rangle \quad \forall \vsf_{\T} \in \V(\T),
\end{equation}
and the control constraints
$
\mathbf{u}_\T = \{\usf_{\T,z}\}_{z \in D} \in \mathcal{U}_{\textrm{ad}}.
$
Similar arguments to those used in section \ref{sec:control_problem} allow us to conclude that the pair $(\bar{\ysf}_{\T},\bar{\mathbf{u}}_\T)$ is optimal for the discrete optimal control problem with point sources if and only if $\bar{\ysf}_{\T}$ solves \eqref{eq:discrete_state_equation} and $\bar{\mathbf{u}}_{\T}$ is such that
\begin{equation}
\label{eq:first_oder_discrete}
  \sum_{z \in D}( \bar{\psf}_{\T}(z) +\lambda \bar{\usf}_{\T,z}) (\usf_{z} - \bar{\usf}_{\T,z})  \geq 0 \quad \forall \mathbf{u} \in \mathcal{U}_{\textrm{ad}},
\end{equation}
where $\bar{{\psf}}_{\T}$ denotes the unique solution to 
\begin{equation}
\label{eq:discrete_adjoint_equation}
\bar{\psf}_{\T} \in \V(\T): \quad (\nabla \vsf_{\T}, \nabla \bar{\psf}_{\T} )_{L^2(\Omega)} = (\bar{\ysf}_{\T} - \ysf_d, \vsf_{\T})_{L^2(\Omega)} \quad \forall \vsf_{\T} \in \V(\T).
\end{equation}

The following a priori error analysis follows from \cite{AOS2}: Let $\epsilon > 0$ and $\Omega_1$ be such that $D \Subset \Omega_1 \Subset \Omega$. Assume that for every $q \in (2,\infty)$, $\ysf_d \in L^q(\Omega)$, $\Omega$ is convex, and the mesh $\T$ is quasiuniform with mesh size $h_{\T}$. If
$n=2$, then we have
\begin{equation}
\label{eq:estimate_1} 
\| \bar{\mathbf{u}} - \mathbf{u}_\T \|_{\R^l} \lesssim h_\T^{2-\epsilon} \left( \| \bu \|_{\R^l} + \| \psf \|_{H^2(\Omega)} + \| \psf \|_{W^{2,r}(\Omega_1)} \right).
\end{equation}
On the other hand, if $n=3$, then
\begin{equation}
\label{eq:estimate_2}
  \| \bar{\mathbf{u}} - \mathbf{u}_\T \|_{\R^l} \lesssim h_\T^{1-\epsilon} \left( \| \bu \|_{\R^l} + \| \psf \|_{H^2(\Omega)} +\| \psf \|_{W^{2,r}(\Omega_1)} \right),
\end{equation}
where $r < n/(n-2)$. The hidden constants in both estimates are independent of the size of the elements in the mesh $\T$, $\# \T$, and the continuous and discrete optimal pairs. We also refer the reader to \cite{MR3225501} for another a priori error analysis.

\section{Pointwise a posteriori error estimation}
\label{sec:a_posteriori_infty}
The a posteriori error estimator proposed in the next section is built on the basis of two error contributions: one associated to the state equation \eqref{eq:state_equation} and another one that involves the adjoint problem \eqref{eq:adjp}. Since the variational inequality \eqref{eq:VIp}, that characterizes the optimal control, involves point evaluations of the optimal adjoint state, it is thus imperative to consider a pointwise error estimator for the adjoint problem \eqref{eq:adjp}. 

The development and analysis of pointwise a posteriori error estimators have been considered in a number of articles. Starting with the pioneering works, in two dimensions, by Eriksson \cite{Eriksson} and Nochetto \cite{rhn}, the theory has been extended to more dimensions and both nonlinear and geometric problems \cite{Camacho01072015,deldia,MR3022214,demlow2014maximum,MR2249676,Nochetto2002}. In most of these works it is assumed that the right hand side of the underlying PDE is bounded. However, in our setting, this does not hold because the adjoint problem \eqref{eq:adjp} has the function $\ysf - \ysf_d \notin L^{\infty}(\Omega)$ as forcing term. In fact, the function $\ysf$, that solves \eqref{eq:state_equation}, belongs to $L^{r}(\Omega)$ for every $r < r^*$, where $r^*$ is defined in \eqref{eq:r_star}.
For this reason, here we develop an a posteriori error analysis for the Laplacian in the maximum norm and with an unbounded forcing term.

We must immediately remark that the results presented in this section are not new \emph{per se}. They are essentially contained in \cite{Camacho01072015,demlow2014maximum} and we develop them not just for the sake of completeness, but also because we could not find them in the form that is necessary for our purposes.

To make matters precise, let $\Omega$ be an open and bounded polytopal domain of $\R^n$ with Lipschitz boundary $\partial \Omega$ and $f \in L^2(\Omega)$.
Let $u$ be the weak solution to:
\begin{equation}
\label{eq:lap_infty}
 u \in H_0^1(\Omega): \quad (\nabla u, \nabla v )_{L^2(\Omega)} = (f,v)_{L^2(\Omega)} \quad \forall v \in H_0^1(\Omega).
\end{equation}
Under these assumptions, Proposition \ref{prop:uisW1p} guarantees the existence of $q>n$ such that $u \in W^{1,q}(\Omega)$. This, in view of the embedding $W^{1,q}(\Omega) \hookrightarrow C(\bar \Omega)$, implies that $u \in C(\bar \Omega)$; it is then appropriate to study a posteriori error estimation in $L^{\infty}(\Omega)$. To accomplish this task, we first need to introduce and set some notation in addition to that of section~\ref{sub:fem}. We define the Galerkin approximation to problem \eqref{eq:lap_infty} by
\begin{equation}
\label{eq:lap_infty_dis}
 u_{\T} \in \V(\T): \quad (\nabla u_{\T}, \nabla v_{\T} )_{L^2(\Omega)} = (f,v_{\T})_{L^2(\Omega)} \quad \forall v_{\T} \in \V(\T).
\end{equation}

We define $\Sides$ as the set of internal ($n-1$)-dimensional interelement boundaries $S$ of $\T$. For $S \in \Sides$, we indicate by $h_S$ the diameter of $S$. For $T\in\T$, let $\Sides_{T}$ denote the subset of $\Sides$ which contains the sides in $\Sides$ which are sides of $T$. We also denote by $\Ne_S$ the subset of $\T$ that contains the two elements that have $S$ as a side. In addition, we define the \emph{stars} or \emph{patches} associated with an element $T \in \T$ as
\begin{equation}
  \Ne_T := \bigcup_{T' \in \T : T \cap T' \neq \emptyset} T'
  \label{eq:NeT}
\end{equation}
and
\begin{equation}
  \Ne_T^* := \bigcup_{T' \in \T : \Sides_T \cap \Sides_{T'} \neq \emptyset} T'.
  \label{eq:NeTstar}
\end{equation}

Given a discrete function $v_{\T} \in \V(\T)$, with $\V(\T)$ defined in \eqref{eq:defFESpace}, we define, for any internal side $S \in \Sides$, the jump or interelement residual $\llbracket \nabla v_{\T} \cdot \nu  \rrbracket$ by
\begin{equation}
\label{eq:jump}
\llbracket \nabla v_{\T} \cdot \nu  \rrbracket = \nu^+ \cdot \nabla v_{\T|T^+} + \nu^- \cdot \nabla v_{\T|T^-} ,
\end{equation}
where $\Ne_S = \{ T^+, T^-\}$ and $\nu^+, \nu^-$ denote the unit normals to $S$ pointing towards $T^+$, $T^{-} \in \T$, respectively, which are such that $T^+ \neq T^-$ and $\partial T^+ \cap \partial T^-=S$.

With these ingredients at hand, we introduce the local a posteriori error indicator
\begin{equation}
 \label{eq:defofEinf}
 \E (u_{\T}; T)= h_T^{2-n/2} \| f \|_{L^{2}(T)} 
 + h_T \| \llbracket \nabla u_{\T}  \cdot \nu\rrbracket \|_{L^{\infty}(\partial T \setminus \partial \Omega)}.
\end{equation}
The global pointwise estimator for problem \eqref{eq:lap_infty} is then defined by
\begin{equation}
\label{eq:defofEinfglobal}
\E_{\infty}(u_{\T}; \T) = \max_{T \in \T } \E (u_{\T}; T).
\end{equation}

We notice that the local indicator \eqref{eq:defofEinf} contains the term $h_T^{2-n/2} \| f \|_{L^{2}(T)}$ instead of the standard consideration in the literature: $h_T^{2} \| f \|_{L^{\infty}(T)}$ \cite{deldia,demlow2014maximum,rhn}. This allows for a pointwise a posteriori error analysis with unbounded right hand sides \cite{Camacho01072015}. In the remainder of this section we will investigate the global reliability and local efficiency of the estimator \eqref{eq:defofEinf}--\eqref{eq:defofEinfglobal}.

\subsection{Reliability}
\label{sec:a_posteriori_infty_rel}
A standard technique for performing an error analysis for finite element approximations in the maximum norm is to represent the pointwise error with the help of a Green's function. For each $x \in \Omega$, the Green's function $\calG(x,y): \Omega \times \Omega \rightarrow \R$ is defined as the solution (in the sense of distributions) to 
\begin{equation}
\label{eq:green}
-\Delta_y \calG = \delta(x-y) \quad y \in \Omega, \qquad  \calG(x,y) = 0 \quad y \in \partial \Omega.
\end{equation}
This definition implies that, if $w \in H_0^1(\Omega) \cap W^{1,q}(\Omega)$, for $q > n$, then the following pointwise representation holds:
\begin{equation}
\label{eq:point_representation}
 w(x) = (\nabla w, \nabla \calG(x, \cdot) )_{L^2(\Omega)}.
\end{equation}
Let us summarize some properties of the Green's function that will be useful.

\begin{prpstn}[properties of $\calG$]
\label{prop:propertiesG}
Let $\Omega \subset \R^n$ be a Lipschitz polytope. The Green's function, $\calG: \Omega \times \Omega \to \R$, defined in \eqref{eq:green}, satisfies:
\begin{enumerate}[1.]
  \item $\GRAD \calG \in L^{\frac{n}{n-1},\infty}(\Omega)$ which, in particular, implies that for all $q \in \left[1,\tfrac{n}{n-1}\right)$, we have that if we denote by $B_R$ the ball of radius $R$ centered at $x \in \Omega$ then
  \begin{equation}
  \label{eq:GreenW1p}
    \| \GRAD \calG \|_{L^q(B_R)} \lesssim R^{1-n+n/q},
  \end{equation}
  where the hidden constant depends on $q$ and $n$ and blows up as $q \uparrow n/(n-1)$.
  
  \item If, again, $B_R$ denotes the ball of radius $R$ and center $x \in \Omega$, then $\calG \in W^{2,1}(\Omega \setminus B_R)$ and satisfies
  \begin{equation}
  \label{eq:GreenW21}
    \| \mathcal{D}^2 \calG \|_{L^1(\Omega\setminus B_R)} \lesssim |\log R^{-1} |.
  \end{equation}
\end{enumerate}
\end{prpstn}
\begin{proof}
For $n = 3$, the estimates on the gradient of $\calG$ can be found in \cite[Theorem 1.1]{MR657523} and \cite[Theorem 4.1]{MR2341783}. In two dimensions, $n =2$, the weak-$L^2$ estimate is given by \cite[Theorem 1.1]{MR2763343}. Using the well known identity, see \cite[Exercise 1.1.11]{Grafakos} and  \cite[estimate (1.12)]{MR657523},
\[
  \| w \|_{L^{s-\epsilon}(U)} \leq \left( \frac{s}\epsilon \right)^{\frac1{s-\epsilon}} | U |^{\frac\epsilon{s(s-\epsilon)}} \|w\|_{L^{s,\infty}(U)}
\]
with $w = \GRAD \calG$, $s = 2$, $s-\epsilon = q$ and $U = B_R$ immediately yields
\[
  \| \GRAD \calG \|_{L^q(B_R)} \lesssim \left( \frac2{2-q} \right)^{\frac1q} R^{\frac{2(2-q)}{2q}} \| \GRAD \calG \|_{L^{2,\infty}(B_R)}.
\]
Since $\tfrac{2(2-q)}{2q} = 1-2 + 2/q$, the estimate above is \eqref{eq:GreenW1p}.

The estimates on the second derivatives of $\calG$ come from \cite[Lemma 2]{demlow2014maximum}.
\end{proof}

Identity \eqref{eq:point_representation}, in conjunction with Galerkin orthogonality and the bounds for the Green's function $\calG$ of Proposition~\ref{prop:propertiesG} are the main ingredients used to obtain a reliability property for $\E_{\infty}$. To state it, and for future reference, we define 
\begin{equation}\label{eq:log}
\ell_\T = \left|\log\left( \max_{T \in \T} \frac{1}{h_{T}}\right)\right|.
\end{equation}

\begin{lmm}[global reliability]
\label{lem:gob_rel_ptwise}
Let $u \in H_0^1(\Omega) \cap L^{\infty}(\Omega)$ and $u_{\T} \in \V(\T)$ be the solutions to problems \eqref{eq:lap_infty} and \eqref{eq:lap_infty_dis}, respectively. Then 
\begin{equation}
 \| u - u_{\T} \|_{L^{\infty}(\Omega)} \lesssim \ell_{\T} \E_{\infty}(u_{\T}; \T),
 \label{eq:global_reliability_linfty}
\end{equation}
where the hidden constant is independent of $u$, $u_{\T}$, the size of the elements in the mesh $\T$ and $\#\T$.
\label{lm:global_reliability_linfty}
\end{lmm}
\begin{proof} 
We follow \cite[Theorem 5.5]{Camacho01072015} and consider $x \in \Omega$ such that $|(u-u_\T)(x)|$ is maximized over $\Omega$. We write $\calG = \calG(x,\cdot)$ for the Green's function of \eqref{eq:green}. Then, invoking the pointwise representation \eqref{eq:point_representation} and Galerkin orthogonality we obtain that
\[
(u - u_\T)(x) = \int_{\Omega} \nabla (u - u_{\T}) \nabla \calG \diff y = \int_{\Omega} \nabla (u - u_{\T}) \nabla( \calG - \calG_{\T}) \diff y,
\]
where $\calG_{\T} \in \V(\T)$ denotes a suitable approximation of $\calG$, for instance, the Scott-Zhang interpolant \cite{SZ:90} or the interpolant based on local averages developed in \cite{NOS2}. Similar arguments to those used to conclude that \eqref{eq:aux_op_2} held can be applied to conclude that \eqref{eq:lap_infty} does in fact hold with $v=\calG$. Consequently,
\[
(u - u_\T)(x) = \sum_{T \in \T} \int_{T} f ( \calG - \calG_{\T}) \diff y + \sum_{S \in \Sides} \int_{S}  \llbracket \nabla u_{\T} \cdot \nu \rrbracket ( \calG - \calG_{\T}) = \textrm{I} + \textrm{II}
\]
upon integrating by parts.
We now proceed to control each term separately.

\noindent \framebox{Bound on $\textrm{I}$}: We begin with a simple application of the Cauchy Schwarz inequality
\[
  |\textrm{I}| \leq \sum_{T \in \T} \| f \|_{L^2(T)} \| \calG -\calG_\T \|_{L^2(T)}.
\]
Next we consider a partition of $\T$ into the sets $\Ne_x = \{ T \in \T: x \in \Ne_T\}$ and $\dot \Ne_x = \T \setminus \Ne_x$. On each one of these subsets we proceed as follows:
\begin{enumerate}[1.]
  \item If $T \in \dot \Ne_x$, then standard approximation results \cite{NOS2,SZ:90} yield
  \[
    \| \calG -\calG_\T \|_{L^2(T)} \lesssim h_T^{2-n/2} \|\mathcal{D}^2 \calG \|_{L^1(\Ne_T)}.
  \]
  Therefore, 
  we have that
  \begin{equation*}
    \sum_{T \in \dot \Ne_x} \| f \|_{L^2(T)} \| \calG - \calG_\T \|_{L^2(T)} \lesssim 
    \max_{ T \in \T} \left\{ h_T^{2-n/2} \| f \|_{L^2(T)} \right\} \sum_{T \in \dot \Ne_x} \| \mathcal{D}^2 \calG \|_{L^1(\Ne_T)} 
    \lesssim \ell_\T \E_\infty(u_\T; \T),
  \end{equation*}
  where, in the last step, we used the finite intersection property of stars and \eqref{eq:GreenW21} with $\displaystyle B_R$ being the largest ball such that $B_R \subset \cup_{T\in\T:\,x\in T}T$ in which case $\min_{T' \in \T} h_{T'}\lesssim R$.
  
  \item If $T \in \Ne_x$, we can estimate the difference $\calG - \calG_\T$, using approximation theory and \eqref{eq:GreenW1p}. In fact, since $2n/(n+2)<n/(n-1)$, we obtain the error estimate
  \[
    \| \calG -\calG_\T \|_{L^2(T)} \lesssim \| \GRAD \calG \|_{L^{2n/(n+2)}(\Ne_T)} \lesssim h_T^{2-n/2}.
  \]
  Therefore, since the quantity $\# \Ne_x$ is uniformly bounded, we conclude that
  \[
    \sum_{T \in \Ne_x} \| f \|_{L^2(T)} \| \calG - \calG_\T \|_{L^2(T)} \lesssim \sum_{T \in \Ne_x} h_T^{2-n/2} \| f \|_{L^2(T)}
    \lesssim
    \E_\infty(u_\T; \T).
  \]
\end{enumerate}
Gathering the estimates for these two cases we obtain the desired bound for $\textrm{I}$.

\noindent \framebox{Bound on $\textrm{II}$}: The ideas are similar to the ones used to control the term $\textrm{I}$. We begin with the estimate
\[
|\textrm{II}| \leq \sum_{S \in \Sides} \| \llbracket \nabla u_{\T} \cdot \nu \rrbracket \|_{L^\infty(S)} \| \calG -\calG_\T \|_{L^1(S)}.
\]
We now consider two cases based on the partition of $\T$ into the sets $\Ne_x$ and $\dot \Ne_x$.
\begin{enumerate}[1.]
  \item If $S \in \Sides$ is a side of $T$ and $T\in \dot \Ne_x$, then a scaled trace inequality yields
\begin{equation}
\label{eq:trace}
   \| \calG -\calG_\T \|_{L^1(S)} \lesssim h_T^{-1} \| \calG -\calG_\T \|_{L^1(T)} + \| \nabla(\calG -\calG_\T) \|_{L^1(T)},
\end{equation}

which, combined with standard approximation results \cite{NOS2,SZ:90}, yields
  \[
    \| \calG -\calG_\T \|_{L^1(S)} \lesssim h_T \|\mathcal{D}^2 \calG \|_{L^1(\Ne_T)}.
  \]
  Therefore, we have that
  \begin{equation*}
    \sum_{S\in\Sides_{T}:~T \in \dot \Ne_x } \| \llbracket \nabla u_{\T}\cdot \nu \rrbracket \|_{L^{\infty}(S)} \| \calG - \calG_\T \|_{L^1(S)} \lesssim 
    \E_\infty(u_\T; \T) \sum_{T \in \dot \Ne_x} \| \mathcal{D}^2 \calG \|_{L^1(\Ne_T)}
    \lesssim \ell_\T \E_\infty(u_\T; \T),
  \end{equation*}
  upon using the same arguments as before.
  
  \item If $S \in \Sides$ is a side of $T$ and $T \in \Ne_x$, we use that $(n+1)/n < n/(n-1)$ and the fact that $\calG \in W^{1,q}(\Omega)$ for $q < n/(n-1)$ to conclude that $\calG \in W^{1,(n+1)/n }(\Omega)$. Then, an application of the scaled trace inequality \eqref{eq:trace} and standard approximation estimates yield
  \[
    \| \calG -\calG_\T \|_{L^1(S)} \lesssim h_T^{n-n^2/(n+1)} \| \GRAD \calG \|_{L^{(n+1)/n}(\Ne_T)}.
  \]
  Since $\calG \in W^{1,(n+1)/n }(\Omega)$, \eqref{eq:GreenW1p} yields
  $\| \nabla \calG\|_{L^{(n+1)/n}(\Ne_T)} \lesssim h_T^{1-n+n^2/(n+1)}$. Then
\begin{equation*}
\sum_{S\in\Sides_{T}:~T \in \Ne_x} \| \llbracket \nabla u_{\T} \cdot \nu \rrbracket \|_{L^{\infty}(S)} \| \calG - \calG_\T \|_{L^1(S)} 
 \lesssim \sum_{T \in \Ne_x} h_T \| \llbracket \nabla u_{\T} \cdot \nu \rrbracket \|_{L^{\infty}(\partial T\setminus\partial\Omega)} 
 \lesssim
\E_\infty(u_\T; \T),
\end{equation*}
\end{enumerate}
where, in the last step, we have again used that the quantity $\# \Ne_x$ is uniformly bounded.

We finally collect the estimates obtained to bound the terms $\textrm{I}$ and $\textrm{II}$ and conclude the desired result \eqref{eq:global_reliability_linfty}.
\end{proof}

\subsection{Efficiency}
\label{sec:a_posteriori_infty_eff}
We now proceed to investigate the local efficiency properties of the error estimator \eqref{eq:defofEinf}--\eqref{eq:defofEinfglobal}. To accomplish this task, for any $g \in L^2(\Omega)$ and $\mathcal{M} \subset \T$ we define 
\begin{equation}
\label{eq:osc_localf}
 \textrm{osc}_{\T}(g;\mathcal{M}) =  \left( \sum_{T \in \mathcal{M}} h^{2(2-n/2)}_{T} \| g - \mathcal{P}_{\T}g \|_{L^{2}(T)}^2 \right)^{\frac{1}{2}},
\end{equation}
where $\mathcal{P}_{\T}$ denotes the $L^2$-projection operator onto piecewise linear functions over $\T$.

\begin{lmm}[local efficiency] 
\label{lm:local_infty}
Let $u \in H_0^1(\Omega)\cap L^\infty(\Omega)$ and $u_{\T} \in \V(\T)$ be the solutions to problems \eqref{eq:state_equation} and \eqref{eq:discrete_state_equation}, respectively. Then
\[
 \E (u_{\T};T) \lesssim \| u - u_{\T} \|_{L^{\infty}(\Ne_T^*)} + \mathrm{osc}_{\T}(f;\Ne_T^*)
\]
for all $T \in \T$, where the hidden constant is independent of $f$, $u$, $u_{\T}$, the size of the elements in the mesh $\T$ and $\#\T$.
\end{lmm}
\begin{proof}
 Consider $v \in H^1_0(\Omega)$ such that $v_{|T} \in C^2(T)$ for all $T \in \T$ as a test function in \eqref{eq:lap_infty}. Since $u \in H^1_0(\Omega)$ solves \eqref{eq:lap_infty}, integration by parts yields
\begin{equation}
\label{eq:infty_res1}
\int_{\Omega} \nabla (u - u_{\T}) \cdot \nabla v = \sum_{T \in \T} \int_{T} f v + \sum_{S \in \Sides } \int_{S } \llbracket \nabla u_{\T} \cdot \nu \rrbracket v
\end{equation}
and
\begin{equation}
\label{eq:infty_res2}
\int_{\Omega} \nabla (u - u_{\T}) \cdot \nabla v = - \sum_{T \in \T} \int_{T} (u - u_{\T}) \Delta v - \sum_{S \in \Sides } \int_{S } \llbracket \nabla v \cdot \nu \rrbracket (u - u_{\T}).
\end{equation}
In light of \eqref{eq:infty_res1}--\eqref{eq:infty_res2}, we proceed to estimate each term on the right hand side of \eqref{eq:defofEinf} separately. 

\noindent \framebox{Step 1.} We bound $h_T^{2-n/2}\| f \|_{L^2(T)}$. An application of the triangle inequality yields
\begin{equation}
\label{eq:triangle_eff}
 h_T^{2-n/2} \| f\|_{L^2(T)} \leq  h_T^{2-n/2} \| f -\mathcal{P}_{\T} f  \|_{L^2(T)} + h_T^{2-n/2} \| \mathcal{P}_{\T} f \|_{L^2(T)}.
\end{equation}
It thus suffices to bound the second term. To do this, we invoke the residual estimation techniques introduced by Verf\"{u}rth in \cite{Verfurth2,Verfurth}. Let $\varphi_T$ be the standard bubble function over $T$. Define $\beta_T:= \varphi_T^2\mathcal{P}_{\T} f $. Then, standard properties of the bubble function yield
\[
\|\mathcal{P}_{\T} f\|^2_{L^{2}(T)} \lesssim \|\varphi_T\mathcal{P}_{\T} f\|^2_{L^{2}(T)}= \int_{T} \mathcal{P}_{\T} f \beta_{T} .
\]
We now proceed to bound the term $ \int_T  \mathcal{P}_{\T} f \beta_{T} $. Set $v=\beta_T$ in \eqref{eq:infty_res1}. This, in view of the fact that $\beta_T|_S = 0$ for every $S \in \Sides$, allows us to obtain that
\begin{equation*}
\int_{T} \mathcal{P}_{\T} f \beta_{T}
 =
\int_{T} (\mathcal{P}_{\T} f - f)\beta_{T} +  \int_{T} f\beta_{T}
=
\int_{T} (\mathcal{P}_{\T} f - f)\beta_{T} + 
\int_{T} \nabla (u - u_{\T})  \nabla \beta_{T}.
\end{equation*}
Now, since $\nabla \beta_T = \varphi_T( 2  \nabla \varphi_T \mathcal{P}_T f  + \varphi_T \nabla \mathcal{P}_T f)$, we conclude that, for every $S \in \Sides$, we have
$
\int_{S } \llbracket \nabla \beta_{T} \cdot \nu \rrbracket (u - u_{\T})=0.
$
Thus, setting $v=\beta_T$ in \eqref{eq:infty_res2} yields that
\begin{equation*}
\int_{T} \mathcal{P}_{\T} f \beta_{T}
 =
\int_{T} (\mathcal{P}_{\T} f - f)\beta_{T}  
- \int_{T} (u - u_{\T}) \Delta \beta_{T}.
\end{equation*}

We control the first term on the right hand side of the previous expression by using
properties of the bubble function: $|\int_{T} (\mathcal{P}_{\T} f - f)\beta_{T}| \lesssim \|\mathcal{P}_{\T} f - f\|_{L^{2}(T)}\|\mathcal{P}_{\T} f\|_{L^{2}(T)}$. Now, to bound the second term we proceed as follows: $|\int_{T} (u - u_{\T}) \Delta \beta_{T}| \leq \|u - u_{\T}\|_{L^{\infty}(T)}\int_{T} |\Delta \beta_{T}|$. It thus suffices to control $\int_{T} |\Delta \beta_{T}|$. To accomplish this task, we use that  $\Delta \mathcal{P}_{\T} f=0$ on $T$, properties of the bubble function $\varphi_T$, and an inverse estimate. In fact,
\begin{align*}
\int_{T} |\Delta \beta_{T}|\leq & 2 \|\mathcal{P}_{\T} f\|_{L^{2}(T)}\left(\|\varphi_T\Delta\varphi_T\|_{L^{2}(T)} + \|\nabla\varphi_T\cdot\nabla \varphi_T\|_{L^{2}(T)}\right)
\\ &+ 4 \|\nabla\mathcal{P}_{\T} f\|_{L^{2}(T)}\|\varphi_T\nabla \varphi_T\|_{L^{2}(T)}
\lesssim
h_{T}^{n/2-2}\|\mathcal{P}_{\T} f\|_{L^{2}(T)}.
\end{align*}
Collecting all the previous findings allows us to conclude that
\begin{align*}
h_T^{2-n/2}\|\mathcal{P}_{\T} f\|_{L^{2}(T)} \lesssim 
h_T^{2-n/2}\|f-\mathcal{P}_{\T} f\|_{L^{2}(T)} +
\|u - u_{\T}\|_{L^{\infty}(T)}.
\end{align*}
Inserting the previous bound into \eqref{eq:triangle_eff} yields 
\begin{equation}\label{eq:ctm1_new}
h_T^{2-n/2} \| f\|_{L^2(T)} \lesssim  h_T^{2-n/2} \| f -\mathcal{P}_{\T} f  \|_{L^2(T)} + 
\|u - u_{\T}\|_{L^{\infty}(T)}.
\end{equation}


\noindent \framebox{Step 2.} We now control the jump terms $h_T \| \llbracket \nabla u_{\T}  \cdot \nu\rrbracket \|_{L^{\infty}(\partial T \setminus \partial \Omega)}$. To achieve this, given $S \in\Sides_T$, we set in \eqref{eq:infty_res1} and \eqref{eq:infty_res2} $v = \varphi_S$, the standard bubble function over $S$. This yields
\begin{equation*}
  \left| \int_S \llbracket  \nabla u_\T \cdot \nu \rrbracket \varphi_S \right|  \leq \sum_{T' \in \Ne_S} \int_{T'} |f| \varphi_S +
  \sum_{T' \in \Ne_S} \int_{T'} |u-u_\T| |\LAP \varphi_S|
+ \sum_{T' \in \Ne_S}\sum_{S' \in \Sides_{T'}}
   \int_{S'} |u-u_\T| |\llbracket \nabla \varphi_S \cdot \nu \rrbracket|.
\end{equation*}
Let us now recall that $\llbracket \nabla u_\T \cdot \nu \rrbracket$ is constant over $S$ and that, for $k=0,1,2$, we have $|\nabla^k \varphi_S| \approx h_S^{-k}$ to obtain
\begin{equation*}
  |S| \| \llbracket \nabla u_\T \cdot \nu \rrbracket \|_{L^\infty(S)}  \lesssim \sum_{T' \in \Ne_S}|T'|^{1/2} \|f \|_{L^2(T')}
  + \sum_{T' \in \Ne_S}\left(h_S^{-2} |T'| + h_S^{-1} \sum_{S' \in \Sides_{T'}}|S'|\right) \| u - u_\T \|_{L^\infty(T')}.
\end{equation*}
Upon multiplying this inequality by $h_T  |S|^{-1}$, it remains to use shape regularity to derive
\[
  h_T  \| \llbracket \nabla u_\T \cdot \nu \rrbracket \|_{L^\infty(S)} \lesssim h_T^{2-n/2} \sum_{T' \in \Ne_S}\| f\|_{L^2(T')} + \max_{T' \in \Ne_S}\| u - u_\T \|_{L^\infty(T')}.
\]
We conclude by using the bound \eqref{eq:ctm1_new} obtained in Step 1.

The  local efficiency is thus proved.
\end{proof}

\section{A posteriori error analysis}
\label{sec:a_posteriori}

The design and analysis of AFEMs to solve the optimal control problem with point sources is motivated by the fact that the error estimates \eqref{eq:estimate_1} and \eqref{eq:estimate_2} are not optimal in terms of approximation; these must be quadratic. These sub-optimal error estimates are expected and are a consequence of the reduced regularity properties of the optimal state $\bar{\ysf}$ solving problem \eqref{eq:state_equation}. In addition, AFEMs are also motivated by the fact that the estimates \eqref{eq:estimate_1} and \eqref{eq:estimate_2} require $\T$ to be quasi--uniform, $\Omega$ to be convex and high integrability assumptions on the desired state $\ysf_d$. In the search for an efficient method to solve the optimal control problem with point sources, in this section we propose and analyze an a posteriori error estimator to drive AFEMs. Our error indicator is built on the basis of the error estimator for elliptic problems involving Dirac measures elaborated in \cite{AGM} and the pointwise a posteriori error estimator allowing unbounded forcing terms investigated in section~\ref{sec:a_posteriori_infty}.

We derive and analyze an a posteriori error estimator for problem \eqref{eq:defofJ}--\eqref{eq:defofu}. The accomplishment of this task is not as simple as it may seem at first: the state equation \eqref{eq:state_equation} involves Dirac measures as a forcing term and the variational inequality \eqref{eq:VIp} that characterizes the optimal control involves the point evaluations of the optimal adjoint state. Consequently, the analysis of such an error estimator involves the interaction of $L^{\infty}(\Omega)$, $\mathbb{R}^l$ and the weighted Sobolev space $H^{1}(\rho,\Omega)$. This is one of the highlights of this work.

\subsection{The error estimator} 
On the basis of the notation introduced in section~\ref{sec:a_posteriori_infty}, we proceed to write an a posteriori error estimator for the optimal control problem with point sources. The error estimator is defined as the sum of two contributions:
\begin{equation}
 \label{eq:defofEocp}
\E_{\textrm{ocp}}^{2}(\bar{\ysf}_{\T},\bar{\psf}_{\T},\bar{\mathbf{u}}_{\T}; \T) = \E_{\ysf}^{2}(\bar{\ysf}_{\T},\bar{\mathbf{u}}_{\T}; \T) + \E_{\psf}^{2}(\bar{\psf}_{\T},\bar{\ysf}_{\T}; \T),
\end{equation}
where $\T \in \Tr$ and $\bar \ysf_{\T}$, $\bar{\mathbf{u}}_{\T}$ and $\bar \psf_{\T}$ denote the discrete optimal variables solving the finite element counterpart of \eqref{eq:defofJ}--\eqref{eq:defofu} described in section~\ref{sub:fem}. We describe each error indicator in \eqref{eq:defofEocp} separately, starting with $\E_{\ysf}$. To define it, we assume that 
\begin{equation}
\label{eq:patch}
\forall T \in \T, \ \# (\Ne_T \cap D) \leq 1,
\end{equation}
\ie for each element $T\in\T$ its patch $\Ne_T$ contains at most one source point $z \in D$. We comment that this assumption is not restrictive since it can always be satisfied by starting with a suitably refined mesh. 

Inspired by \cite{AGM,AOSR}, we then define the local error indicator $\E_{\ysf}^2(\bar{\ysf}_{\T},\bar{\mathbf{u}}_{\T}; T)$ as
\begin{equation}
 \label{eq:defofEy}
 \E_{\ysf}^2(\bar{\ysf}_{\T},\bar{\mathbf{u}}_{\T}; T)= h_T D_T^{\alpha} \| \llbracket \nabla \bar{\ysf}_{\T} \cdot \nu \rrbracket\|^2_{L^2(\partial T\setminus \partial \Omega)} +  \sum_{z \in D \cap T} h_T^{\alpha + 2 - n} |\bar{\usf}_{\T,z}|^2,
\end{equation}
where, as usual, $\sum_{\emptyset} = 0$,  the interelement residual $\| \llbracket \nabla \bar{\ysf}_{\T} \cdot \nu \rrbracket\|$ is defined by \eqref{eq:jump}, $\alpha \in (n-2,2)$ and 
\begin{equation}
\label{eq:DT}
  D_T := \min_{z \in D} \left\{  \max_{x \in T} |x-z| \right\}.
\end{equation}
The global error estimator $\E_{\ysf}(\bar{\ysf}_{\T},\bar{\mathbf{u}}_{\T}; \T)$ is thus defined by
\begin{equation}
\label{eq:defofEyglobal}
\E_{\ysf}(\bar{\ysf}_{\T},\bar{\mathbf{u}}_{\T}; \T) = \left(  \sum_{T \in \T}  \E^2_{\ysf}(\bar{\ysf}_{\T},\bar{\mathbf{u}}_{\T}; T) \right)^{\frac{1}{2}}. 
\end{equation}
Notice that the range of $\alpha$ is as in Lemma~\ref{lem:restalphagen}.

The second error contribution in \eqref{eq:defofEocp} is based on the maximum norm error indicator that we developed and analyzed in section~\ref{sec:a_posteriori_infty}. Locally it is defined by
\begin{equation}
 \label{eq:defofEp}
 \E_{\psf}(\bar{\psf}_{\T},\bar{\ysf}_{\T}; T)= h_T^{2-\tfrac{n}{2}} \| \bar{\ysf}_{\T} - \ysf_d \|_{L^{2}(T)} 
 + h_T \| \llbracket \nabla \bar{\psf}_{\T} \cdot \nu \rrbracket  \|_{L^{\infty}(\partial T \setminus \partial \Omega)}.
\end{equation}
The global pointwise error estimator $\E_{\psf}(\bar{\psf}_{\T},\bar{\ysf}_{\T}; \T)$ is then defined by
\begin{equation}
\label{eq:defofEpglobal}
\E_{\psf}(\bar{\psf}_{\T},\bar{\ysf}_{\T}; \T) = \max_{T \in \T } \E_{\psf}(\bar{\psf}_{\T},\bar{\ysf}_{\T}; T).
\end{equation}

\subsection{Error estimator: reliability}
\label{subsub:reliable}
In this section we follow the arguments of \cite[Theorem 2]{AOSR} and derive a global reliability property for the error estimator $\E_{\textrm{ocp}}$.

\begin{thrm}[global reliability property of $\E_{\textrm{ocp}}$]
\label{thm:globrelaible}
Let $(\bar{\mathbf{u}},\bar{\ysf},\bar{\psf}) \in \mathcal{U}_{\textrm{ad}} \times H_0^1(\rho, \Omega) \times H_0^1(\Omega)$ be the solution to the optimality system \eqref{eq:optimal_system} associated with the optimal control problem with point sources and $(\bar{\mathbf{u}}_{\T},\bar{\ysf}_{\T},\bar{\psf}_{\T}) \in \mathcal{U}_{\textrm{ad}}\times \V(\T) \times \V(\T)$ be its numerical approximation given by \eqref{eq:discrete_state_equation}--\eqref{eq:discrete_adjoint_equation}. If $\alpha \in (n-2,2)$ and $\ysf_d \in L^r(\Omega)$, for every $r < 3$, then 
\begin{align}
\| \bar{\mathbf{u}} - \bar{\mathbf{u}}_{\T} \|_{\R^l}^{2}  + \| \nabla( \bar{\ysf} - \bar{\ysf}_{\T}) \|_{L^2(\rho,\Omega)}^{2} 
+ \| \bar{\psf} - \bar{\psf}_{\T} \|_{L^{\infty}(\Omega)}^{2}  
& \lesssim \E_{\ysf}^{2}(\bar{\ysf}_{\T},\bar{\mathbf{u}}_{\T}; \T)+ \ell_{\T}^{2} \E_{\psf}^{2}(\bar{\psf}_{\T},\bar{\ysf}_{\T}; \T) \nonumber
\\
& \lesssim (1 + \ell_{\T}^{2}) \E_{\textrm{\emph{ocp}}}^{2} (\bar{\ysf}_{\T},\bar{\psf}_{\T},\bar{\mathbf{u}}_{\T}; \T), \label{eq:reliability}
\end{align}
where $\ell_\T$ is defined in \eqref{eq:log} and the hidden constant is independent of the continuous and discrete optimal variables, the size of the elements in the mesh $\T$ and $\#\T$.
\label{TH:global_reliability}
\end{thrm}
\begin{proof}
We proceed in six steps.

\noindent \framebox{Step 1.} First, we notice that the discrete structure of the set $\mathcal{U}_{\textrm{ad}}$ allows us to consider $\usf_{z} = \bar{\usf}_{z}$ in \eqref{eq:first_oder_discrete}.  Second, we set $\usf_{z} = \bar{\usf}_{\T,z}$ in \eqref{eq:VIp}. Adding the obtained variational inequalities we obtain that
\begin{equation}
\label{eq:usf-usfT}
 \lambda \| \bar{\mathbf{u}} - \bar{\mathbf{u}}_{\T} \|^2_{\R^l} \leq \sum_{z \in D} (\bar{\psf}(z) - \bar{\psf}_{\T}(z))(\bar{\usf}_{\T,z} - \bar{\usf}_{z}).
\end{equation}

\noindent \framebox{Step 2.} The goal of this step is to bound the right hand side of \eqref{eq:usf-usfT}. To accomplish this task, we define an auxiliary adjoint state via the following weak problem:
\begin{equation}
\label{eq:adjq}
\qsf \in H_0^1(\Omega): \quad (\nabla \wsf,\nabla \qsf)_{L^2(\Omega)} = (\bar{\ysf}_{\T} - \ysf_{d}, \wsf )_{L^2(\Omega)}\quad \forall \wsf \in H^1_0(\Omega).
\end{equation}
With this auxiliary adjoint state at hand, we write $\bar{\psf} - \bar{\psf}_{\T} = (\bar{\psf}  - \qsf) + (\qsf - \bar{\psf}_{\T})$. Then, on the basis of \eqref{eq:usf-usfT}, we arrive at 
\begin{equation}
\label{I+II}
 \lambda \| \bar{\mathbf{u}} - \bar{\mathbf{u}}_{\T} \|^2_{\R^l} \leq 
 \sum_{z \in D} \left[ (\bar{\psf}(z) - \qsf(z))
 + (\qsf(z) - \bar{\psf}_{\T}(z))\right](\bar{\usf}_{\T,z} - \bar{\usf}_{z}) 
  = \textrm{I} + \textrm{II}.
\end{equation}
The rest of this step is dedicated to control the term $\textrm{II}$. To do this, we exploit that $\bar{\psf}_{\T}$ is the Galerkin approximation of the state $\qsf$ that solves \eqref{eq:adjq}. This property, an application of the Cauchy-Schwarz and Young's inequalities, and the global reliability of the error indicator $\E_{\psf}$ obtained in Lemma~\ref{lem:gob_rel_ptwise} yield the estimate 
\begin{equation}
\label{II}
 |\textrm{II}| \leq l^{1/2} \| \qsf  - \bar{\psf}_{\T} \|_{L^{\infty}(\Omega)}\| \bar{\mathbf{u}} - \bar{\mathbf{u}}_{\T} \|_{\R^l} \leq C \ell_\T^{2} \E^2_{\psf}(\bar{\psf}_{\T},\bar{\ysf}_{\T}; \T ) + \frac{\lambda}{4}\| \bar{\mathbf{u}} - \bar{\mathbf{u}}_{\T} \|^2_{\R^l},
\end{equation}
where $\ell_{\T}$ is defined in \eqref{eq:log}. The constant $C$ in independent of the optimal continuous and discrete optimal variables and the size of the elements in the mesh $\T$.

\noindent \framebox{Step 3.} The goal of this step is to bound the term $\textrm{I}$ in \eqref{I+II}. To accomplish this task we introduce two auxiliary states. We first define
\begin{equation}
\label{eq:state_equation_aux}
\tilde \ysf \in H^1_0(\rho, \Omega): \quad  (\nabla \tilde \ysf, \nabla \vsf)_{L^2(\Omega)} = \sum_{z \in D } \bar \usf_{\T,z} \langle \delta_z, \vsf \rangle \quad \forall \vsf \in H^1_0(\rho^{-1},\Omega),
\end{equation}
and $\rsf \in H^1_0(\Omega)$ that solves
\[
(\nabla \wsf,\nabla \rsf)_{L^2(\Omega)} = ( \tilde \ysf - \ysf_d, \wsf)_{L^2(\Omega)} \quad \forall \wsf \in H^1_0(\Omega).
\]
With this notation, we then write $\mathrm{I}$ as the sum of two terms:
\[
  \mathrm{I} = \mathrm{I}_a + \mathrm{I}_b := \sum_{z \in D} (\bar{\psf}(z) - \rsf(z))(\bar{\usf}_{\T,z} - \bar{\usf}_z )
    + \sum_{z \in D} (\rsf(z) - \qsf(z) )(\bar{\usf}_{\T,z} - \bar{\usf}_z)
\]
and bound each one of them separately.

The bound of $\mathrm{I}_a$ borrows from ideas in \cite{AOSR}. 
We first observe that the difference $\bar{\psf} - \rsf \in H^1_0(\Omega)$ satisfies
\begin{equation}
  (\nabla \vsf, \nabla (\bar{\psf} - \rsf) )_{L^2(\Omega)} = ( \bar\ysf - \tilde\ysf, \vsf)_{L^2(\Omega)} \quad \forall \vsf \in H^1_0(\Omega),
  \label{eq:aux_rel}
\end{equation}
with forcing term $\bar\ysf - \tilde \ysf \in L^r(\Omega)$, for every $r<r^{*}$; $r^{*}$ being defined in \eqref{eq:r_star}. We thus apply the results of Proposition \ref{pro:weighted_integra} to conclude that $\bar{\psf} - \rsf \in H_0^1(\rho^{-1},\Omega)$. In addition, invoking Proposition \ref{prop:uisW1p} we conclude that $\bar \psf - \rsf$ is H\"older continuous in $\bar \Omega$ and thus that its point evaluations are well--defined. On the basis of these results, in the equation that $\bar\ysf - \tilde\ysf$ solves, \ie
\[
  (\nabla (\bar\ysf - \tilde\ysf), \nabla \vsf )_{L^2(\Omega)} = \sum_{z \in D} (\bar{\usf}_z - \bar{\usf}_{\T,z}) \langle \delta_z, \vsf \rangle \quad \forall \vsf \in H^1_0(\rho^{-1},\Omega),
\]
that is obtained from \eqref{eq:optimal_system} and \eqref{eq:state_equation_aux}, it is admissible to set $\vsf =\bar{\psf} - \rsf $. This yields
\[
  (\nabla (\bar\ysf - \tilde\ysf), \nabla (\bar \psf - \rsf) )_{L^2(\Omega)} = \sum_{z \in D} (\bar{\usf}_z - \bar{\usf}_{\T,z}) (\bar \psf(z) - \rsf(z))
\]
and hence $\mathrm{I}_a=-(\nabla (\bar \ysf - \tilde\ysf), \nabla (\bar \psf - \rsf))_{L^2(\Omega)}$. On the other hand, applying a similar approximation argument to that used in the proof of Theorem \ref{thm:optim_conds} to \eqref{eq:aux_rel} allows us to arrive at
\[
 (\nabla (\bar \ysf - \tilde\ysf), \nabla (\bar \psf - \rsf))_{L^2(\Omega)} = (\bar \ysf - \tilde \ysf, \bar \ysf - \tilde\ysf)_{L^2(\Omega)}.
\]
We have thus obtained that
\[
 \mathrm{I}_a = - \|  \bar \ysf - \tilde \ysf \|^2_{L^2(\Omega)} \leq 0.
\]

To estimate $\mathrm{I}_b$ we observe two things. First, since $\Omega$ is Lipschitz,
we have, using Proposition \ref{prop:uisW1p}, the existence of $q>n$ for which
\[
  \| \rsf - \qsf  \|_{W^{1,q}(\Omega)} \lesssim \| \tilde{\ysf} - \bar{\ysf}_{\T} \|_{L^2(\Omega)}.
\]
This, in view of the continuous embedding $W^{1,q}(\Omega) \hookrightarrow C(\bar \Omega)$ for $q > n$, implies that
\begin{equation}
\label{eq:W1pestimate}
 \| \rsf - \qsf  \|_{L^{\infty}(\Omega)} \lesssim \| \tilde{\ysf} - \bar{\ysf}_{\T} \|_{L^2(\Omega)}.
\end{equation}
Second, we observe that $\bar{\ysf}_\T$ is the Galerkin approximation of $\tilde\ysf$ that solves \eqref{eq:state_equation_aux}. Thus, an adaptation of the arguments developed in \cite[Theorem 5.1]{AGM}, combined with the embedding of Lemma \ref{lem:restalphagen}, yields that
\begin{equation}
\label{eq:aux1_reliability}
\|  \tilde\ysf - \bar{\ysf}_{\T}\|_{L^{2}(\Omega)} \lesssim \| \nabla( \tilde\ysf - \bar{\ysf}_{\T} )\|_{L^{2}(\rho,\Omega)} \lesssim   \E_{\ysf}(\bar{\ysf}_{\T},\bar{\mathbf{u}}_{\T};\T),
\end{equation}
where $\E_{\ysf}$ denotes the a posteriori error estimator defined by \eqref{eq:defofEy} and \eqref{eq:defofEyglobal}. For brevity we skip details and only remark that this estimate is valid because of assumption \eqref{eq:patch} and the given range of $\alpha$.

Collecting \eqref{eq:W1pestimate}, \eqref{eq:aux1_reliability} and the derived estimate for $\mathrm{I}_a$ we obtain
\[
 \mathrm{I} \leq \mathrm{I}_b \leq l^{1/2} \| \rsf - \qsf  \|_{L^{\infty}(\Omega)} \|\bar{\mathbf{u}} - \bar{\mathbf{u}}_{\T} \|_{\R^l} \leq \frac{\lambda}{4}\| \bar{\mathbf{u}} - \bar{\mathbf{u}}_{\T} \|^2_{\R^l}  + C \E_{\ysf}^2(\bar{\ysf}_{\T},\bar{\mathbf{u}}_{\T};\T),
\]
again with a constant $C$ that is independent of the continuous and discrete optimal variables and the size of the elements of the mesh $\T$. This, in conjunction with the estimate \eqref{II} for $\mathrm{II}$, yields 
\begin{equation}
\label{eq:a_post_u}
  \| \bar{\mathbf{u}} - \bar{\mathbf{u}}_{\T} \|_{\R^l}^2 \lesssim \E_{\ysf}^2(\bar{\ysf}_{\T},\bar{\mathbf{u}}_{\T};\T) + \ell_{\T}^2 \E_{\psf}^2(\bar{\psf}_{\T},\bar{\ysf}_{\T};\T).
\end{equation}

\noindent \framebox{Step 4.} In this step we bound the error $\bar{\ysf} - \bar{\ysf}_{\T}$, in the $H^{1}(\rho,\Omega)$-seminorm, in terms of the estimator $\E_{\textrm{ocp}}$. We follow the ideas developed in Step 3 and write $\bar{\ysf} - \bar{\ysf}_{\T} = ( \bar{\ysf} - \tilde\ysf ) + (\tilde\ysf- \bar{\ysf}_{\T})$ and estimate each term separately. The first term is controlled in light of the well-posedness of the state equation \eqref{eq:state_equation} \cite[Theorem 1]{AGM}:
\[
  \| \nabla(\bar{\ysf} - \tilde\ysf)  \|_{L^2(\rho,\Omega)} \lesssim \| \bar{\mathbf{u}}- \bar{\mathbf{u}}_{\T} \|_{\R^l}. 
\]
The second term, \ie~ the difference $\tilde\ysf- \bar{\ysf}_{\T}$ is controlled by invoking \eqref{eq:aux1_reliability}.
Then, in view of \eqref{eq:a_post_u}, the collection of the derived estimates implies that
\begin{equation}
\label{eq:aux2_reliability}
  \| \nabla(\bar{\ysf} - \bar{\ysf}_{\T} ) \|_{L^2(\rho,\Omega)}^2
    \lesssim
  \E_{\ysf}^2(\bar{\ysf}_{\T},\bar{\mathbf{u}}_{\T};\T) + \ell_{\T}^2 \E_{\psf}^2(\bar{\psf}_{\T},\bar{\ysf}_{\T};\T).
\end{equation}

\noindent \framebox{Step 5.} In this step we bound the term $\bar{\psf} - \bar{\psf}_{\T}$. In the spirit of Step 4, we write  $\bar{\psf} - \bar{\psf}_{\T} = ( \bar{\psf} - \qsf ) + ( \qsf - \bar{\psf}_{\T})$ with $\qsf$ solving \eqref{eq:adjq}. The control of the first term follows the arguments detailed in 
Step 3:
\[
\| \bar{\psf} - \qsf \|_{L^{\infty}(\Omega)} \lesssim \|  \bar\ysf - \bar{\ysf}_{\T}\|_{L^{2}(\Omega)} \lesssim \| \nabla( \bar\ysf - \bar{\ysf}_{\T} )\|_{L^{2}(\rho,\Omega)}
\]
which can be bounded using \eqref{eq:aux2_reliability}. On the other hand, from Lemma~\ref{lem:gob_rel_ptwise} we conclude that
$
\| \qsf - \bar{\psf}_{\T}\|_{L^{\infty}(\Omega)} \lesssim \ell_{\T}\E_{\psf}(\bar{\psf}_{\T},\bar{\ysf}_{\T} ;\T).
$
Collecting the estimates we arrive at
\begin{equation}
\label{eq:aux3_reliability}
\| \bar{\psf} - \bar{\psf}_{\T} \|_{L^{\infty}(\Omega)}^2 \lesssim \E_{\ysf}^2(\bar{\ysf}_{\T},\bar{\mathbf{u}}_{\T} ;\T) + \ell_{\T}^2\E_{\psf}^2(\bar{\psf}_{\T},\bar{\ysf}_{\T} ;\T).
\end{equation}

\noindent \framebox{Step 6.} The collection of the estimates \eqref{eq:a_post_u}, \eqref{eq:aux2_reliability} and \eqref{eq:aux3_reliability} yields the desired a posteriori error estimate \eqref{eq:reliability}. This concludes the proof.
\end{proof}

\subsection{Error estimator: efficiency}
\label{subsub:efficient}
In this section we analyze the efficiency properties of the error estimator $\E_{\textrm{ocp}}$ defined in \eqref{eq:defofEocp}. To accomplish this task, we examine each of its contributions separately. We start with the indicator $\E_{\ysf}^2 (\bar{\ysf}_{\T},\bar{\mathbf{u}}_{\T}; T)$ defined by \eqref{eq:defofEy}. A key ingredient in its efficiency analysis is an abstract estimate for the \emph{residual} $\Res_{\ysf} = \Res_{\ysf}(\bar{\ysf}_{\T}) \in H_0^1(\rho^{-1},\Omega)'$ which, for all $\vsf \in H^1_0(\rho^{-1},\Omega)$, is defined by
\begin{equation}
\label{Reseq}
\langle \Res_{\ysf} (\bar{\ysf}_{\T}),\vsf \rangle := (\nabla (\bar{\ysf} - \bar{\ysf}_{\T}),\nabla \vsf)_{L^2(\Omega)} = \sum_{z \in D} \bar{\usf}_z \langle \delta_z, \vsf \rangle - (\nabla \bar{\ysf}_{\T}, \nabla \vsf)_{L^2(\Omega)}.
\end{equation}
The aforementioned abstract estimate reads as follows: If $\Orara$ denotes a subdomain of $\Omega$ and $\vsf \in H^1_0(\rho^{-1},\Orara)$, then 
\begin{equation}
\label{eq:residual}
 | \langle \Res_{\ysf} (\bar{\ysf}_{\T}),\vsf \rangle | \leq \| \nabla ( \bar \ysf - \bar{\ysf}_{\T}) \|_{L^2(\rho,\Orara)}
  \| \nabla \vsf \|_{L^2(\rho^{-1},\Orara)}.
\end{equation}

We now utilize standard residual estimation techniques \cite{Verfurth2,Verfurth}, that have, as a key element in the analysis, the existence of a suitable bubble function. Given $S \in \Sides$, we introduce a bubble function $\psi_S$ (whose construction we owe to \cite{AGM}) that satisfies the following properties: $\psi_S(z) = 0$ for all $z \in D$,
\begin{equation}
\label{eq:psi_S}
 |S| \lesssim \int_{S} \psi_{S}, \quad \|\nabla \psi_{S} \|_{L^{2}(R_{S})} \lesssim h_{T}^{-1/2}|S|^{1/2},
\end{equation}
where $R_S = \supp \psi_S$. Moreover, if $\Ne_S = \{ T, T' \}$, there are simplices $T_{*} \subset T$ and $T'_{*} \subset T'$ such that $R_S \subset T_{*} \cup T_{*}'$.
We refer the reader to \cite[Section 5.2]{AGM} for details. We comment that, under assumption \eqref{eq:patch}, the construction of \cite[Section 5.2]{AGM} guarantees
\begin{equation}\label{DTbound}
  D_T \lesssim \min_{z \in D} \left\{ \min_{x \in T_*} |x-z| \right\}\mbox{ and }D_T \lesssim \min_{z \in D} \left\{ \min_{x \in T_*'} |x-z| \right\}.
\end{equation}

With all these ingredients at hand, we are ready to prove the local efficiency of $\E_{\ysf}^2 (\bar{\ysf}_{\T},\bar{\mathbf{u}}_{\T}; T)$. The proof is based on the arguments of \cite[Theorem 5.3]{AGM} and \cite[Lemma 5]{AOSR}.

\begin{lmm}[local efficiency of $\E_{\ysf}$] 
Let $( \bar{ \mathbf{u} },\bar{\ysf},\bar{\psf} ) \in \mathcal{U}_{\textrm{ad}} \times H_0^1(\rho,\Omega) \times H_0^1(\Omega)$ be the solution to the optimality system \eqref{eq:optimal_system} associated with the optimal control with point sources and $(\bar{\mathbf{u}}_{\T},\bar{\ysf}_{\T},\bar{\psf}_{\T}) \in \mathcal{U}_{\textrm{ad}}\times \V(\T) \times \V(\T)$ be its numerical approximation given by \eqref{eq:discrete_state_equation}--\eqref{eq:discrete_adjoint_equation}. If $\alpha \in (n-2,2)$, then 
\begin{equation}
\E_{\ysf}^2 (\bar{\ysf}_{\T},\bar{\mathbf{u}}_{\T}; T) \lesssim \| \nabla( \bar{\ysf} - \bar{\ysf}_{\T}) \|_{L^{2}(\rho,\Ne_T)}^2 
+
h_T^{\alpha + 2 - n} \sum_{z\in T \cap D}|\bar{\usf}_{z} - \bar{\usf}_{\T,z}|^2,
\label{eq:efficiencyy}
\end{equation}
where $\Ne_T$ is defined as in \eqref{eq:NeT} and the hidden constant is independent of the optimal variables, their approximations, the size of the elements in the mesh $\T$ and $\#\T$. 
\label{le:local_eff_y}
\end{lmm}
\begin{proof}  Let $T \in \T$ and $S \in \Sides$ be a side of $T$. We start the proof by bounding the term $h_T D_T^{\alpha} \| \llbracket  \nabla \bar{\ysf}_{\T} \cdot \nu \rrbracket\|^2_{L^2(\partial T\setminus \partial \Omega)}$ in \eqref{eq:defofEy}. To do this, we first invoke the bubble function $\psi_S$ and property \eqref{eq:psi_S} to obtain that
\begin{equation}
\label{eq:aux_effy}
 \| \llbracket   \nabla \bar{\ysf}_{\T} \cdot \nu \rrbracket\|^2_{L^2(S)} \lesssim 
 \int_{S}  \llbracket  \nabla \bar{\ysf}_{\T} \cdot \nu \rrbracket^2  \psi_S
 = \int_{S}  \llbracket \nabla \bar{\ysf}_{\T} \cdot \nu \rrbracket  \phi_S,
\end{equation}
where $\phi_S :=  \llbracket \nabla \bar{\ysf}_{\T} \cdot \nu \rrbracket \psi_S$. We thus utilize the properties $\supp \psi_{S} \subset T_* \cup T'_* \subset \Ne_S$ and $\psi_S (z) = 0$ for all $z \in D$, to obtain that $\int_S \llbracket \nabla \bar{\ysf}_{\T} \cdot \nu \rrbracket  \phi_S  = \langle \Res_{\ysf}(\bar{\ysf}_{\T}), \phi_S \rangle $ upon letting $\vsf=\phi_S$ in \eqref{Reseq} and integrating by parts. To bound the term $\langle \Res_{\ysf}(\bar{\ysf}_{\T}), \phi_S \rangle$, we first use \eqref{DTbound} along with the arguments used to arrive at \cite[equation (5.9)]{AGM} to conclude that
\[
\| \nabla \phi_S \|_{L^2(\rho^{-1},R_S)} \lesssim h_T^{-\frac{1}{2}}D_T^{-\frac{\alpha}{2}}\| \llbracket \nabla \bar{\ysf}_{\T} \cdot \nu \rrbracket\|_{L^2(S)}.
\]
Thus, in view of the abstract estimate \eqref{eq:residual} with $\Orara = R_S$, we arrive at
\[
| \langle \Res_{\ysf}(\bar{\ysf}_{\T}), \phi_S \rangle | \lesssim h_T^{-\frac{1}{2}} D_T^{-\frac{\alpha}{2}} \| \nabla ( \bar \ysf - \bar{\ysf}_{\T}) \|_{L^2(\rho,R_S)}\| \llbracket \nabla \bar{\ysf}_{\T} \cdot \nu \rrbracket\|_{L^2(S)}.
\]
This, in light of \eqref{eq:aux_effy}, immediately yields the estimate 
\begin{align}
\label{eq:jump_y}
h_T D_T^{\alpha} \| \llbracket \nabla \bar{\ysf}_{\T} \cdot \nu \rrbracket\|^2_{L^2(S)} \lesssim \sum_{T'\in\Ne_S}\| \nabla ( \bar \ysf - \bar{\ysf}_{\T}) \|_{L^2(\rho,T')}^2. 
\end{align}

It only remains to bound the term $\sum_{z \in D \cap T} h_T^{\alpha + 2 - n} |\bar{\usf}_{\T,z}|^2$ in \eqref{eq:defofEy}. In view of assumption \eqref{eq:patch}, we have that $T \cap D$ is either empty or consists of exactly one point. If $T \cap D = \emptyset$, then the desired estimate \eqref{eq:efficiencyy} follows immediately from \eqref{eq:jump_y}. If $T \cap D = \{ z \}$, then the estimator $\E_{\ysf}$ contains the term $h_T^{\alpha + 2 - n} |\bar{\usf}_{\T,z}|^2$. We thus proceed as follows: standard inequalities yield that
\begin{equation}
\label{eq:paso_vaca}
 | \bar{\usf}_{\T,z} |^2 \lesssim |\bar{\usf}_{\T,z} - \bar{\usf}_{z}|^2 + |\bar{\usf}_{z}|^2.
\end{equation}
The control of $|\bar{\usf}_{z}|$ follows from the arguments of \cite[Lemma 5]{AOSR} and \cite[Theorem 5.3]{AGM}, which rely on the existence of a suitable smooth function $\chi$ satisfying
\begin{equation}
\label{eq:chi}
\chi(z) = 1, \quad \| \chi \|_{L^{\infty}(\Omega)} = 1, \quad \| \nabla \chi \|_{L^{\infty}(\Omega)} = h_T^{-1}, \quad \supp \chi \subset \Ne_T.
\end{equation}
In fact, utilizing the first equation in \eqref{eq:optimal_system} in conjunction with $\chi(z) = 1$, $\textrm{supp } \chi \subset \Ne_T$, assumption \eqref{eq:patch} and integration by parts, we arrive at
\begin{align*}
 |\bar{\usf}_{z}| & = |\bar{\usf}_{z} \chi(z) | \leq | (\nabla (\bar \ysf - \bar \ysf_\T, \nabla \chi)_{L^2(\Omega)}| + |(\nabla \bar \ysf_\T, \nabla \chi)_{L^2(\Omega)}|
 \\ & \leq \|\nabla(\bar{\ysf} - \bar{\ysf}_{\T}) \|_{L^2(\rho,\Ne_{T})} \| \nabla \chi \|_{L^2(\rho^{-1},\Ne_{T})} 
+  \sum_{ \substack{T' \in \T:\\ T' \subset \Ne_T} } 
 \sum_{ \substack{S \in \Sides_{T'}:\\ S \not\subset \partial\Ne_T} } 
 \| \llbracket  \nabla \bar{\ysf}_{\T} \cdot \nu \rrbracket \|_{L^2(S)} \| \chi \|_{L^2(S)}.
\end{align*}
Then, we utilize $\| \chi \|_{L^2(S)} \lesssim h_{T}^{\frac{n-1}{2}}$ and $\| \nabla \chi \|_{L^2(\rho^{-1},\Ne_{T})} \lesssim h_T^{\frac{n-2}{2} - \frac{\alpha}{2}}$ (see \cite[Theorem 5.3]{AGM} for details) and conclude that
\[
 h_T^{\alpha + 2 - n}  |\bar{\usf}_{z}|^2 \lesssim \|\nabla(\bar{\ysf} - \bar{\ysf}_{\T}) \|^2_{L^2(\rho,\Ne_{T})} + 
 \sum_{ \substack{T' \in \T:\\ T' \subset \Ne_T} } 
 \sum_{ \substack{S \in \Sides_{T'}:\\ S \not\subset \partial\Ne_T} } 
 h_{T'} D_{T'}^{\alpha} \| \llbracket \nabla \bar{\ysf}_{\T} \cdot \nu \rrbracket\|^2_{L^2(S)}.
\]
This, in conjunction with \eqref{eq:jump_y} and \eqref{eq:paso_vaca}, yields the desired estimate \eqref{eq:efficiencyy}. 
\end{proof}

\begin{rmrk}[range of $\alpha$] \rm
\label{rk:alpha}
Since $\alpha \in (n-2,2)$, we immediately deduce that
$
  \alpha + 2 -n >0.
$
Consequently, \eqref{eq:efficiencyy} is indeed an efficiency bound.
\end{rmrk}

We now continue with the study of the local efficiency properties of the indicator $\E_{\psf}$ defined by \eqref{eq:defofEp}.

\begin{lmm}[local efficiency of $\E_{\psf}$] 
Let $( \bar{ \mathbf{u} },\bar{\ysf},\bar{\psf} ) \in \mathcal{U}_{\textrm{ad}} \times H_0^1(\rho,\Omega) \times H_0^1(\Omega)$ be the solution to the optimality system \eqref{eq:optimal_system} associated with the optimal control with point sources and $(\bar{\mathbf{u}}_{\T},\bar{\ysf}_{\T},\bar{\psf}_{\T}) \in \mathcal{U}_{\textrm{ad}}\times \V(\T) \times \V(\T)$ be its numerical approximation given by \eqref{eq:discrete_state_equation}--\eqref{eq:discrete_adjoint_equation}. If $\alpha \in (n-2,2)$, then 
\begin{equation}
\label{eq:efficiencyp}
\E_{\psf} (\bar{\psf}_{\T},\bar{\ysf}_{\T}; T) \lesssim \|  \bar{\psf} - \bar{\psf}_{\T} \|_{L^{\infty}(\Ne_T^*)} + h_{T}^{2-n/2} \| \bar{\ysf} - \bar{\ysf}_\T \|_{L^2(\Ne_T^*)} 
+ \mathrm{osc}_{\T}(\ysf_d;\Ne_T^*),
\end{equation}
where $\Ne_T^*$ is defined as in \eqref{eq:NeTstar} and the hidden constant is independent of the optimal variables, their approximations, the size of the elements in the mesh $\T$ and $\#\T$. 
\label{le:local_eff_p}
\end{lmm}
\begin{proof}
The proof closely follows the arguments of Lemma \ref{lm:local_infty}. Let $\vsf \in H^1_0(\Omega)$ be such that $\vsf_{|T} \in C^2(T)$ for all $T \in \T$. Using \eqref{eq:optimal_system} and integration by parts we obtain
\begin{equation*}
\label{eq:res1}
\int_{\Omega} \nabla \vsf \cdot \nabla (\bar{\psf} - \bar{\psf}_{\T}) = \sum_{T \in \T} \int_{T} \left(\bar{\ysf} - \ysf_d \right) \vsf  + \sum_{S \in \Sides } \int_{S } \llbracket \nabla \bar{\psf}_{\T} \cdot \nu \rrbracket  \vsf.
\end{equation*}
Since on each $T \in \T$ we have that $\vsf \in C^2(T)$, integration by parts also yields
\begin{equation*}
\label{eq:res2}
\int_{\Omega} \nabla \vsf \cdot \nabla (\bar{\psf} - \bar{\psf}_{\T}) = - \sum_{T \in \T} \int_{T} \Delta \vsf (\bar{\psf} - \bar{\psf}_{\T})   - \sum_{S \in \Sides } \int_{S} \llbracket \nabla \vsf \cdot \nu \rrbracket (\bar{\psf} - \bar{\psf}_{\T}).
\end{equation*}
In conclusion, since the left hand sides of the previous expressions coincide, we arrive at the identity
\begin{equation}\label{eq:error_adjoint}
\sum_{T \in \T} \int_{T} \left(\bar{\ysf} - \ysf_d \right) \vsf  + \sum_{S \in \Sides } \int_{S } \llbracket \nabla \bar{\psf}_{\T} \cdot \nu \rrbracket  \vsf
=
- \sum_{T \in \T} \int_{T} \Delta \vsf (\bar{\psf} - \bar{\psf}_{\T})   - \sum_{S \in \Sides } \int_{S} \llbracket \nabla \vsf \cdot \nu \rrbracket (\bar{\psf} - \bar{\psf}_{\T}),
\end{equation}
for every $\vsf \in H^1_0(\Omega)$ such that $\vsf_{|T} \in C^2(T)$ for all $T \in \T$. We now proceed, on the basis of \eqref{eq:defofEp}, in two steps.

\noindent \framebox{Step 1.} Let $T\in\T$. A simple application of the triangle inequality yields
\begin{equation}
\label{eq:aux_EO}
h_{T}^{2-n/2}
\|\bar{\ysf}_{\T}-\ysf_{d}\|_{L^{2}(T)}
\leq 
h_{T}^{2-n/2}
\|\bar{\ysf}_{\T}-\mathcal{P}_{\T}\ysf_{d}\|_{L^{2}(T)} +
 h_{T}^{2-n/2}
\|\ysf_{d}-\mathcal{P}_{\T}\ysf_{d}\|_{L^{2}(T)}.
\end{equation}
We recall that $\mathcal{P}_{\T}$ denotes the $L^2$-projection operator onto piecewise linear functions over $\T$. Now, by letting $\vsf=\beta_{T}=\left( \bar{\ysf}_{\T} - \mathcal{P}_{\T}\ysf_d  \right)\varphi_{T}^2$ in \eqref{eq:error_adjoint}, where $\varphi_{T}$ is the standard bubble function over $T$ \cite{Verfurth2,Verfurth}, we have that
\begin{align}
\int_{T}  \left( \bar{\ysf}_{\T} - \mathcal{P}_{\T}\ysf_d  \right) \beta_{T} 
& = 
\int_{T} \left[  (\bar{\ysf}-\ysf_{d})\beta_{T} - (\bar{\ysf}-\bar{\ysf}_{\T})\beta_{T}  + 
          \left( \ysf_d - \mathcal{P}_{\T}\ysf_d   \right) \beta_{T}\right] \nonumber \\
& =
-\int_{T} \Delta \beta_{T} (\bar{\psf} - \bar{\psf}_{\T}) - \int_{T}(\bar{\ysf}-\bar{\ysf}_{\T})\beta_{T} + \int_{T}\left( \ysf_d - \mathcal{P}_{\T}\ysf_d   \right) \beta_{T} := \textrm{I} + \textrm{II} + \textrm{III}.
\label{eq:eff_aux}
\end{align}
Notice that 
$
\int_{S } \llbracket \nabla \beta_{T} \cdot \nu \rrbracket (\bar \psf - \bar p_{\T}) = 0
$
for all $S\in\Sides$.
We now bound each term on the right hand side of \eqref{eq:eff_aux} separately. Since $\Delta(\bar{\ysf}_{\T}-\mathcal{P}_{\T}\ysf_{d})=0$, we have that
\[
\Delta \beta_{T} 
=4\nabla(\bar{\ysf}_{\T}-\mathcal{P}_{\T}\ysf_{d})\cdot\nabla \varphi_T \varphi_T
+
2 (\bar{\ysf}_{\T}-\mathcal{P}_{\T}\ysf_{d})(\varphi_T\Delta\varphi_T + \nabla\varphi_T\cdot\nabla \varphi_T).
\]
This, 
the properties of the bubble function $\varphi_T$ and an inverse inequality imply that
\begin{align*}
| \textrm{I} | & \lesssim 
\left(
h_{T}^{n/2-1}\|\nabla(\bar{\ysf}_{\T}-\mathcal{P}_{\T}\ysf_{d})\|_{L^{2}(T)} + 
h_{T}^{n/2-2}\|\bar{\ysf}_{\T}-\mathcal{P}_{\T}\ysf_{d}\|_{L^{2}(T)}
\right)\|\bar{\psf}-\bar{\psf}_{\T}\|_{L^{\infty}(T)}
\\
& \lesssim  h_{T}^{n/2-2}
\|\bar{\ysf}_{\T}-\mathcal{P}_{\T}\ysf_{d}\|_{L^{2}(T)} 
\|\bar{\psf}-\bar{\psf}_{\T}\|_{L^{\infty}(T)}.
\end{align*}
The terms \textrm{II} and \textrm{III} are bounded as follows:
\begin{align*}
 |\textrm{II}| 
\lesssim
\|\bar{\ysf}-\bar{\ysf}_{\T}\|_{L^{2}(T)}\|\bar{\ysf}_{\T} - \mathcal{P}_{\T}\ysf_d\|_{L^{2}(T)}
\end{align*}
and
\begin{align*}
|\textrm{III}|
\lesssim
\|\ysf_d - \mathcal{P}_{\T}\ysf_d\|_{L^{2}(T)}
\|\bar{\ysf}_{\T} - \mathcal{P}_{\T}\ysf_d\|_{L^{2}(T)}.
\end{align*}
In view of the fact that,
$
\|\bar{\ysf}_{\T} - \mathcal{P}_{\T}\ysf_d\|_{L^{2}(T)}^2  
\lesssim
\int_{T}  \left( \bar{\ysf}_{\T} - \mathcal{P}_{\T}\ysf_d  \right) \beta_{T} ,
$
the previous findings allow us to state that
\begin{equation*}
h_{T}^{2-n/2}
\| \bar{\ysf}_{\T} - \mathcal{P}_{\T}\ysf_d  \|_{L^{2}(T)}\\
\lesssim
\|\bar{\psf}-\bar{\psf}_{\T}\|_{L^{\infty}(T)} 
+ h_{T}^{2-n/2}
\left(
\|\bar{\ysf}-\bar{\ysf}_{\T}\|_{L^{2}(T)}
+
\|\ysf_d - \mathcal{P}_{\T}\ysf_d\|_{L^{2}(T)}
\right).
\end{equation*}
Consequently, \eqref{eq:aux_EO} allows us to conclude that
\begin{equation}\label{peffstep1}
h_{T}^{2-n/2} \|\bar{\ysf}_{\T} - \ysf_d   \|_{L^{2}(T)}
\lesssim
\|\bar{\psf}-\bar{\psf}_{\T}\|_{L^{\infty}(T)} +
 h_{T}^{2-n/2}
\|\bar{\ysf}-\bar{\ysf}_{\T}\|_{L^{2}(T)}
+
\textrm{osc}_{\T}(\ysf_d;T),
\end{equation}
where $\textrm{osc}_{\T}(\ysf_d;T)$ is defined as in \eqref{eq:osc_localf}.

\noindent \framebox{Step 2.} Let $T \in \T$ and $S \in \Sides_T$. The goal of this step is control the term $h_T \| \llbracket   \nabla \bar{\psf}_{\T} \cdot \nu \rrbracket  \|_{L^{\infty}(S)}$ in \eqref{eq:defofEp}. To do this, we use the property 
\[
|S|\| \llbracket \nabla \bar{\psf}_{\T} \cdot \nu \rrbracket \|_{L^{\infty}(S)} \lesssim \left|\int_S \llbracket \nabla \bar{\psf}_{\T} \cdot \nu \rrbracket \varphi_S \right|,
\]
of $\varphi_S$, the standard bubble function over $S$ \cite{Verfurth2,Verfurth}. We control the right hand side of the previous expression by letting $\vsf=\varphi_S$ in \eqref{eq:error_adjoint}. This yields
\begin{align*}
  \left| \int_S \llbracket  \nabla \bar{\psf}_\T \cdot \nu \rrbracket \varphi_S \right|
   \leq & \sum_{T' \in \Ne_S} \int_{T'} |\bar\ysf - \ysf_d| \varphi_S +
  \sum_{T' \in \Ne_S} \int_{T'} |\bar{\psf}-\bar{\psf}_\T| |\LAP \varphi_S|
  + \sum_{T' \in \Ne_S}\sum_{S' \in \Sides_{T'}}
   \int_{S'} |\bar{\psf}-\bar{\psf}_\T| |\llbracket \nabla \varphi_S \cdot \nu \rrbracket|
   \\
   \lesssim & \sum_{T' \in \Ne_S}|T'|^{1/2}\left( \|\bar \ysf - \bar \ysf_{\T} \|_{L^2(T')}+ \| \bar \ysf_{\T} - \ysf_d\|_{L^2(T')}\right) \\
  &+ \sum_{T' \in \Ne_S}\left(h_S^{-2} |T'| + h_S^{-1} \sum_{S' \in \Sides_{T'}}|S'|\right) \| \bar{\psf}-\bar{\psf}_\T \|_{L^\infty(T')}.
\end{align*}
Combining this estimate with \eqref{peffstep1} yields the bound
\begin{align*}
 h_T \| \llbracket  \nabla \bar{\psf}_{\T} \cdot \nu \rrbracket  \|_{L^{\infty}(S)} &\lesssim h_T^{2-n/2} \|\bar \ysf - \bar \ysf_{\T} \|_{L^{2}(\Ne_S)} +  \| \bar \psf -  \bar \psf_{\T} \|_{L^{\infty}(\Ne_S)} + \textrm{osc}_{\T}( \ysf_d;\Ne_S).
\end{align*}

We finally combine the results of Step 1 and 2 and arrive at the desired estimate \eqref{eq:efficiencyp}.
\end{proof}

The results of Lemmas \ref{le:local_eff_y} and \ref{le:local_eff_p} immediately yield the following result.

\begin{thrm}[local efficiency of $\E_{\textrm{ocp}}$] 
Let $( \bar{ \mathbf{u} },\bar{\ysf},\bar{\psf} ) \in \mathcal{U}_{\textrm{ad}} \times H_0^1(\rho,\Omega) \times H_0^1(\Omega)$ be the solution to the optimality system \eqref{eq:optimal_system} associated with the optimal control with point sources and $(\bar{\mathbf{u}}_{\T},\bar{\ysf}_{\T},\bar{\psf}_{\T}) \in \mathcal{U}_{\textrm{ad}}\times \V(\T) \times \V(\T)$ be its numerical approximation given by \eqref{eq:discrete_state_equation}--\eqref{eq:discrete_adjoint_equation}. If $\alpha \in (n-2,2)$, then 
\begin{align*}
\E_{\ysf}^2 (\bar{\ysf}_{\T},\bar{\mathbf{u}}_{\T}; T) + \E_{\psf}^2 (\bar{\psf}_{\T},\bar{\ysf}_{\T}; T)
\lesssim &
\| \nabla( \bar{\ysf} - \bar{\ysf}_{\T}) \|_{L^{2}(\rho,\Ne_T)}^2 + \|  \bar{\psf} - \bar{\psf}_{\T} \|_{L^{\infty}(\Ne_T^*)}^2
\\
& +h_{T}^{\alpha+2-n}\sum_{z\in T \cap D}|\bar{\usf}_{z} - \bar{\usf}_{\T,z}|^2
+ h_T^{4 - n}\| \bar{\ysf} - \bar{\ysf}_{\T} \|_{L^{2}(\Ne_T^*)}^2 + \mathrm{osc}_{\T}^2( \ysf_d;\Ne_T^*),
\end{align*}
where $\Ne_T$ and $\Ne_T^{*}$ are given by \eqref{eq:NeT} and \eqref{eq:NeTstar}, respectively and the hidden constant is independent of the optimal variables, their approximations, the size of the elements in the mesh $\T$ and $\#\T$.
\label{th:localeff_ocp}
\end{thrm}

Our final result gives the global efficiency property of the estimator.

\begin{thrm}[global efficiency of $\E_{\textrm{ocp}}$] 
Let $( \bar{ \mathbf{u} },\bar{\ysf},\bar{\psf} ) \in \mathcal{U}_{\textrm{ad}} \times H_0^1(\rho,\Omega) \times H_0^1(\Omega)$ be the solution to the optimality system \eqref{eq:optimal_system} associated with the optimal control with point sources and $(\bar{\mathbf{u}}_{\T},\bar{\ysf}_{\T},\bar{\psf}_{\T}) \in \mathcal{U}_{\textrm{ad}}\times \V(\T) \times \V(\T)$ be its numerical approximation given by \eqref{eq:discrete_state_equation}--\eqref{eq:discrete_adjoint_equation}. If $\alpha \in (n-2,2)$, then 
\begin{equation}\label{eq:efficiencyocp}
\E_{\mathrm{ocp}}^2 (\bar{\ysf}_{\T},\bar{\psf}_{\T},\bar{\mathbf{u}}_{\T}; \T) \lesssim 
\|  \bar{\mathbf{u}} - \bar{\mathbf{u}}_{\T} \|_{\mathbb{R}^l}^2+\| \nabla( \bar{\ysf} - \bar{\ysf}_{\T}) \|_{L^{2}(\rho,\Omega)}^2 +\|  \bar{\psf} - \bar{\psf}_{\T} \|_{L^{\infty}(\Omega)}^2 
 + \max_{T\in\T}\mathrm{osc}_{\T}^2( \ysf_d;\Ne_T^*)
\end{equation}
where $\Ne_T^{*}$ is defined as in \eqref{eq:NeTstar} and the hidden constant is independent of the optimal variables, their approximations, the size of the elements in the mesh $\T$ and $\#\T$.
\label{th:eff_ocp}
\end{thrm}
\begin{proof}
Since assumption \eqref{eq:patch} implies that, for all $T\in\T$, $\#(T \cap D)\leq 1$, we arrive at
\begin{align*}
\sum_{T\in\T}\sum_{z\in T \cap D}h_T^{\alpha + 2 - n}|\bar{\usf}_{z} - \bar{\usf}_{\T,z}|^2\leq&\mathrm{diam}(\Omega)^{\alpha + 2 - n}\sum_{T\in\T}\sum_{z\in T \cap D}|\bar{\usf}_{z} - \bar{\usf}_{\T,z}|^2
\\
\leq&\mathrm{diam}(\Omega)^{\alpha + 2 - n}\left(\max_{z\in D}\#\T_z\right)\sum_{z\in D}|\bar{\usf}_{z} - \bar{\usf}_{\T,z}|^2
\\
=&\mathrm{diam}(\Omega)^{\alpha + 2 - n}\left( \max_{z\in D}\#\T_z \right)\|  \bar{\mathbf{u}} - \bar{\mathbf{u}}_{\T} \|^2_{\mathbb{R}^l}
\end{align*}
where $\T_z=\{T\in\T:\,z\in T\}$. In view of the definition of the estimator $\E_{\ysf}$ given by \eqref{eq:defofEyglobal}, the estimate \eqref{eq:efficiencyy} and the the finite overlapping property of stars, we can then conclude that
\begin{equation}
\label{eq:global_ocp_aux}
\E_{\ysf}^2 (\bar{\ysf}_{\T},\bar{\mathbf{u}}_{\T}; \T) \lesssim \| \nabla( \bar{\ysf} - \bar{\ysf}_{\T}) \|_{L^{2}(\rho,\Omega)}^2 
+ \mathrm{diam}(\Omega)^{\alpha + 2 - n} \|  \bar{\mathbf{u}} - \bar{\mathbf{u}}_{\T} \|_{\mathbb{R}^l}^2.
\end{equation}
On the other hand, the definition of the estimator $\E_{\psf}$ given by \eqref{eq:defofEpglobal} and the estimate \eqref{eq:efficiencyp} provide the bound
\begin{equation*}
\E_{\psf}^2 (\bar{\psf}_{\T},\bar{\ysf}_{\T}; \T)  \lesssim \|  \bar{\psf} - \bar{\psf}_{\T} \|_{L^{\infty}(\Omega)}^2 + \mathrm{diam}(\Omega)^{4-n} \| \bar{\ysf} - \bar{\ysf}_\T \|_{L^2(\Omega)}^2 
+ \max_{T \in \T} \mathrm{osc}_{\T}^2(\ysf_d;\Ne_T^*).
\end{equation*}
By combining this estimate with the weighted Poincar\'e inequality 
$\| \bar{\ysf} - \bar{\ysf}_\T \|_{L^2(\Omega)}\lesssim\| \nabla( \bar{\ysf} - \bar{\ysf}_{\T}) \|_{L^{2}(\rho,\Omega)}$
of Lemma \ref{lem:restalphagen}, that holds for $\alpha \in (n-2,2)$, and \eqref{eq:global_ocp_aux}, we arrive at \eqref{eq:efficiencyocp}. This concludes the proof.
\end{proof}

\section{Numerical examples}
\label{sec:numex}
We conduct a series of numerical examples that illustrate the performance of the error estimator. In some of these examples, we go beyond the presented theory and perform numerical experiments where we violate the assumption of homogeneous Dirichlet boundary conditions. These have been carried out with the help of a code that we implemented using \texttt{C++}. All matrices have been assembled exactly. The right hand sides and approximation errors are computed by a quadrature formula which is exact for polynomials of degree 19 for two dimensional domains and degree 14 for three dimensional domains. All linear systems were solved using the multifrontal massively parallel sparse direct solver (MUMPS) \cite{MUMPS1,NUMPS2}. 

\begin{figure}[!htbp]
\begin{center}
\scalebox{0.2}{\includegraphics{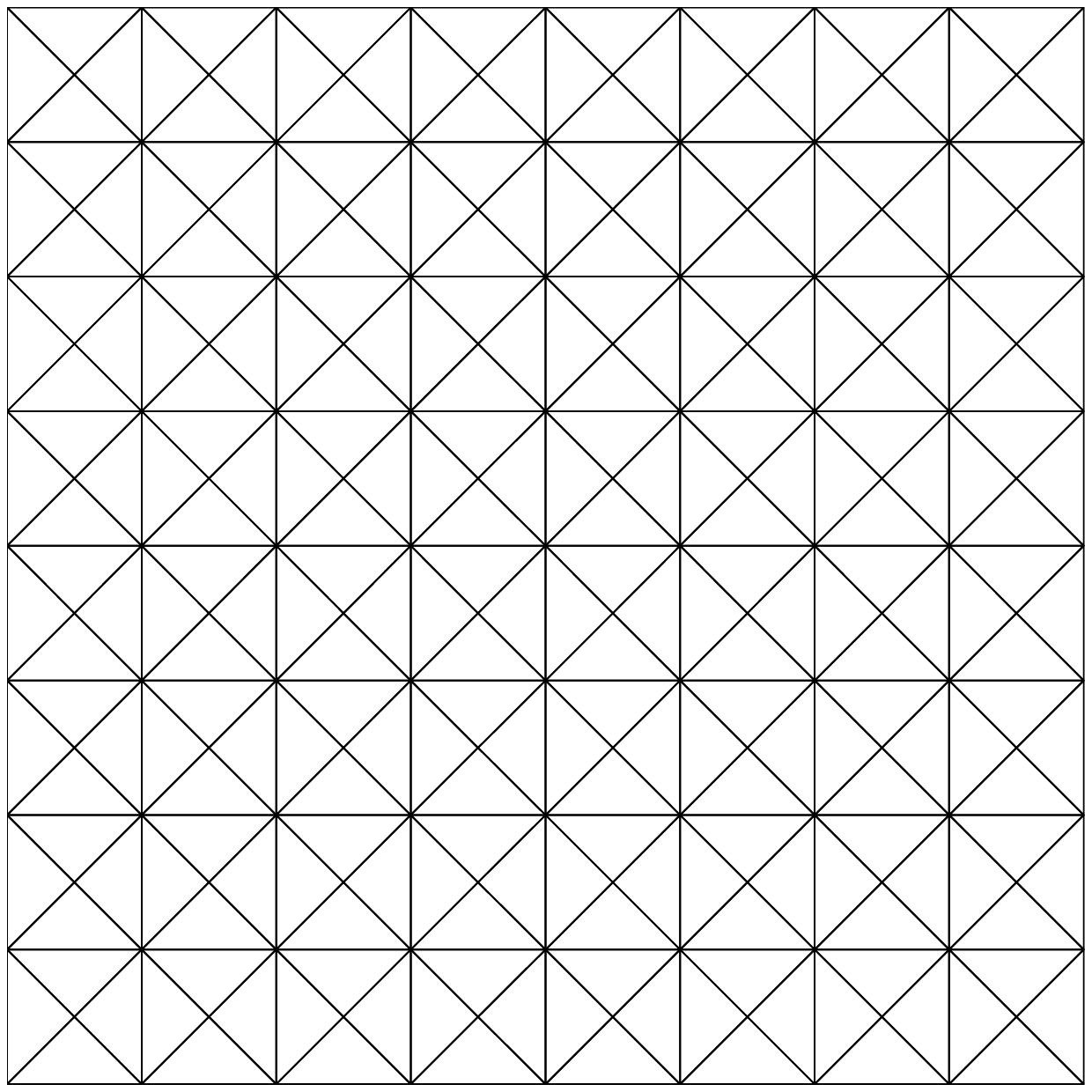}}
\scalebox{0.2}{\includegraphics{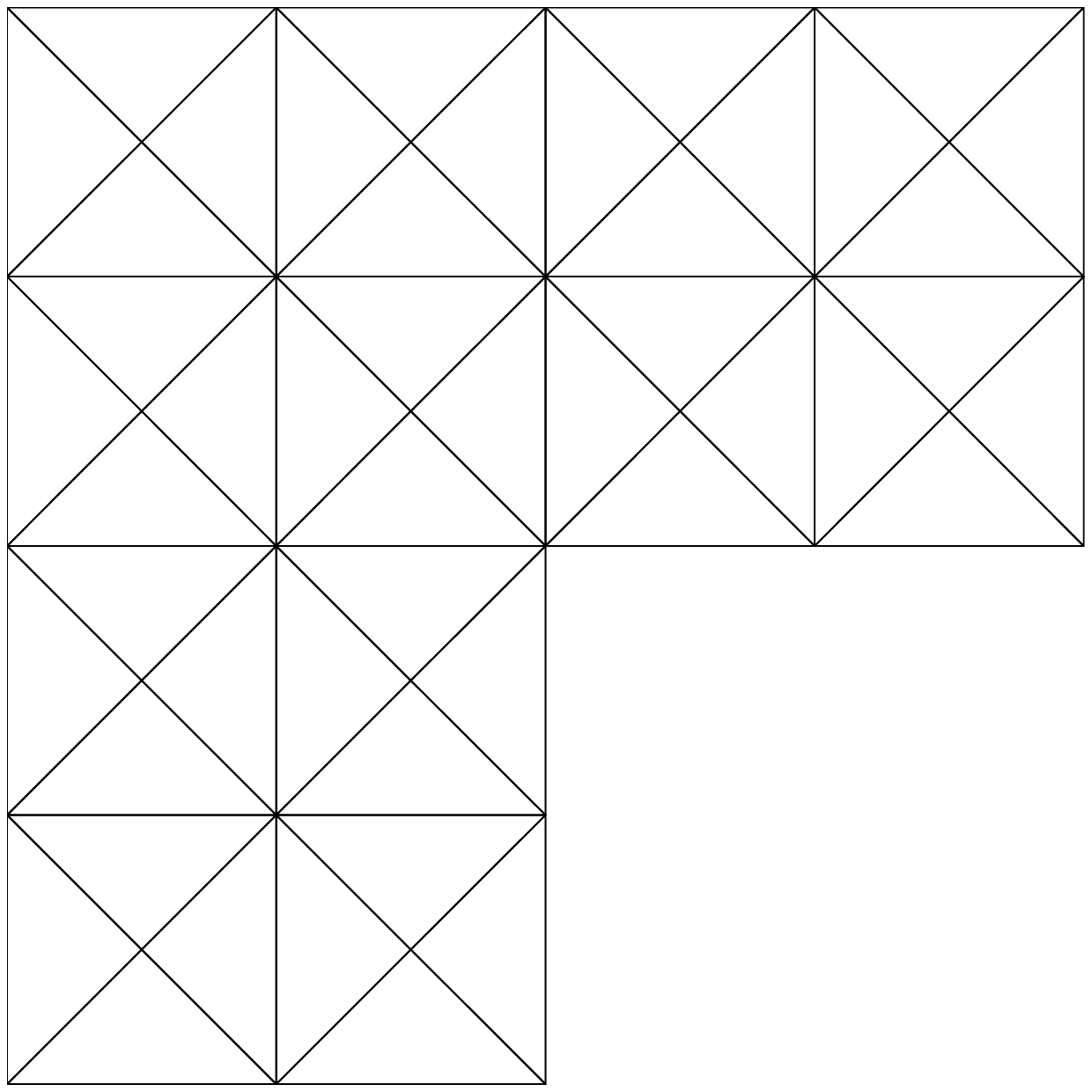}}
\scalebox{0.2}{\includegraphics{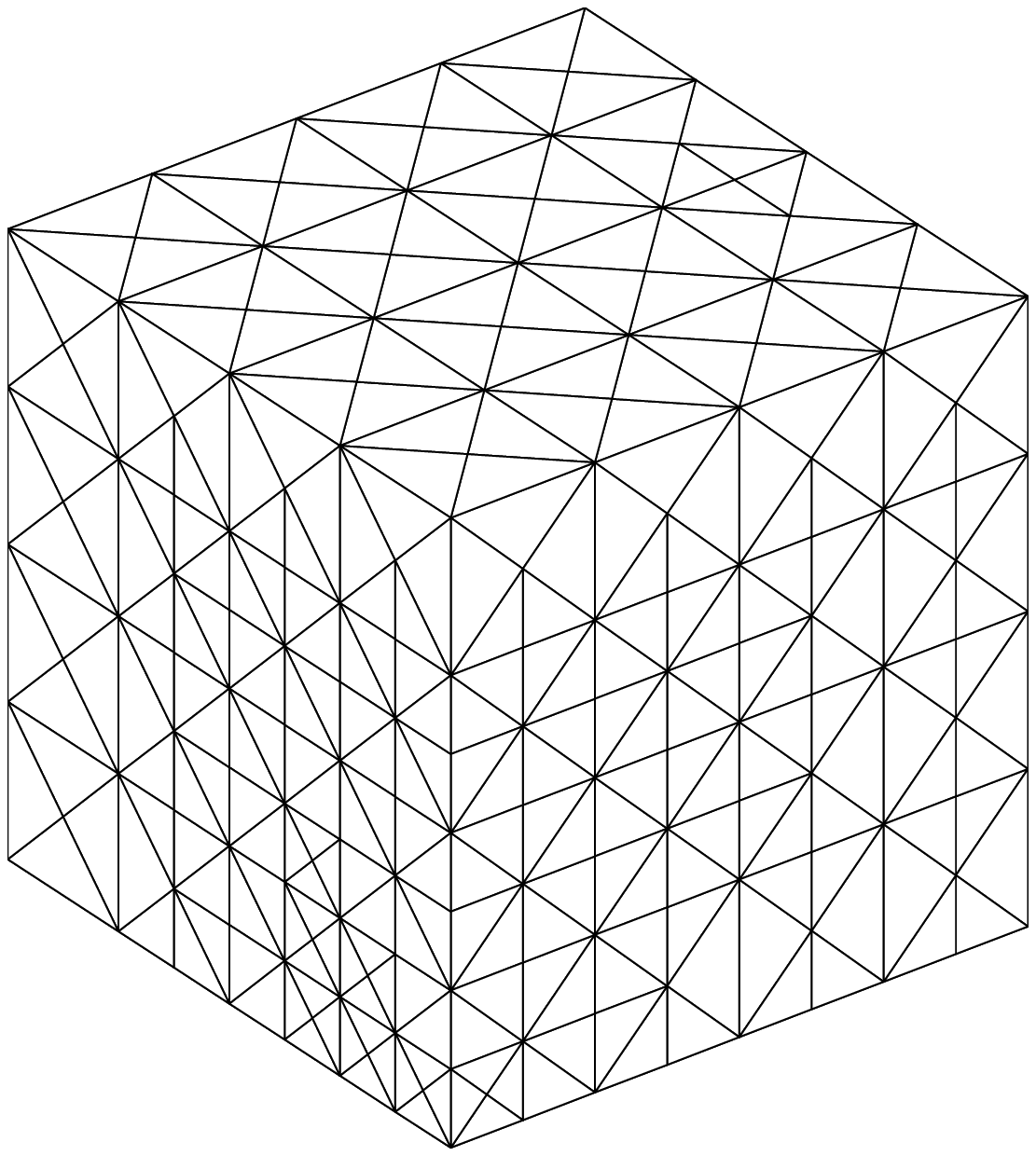}}
\scalebox{0.2}{\includegraphics{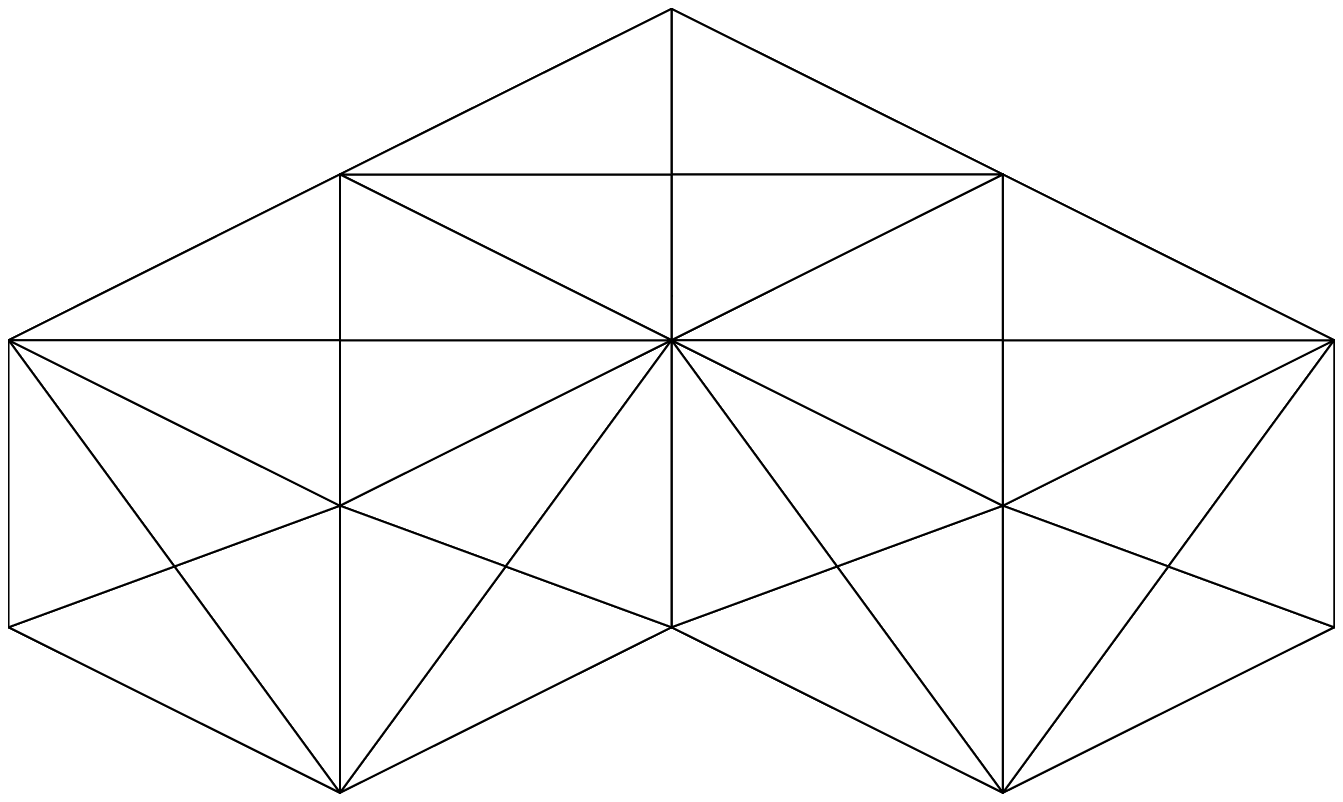}}
\end{center}
\caption{The initial meshes used when the domain $\Omega$ is a square (Examples 1 and 2), two dimensional L-shape (Example 3), cube (Examples 4 and 5) and three dimensional L-shape (Example 6).}
\label{FigMesh}
\end{figure}

For a given partition $\T$ we seek $(\bar{\ysf}_\T,\bar{\psf}_\T,\bar{\usf}_\T)\in \V(\T)\times \V(\T)\times \R^{l}$ that solves \eqref{eq:discrete_state_equation}, \eqref{eq:first_oder_discrete} and \eqref{eq:discrete_adjoint_equation}. We solve the nonlinear system of equations using a Newton-type primal-dual active set strategy \cite[\S 2.12.4]{Tbook}. Once a discrete solution is obtained, we calculate the error estimator and use the local error indicators to drive the adaptive procedure described in \textbf{Algorithm 1}. On the basis of \eqref{eq:reliability} we use the a posteriori error estimator defined by \eqref{eq:defofEocp} with corresponding local error indicators
\begin{equation}
\E_{\mathsf{ocp};T}=
\E_{\ysf}^{2}(\bar{\ysf}_{\T},\bar{\usf}_{\T};T)+
\E_{\psf}^{2}(\bar{\psf}_{\T},\bar{\ysf}_{\T};T),
\end{equation}
which are defined in terms of \eqref{eq:defofEy} and \eqref{eq:defofEp}. The total number of degrees of freedom $\textrm{Ndof} = 2 \dim(\V(\T)) + l$, where $l = \# D$. The initial meshes for our numerical examples are shown in Figure~\ref{FigMesh}.

We consider problems with homogeneous Dirichlet boundary conditions whose exact solutions are not known, and problems with inhomogeneous Dirichlet boundary conditions whoses exact solutions are known. For the numerical examples for which the exact solutions are known, we took the optimal state variable to be a linear combination of fundamental solutions for the Laplacian, that is
\begin{equation}
\label{exact_adjoint}
\bar{\ysf}(x) = \begin{dcases}
                  -\frac{1}{2\pi}\sum_{z\in D}\varrho_{z}\log|x-z|, & \textrm{if}~\Omega\subset\mathbb{R}^{2},\\
                  \frac{1}{4\pi}\sum_{z\in D}\varrho_{z}\frac{1}{|x-z|}, & \textrm{if}~\Omega\subset\mathbb{R}^{3},
                \end{dcases}
\end{equation}
with $\varrho_{z}\in\mathbb{R}$ for all $z\in D$. Upon fixing an exact adjoint state and the constraints $\asf$ and $\bsf$, the exact optimal control is computed using the projection formula \eqref{eq:Pi}. The $\varrho_{z}$ in \eqref{exact_adjoint} are then computed using \eqref{eq:defofPDE}. Finally, the desired state is computed using \eqref{eq:adjp}. The error is measured in the norm
\[
  \|(e_{\bar{\ysf}},e_{\bar{\psf}},e_{\bar{\usf}})\|_{\Omega}^{2} = 
  \|\nabla e_{\bar{\ysf}}\|_{L^{2}(\rho,\Omega)}^{2}+
  \|e_{\bar{\psf}}\|_{L^{\infty}(\Omega)}^{2}+
  \|e_{\bar{\usf}}\|_{\R^{l}}^{2},
\]
where $e_{\bar{\ysf}}=\bar{\ysf}-\bar{\ysf}_{\T}$, $e_{\bar{\psf}}=\bar{\psf}-\bar{\psf}_{\T}$ and $e_{\bar{\usf}}=\bar{\usf}-\bar{\usf}_{\T}$.

\begin{table}[!htbp]
\begin{flushleft}
\begin{tabular}{l l} 
\multicolumn{2}{l}{\textbf{Algorithm 1:  Adaptive Primal-Dual Active Set Algorithm.}} 
\vspace{0.15cm}\\
\toprule
\multicolumn{2}{l}{\textbf{Input:} Initial mesh $\T_{0}$, set of source points $D$, desired state $\ysf_{d}$, constraints $\asf$ and $\bsf$,}\\
\multicolumn{2}{l}{and regularization parameter $\lambda$.}
\vspace{0.1cm} \\
\textbf{Set:}  & $i=0$.
\vspace{0.1cm}\\
\multicolumn{2}{l}{\textbf{Active set strategy:}}
\vspace{0.1cm}\\
\textbf{1:}    &  Compute $[\bar{\ysf}_{\T},\bar{\psf}_{\T},\bar{\usf}_{\T}]=\textrm{\textbf{Active-Set}}[\T_{i},D,\ysf_{d},\asf,\bsf,\lambda]$. \\
               &  $\textrm{\textbf{Active-Set}}$ implements the active set strategy of \cite[\S 2.12.4]{Tbook}. \\
\multicolumn{2}{l}{\textbf{Adaptive loop:}}
\vspace{0.1cm}\\
\textbf{2:}    &  For each $T \in \T$ compute the local error indicator $\E_{\mathsf{ocp};T}$.
 \\
\textbf{3:}    & Mark an element $T$ for refinement if $\displaystyle\E_{\mathsf{ocp};T}^{2}> 0.5\max_{T'\in\T}\E_{\mathsf{ocp};T'}^{2}$.\\
\textbf{4:}    & From step \textbf{3}, construct a new mesh, using a longest edge bisection algorithm. \\& Set $i\leftarrow i+1$, and go to step \textbf{1}.\\
\bottomrule
\end{tabular}
\vspace{-0.3cm}
\end{flushleft}
\end{table}

\subsection{Two-dimensional examples}

\begin{figure}[!htbp]
\begin{center}
\psfrag{- lambda=1.0}{\LARGE $\lambda=1.0$}
\psfrag{- lambda=0.1}{\LARGE $\lambda=0.1$}
\psfrag{- lambda=0.01}{\LARGE $\lambda=0.01$}
\psfrag{- lambda=0.001}{\LARGE $\lambda=0.001$}
\psfrag{- lambda=0.0001}{\LARGE $\lambda=0.0001$}
\psfrag{Ndofs}{\huge Ndof}
\psfrag{O(1)}{\footnotesize {\LARGE $\textrm{Ndof}^{-1/2}$}}
\psfrag{O(2)}{\footnotesize{\LARGE $\textrm{Ndof}^{-1}$}}
\psfrag{- total error estimator}{\LARGE $\E_{\mathsf{ocp};\T}$}
\psfrag{- error estimator yh}{\LARGE $\E_{\ysf}$}
\psfrag{- error estimator ph}{\LARGE $\E_{\psf}$}
\psfrag{example A - total-estimator}{\huge  $\E_{\mathsf{ocp};\T}$ varying $\lambda$}
\psfrag{example A - lambda=----}{\huge {Example 1 - $\lambda=1$}}
\subfigure[]{\scalebox{.4}{\includegraphics{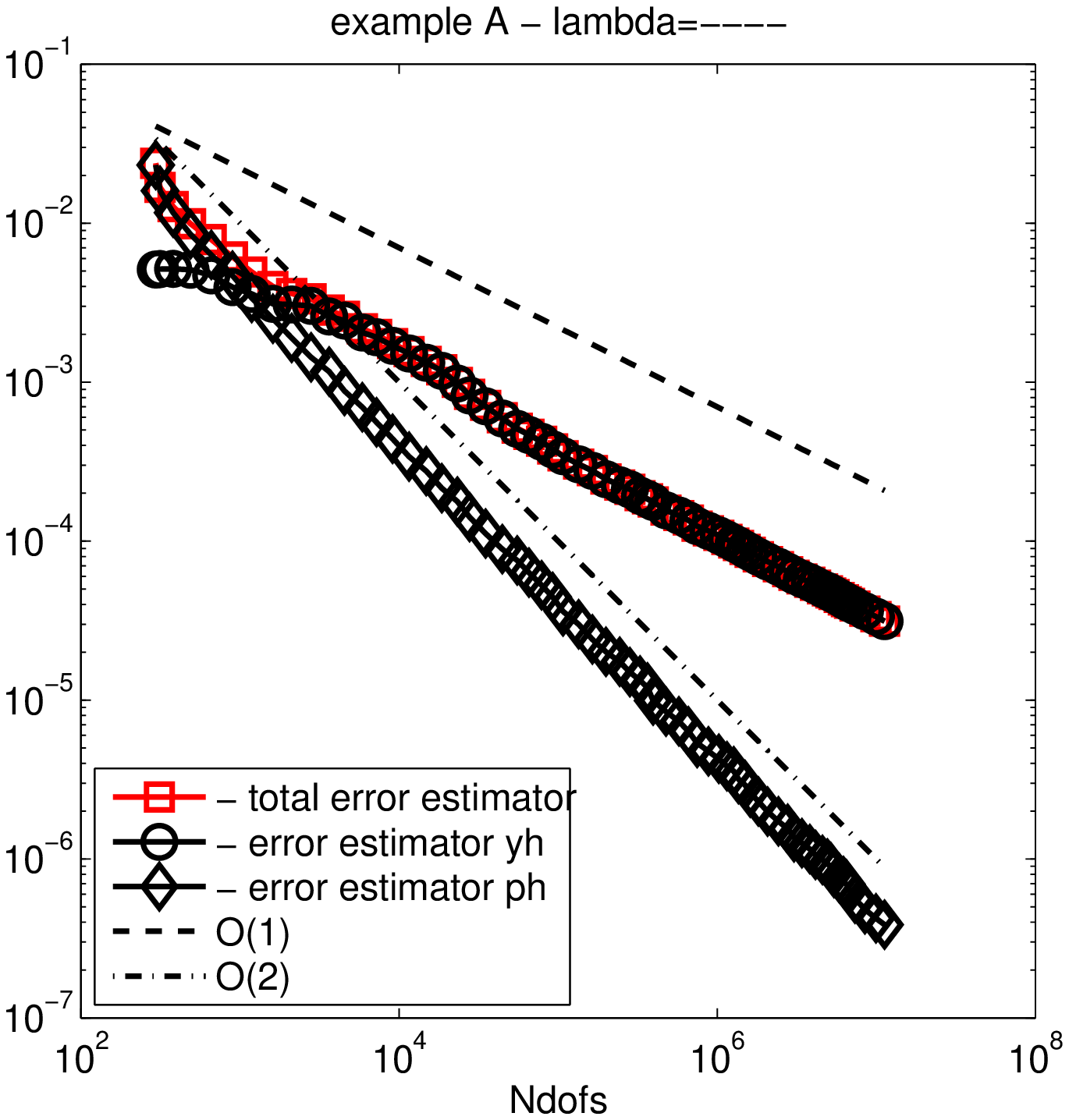}}}
\psfrag{example A - lambda=----}{\huge {Example 1 - $\lambda=10^{-1}$}}
\subfigure[]{\scalebox{.4}{\includegraphics{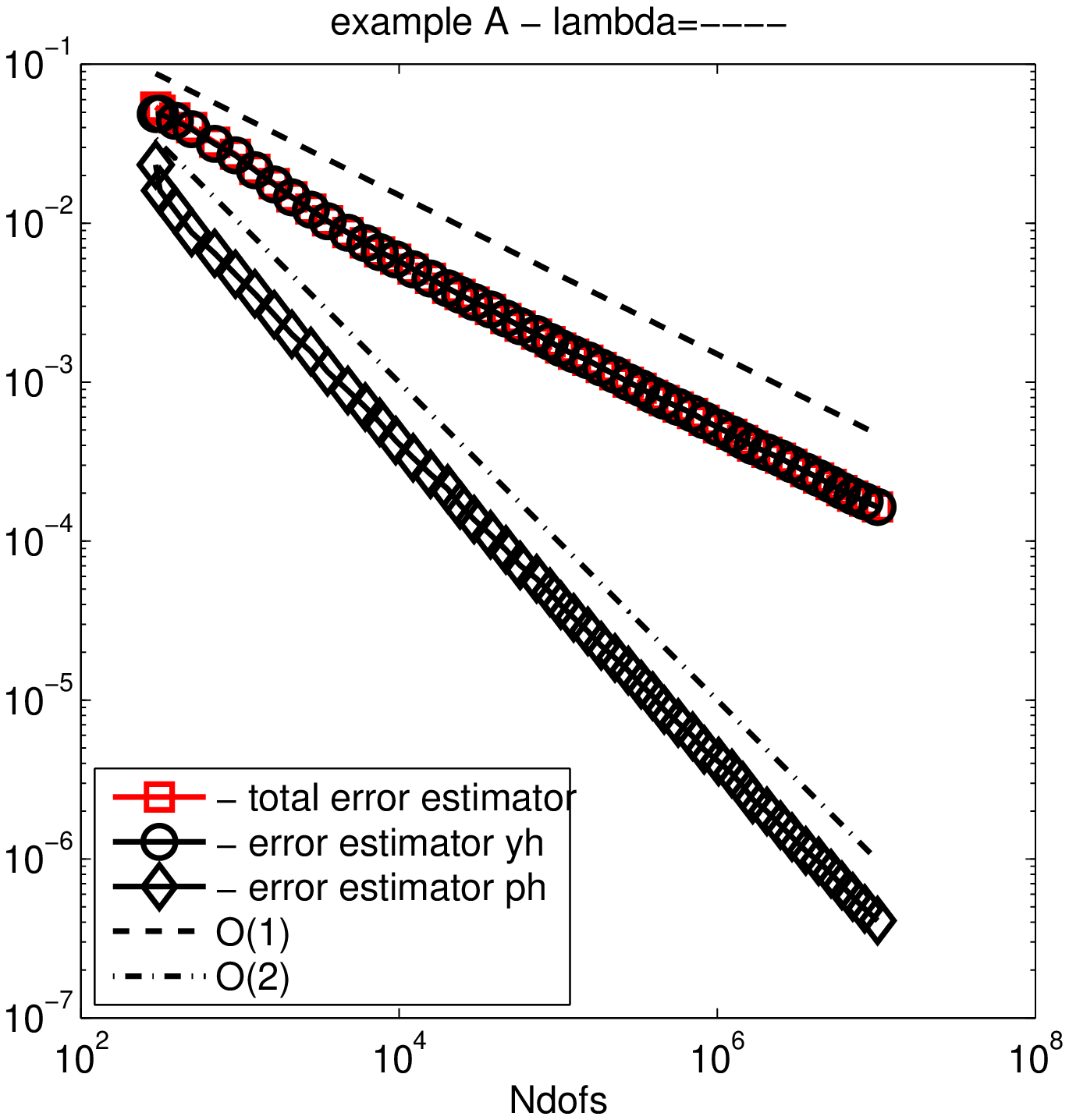}}}
\psfrag{example A - lambda=----}{\huge {Example 1 - $\lambda=10^{-2}$}}
\subfigure[]{\scalebox{.4}{\includegraphics{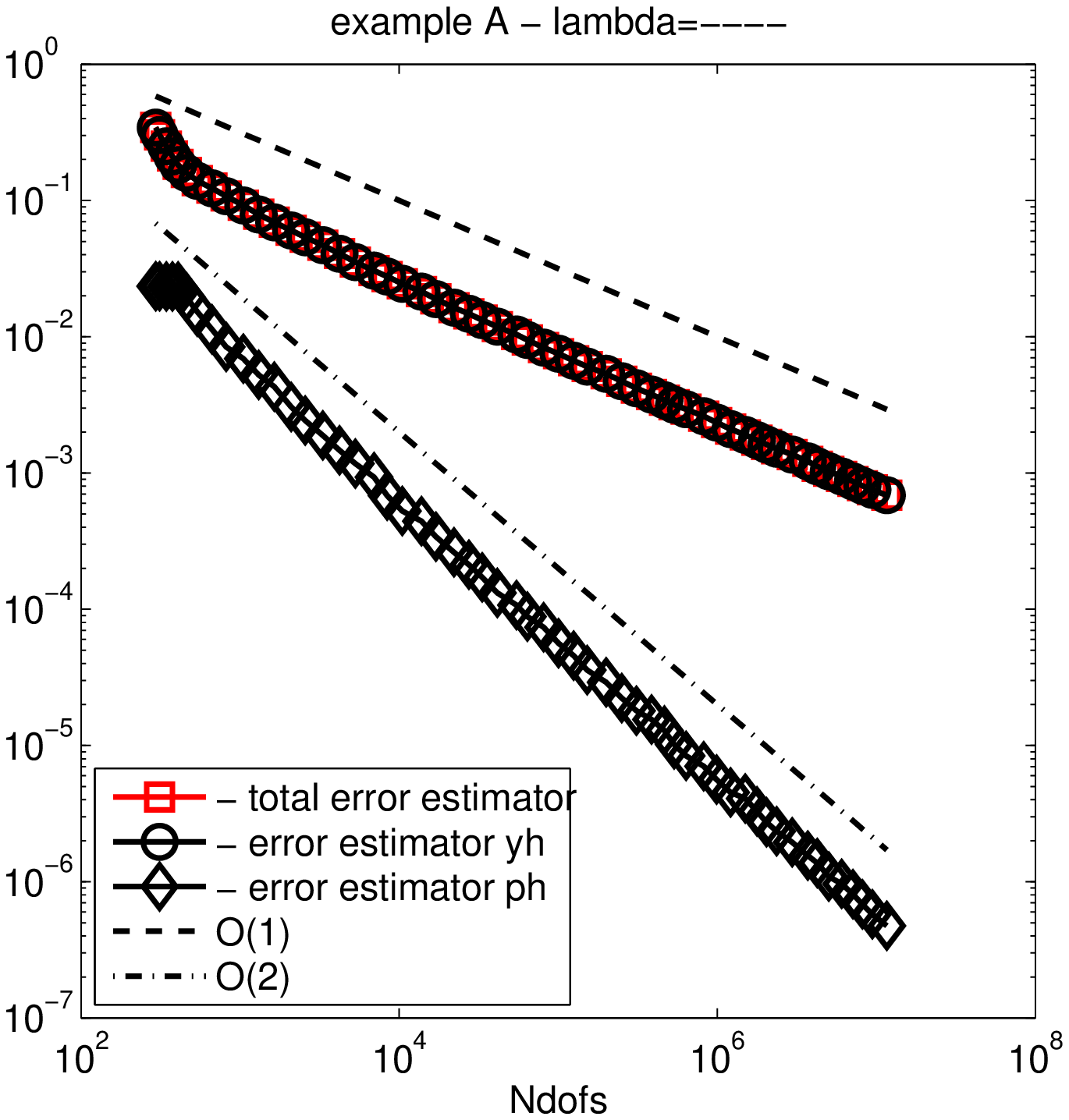}}}
\psfrag{example A - lambda=----}{\huge {Example 1 - $\lambda=10^{-3}$}}
\subfigure[]{\scalebox{.4}{\includegraphics{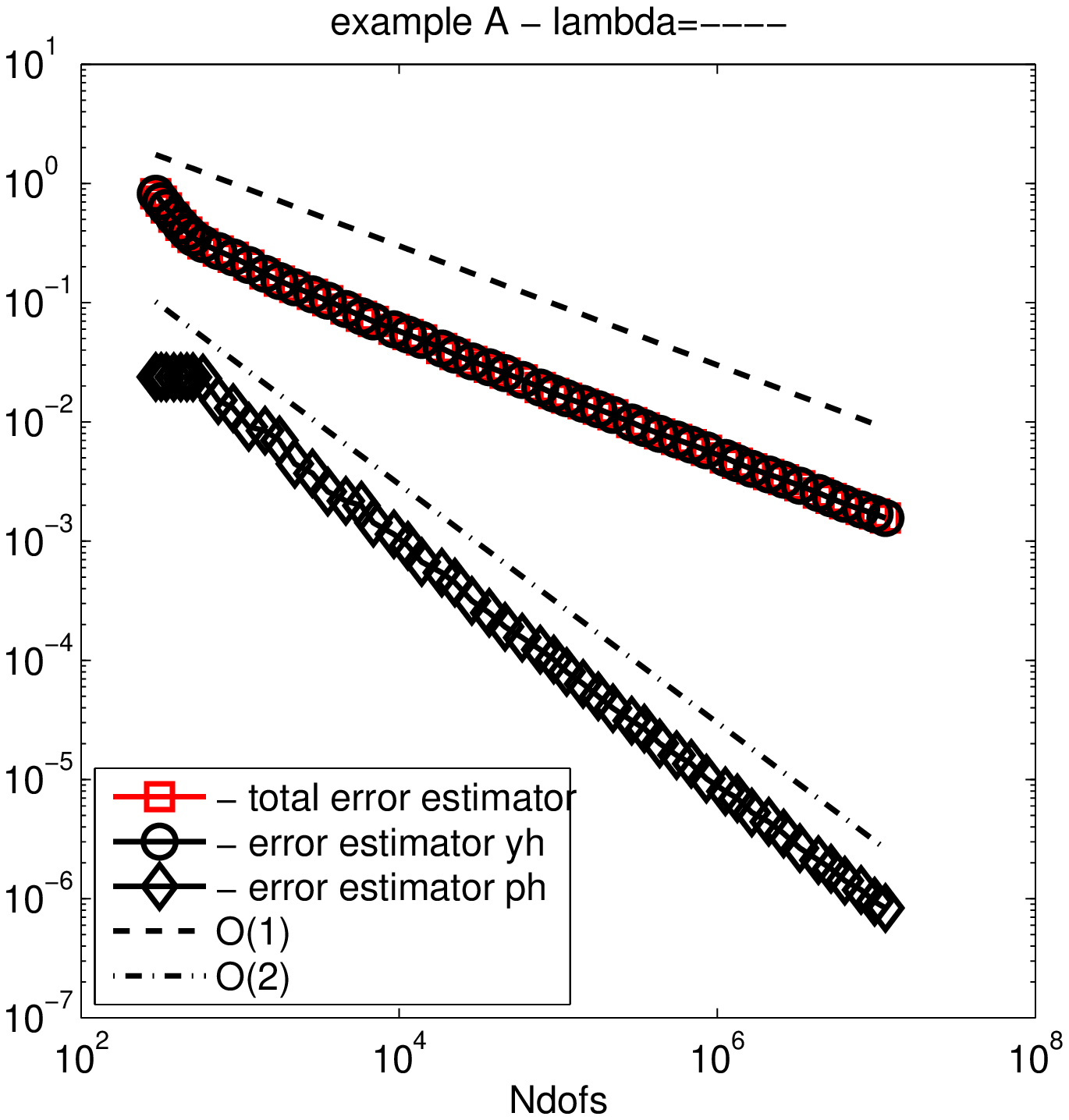}}}
\psfrag{example A - lambda=----}{\huge {Example 1 - $\lambda=10^{-4}$}}
\subfigure[]{\scalebox{.4}{\includegraphics{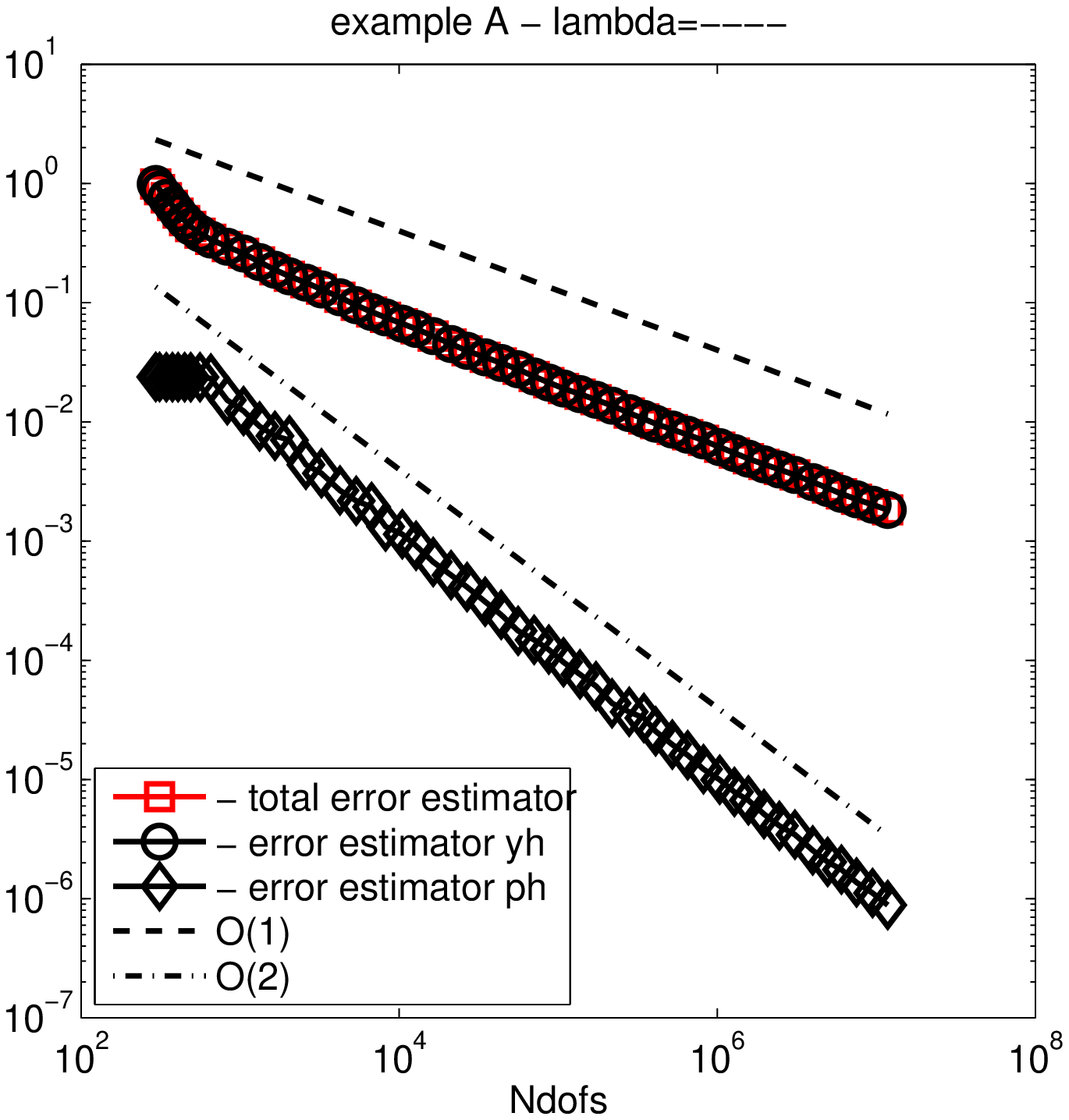}}}
\psfrag{example A - total}{\huge {$\E_{\mathsf{ocp};\T}$ varying $\lambda$}}
\subfigure[]{\scalebox{.4}{\includegraphics{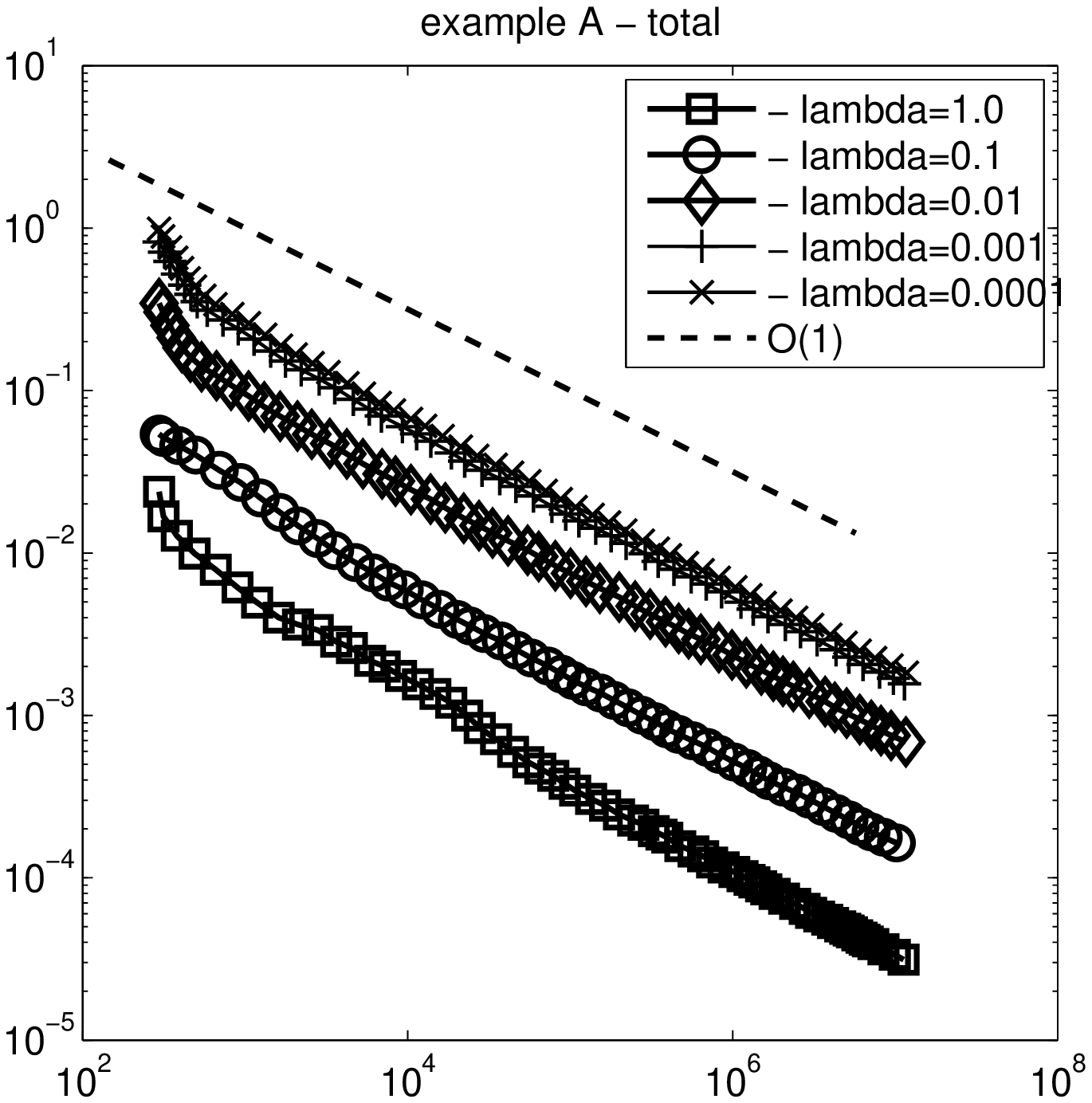}}}
\end{center}
\caption{Example 1: Experimental rates of convergence for the error estimator $\E_{\mathsf{ocp};\T}$ and its contributions $\E_{\ysf}$ and $\E_{\psf}$ for $\lambda\in\{1,10^{-1},10^{-2},10^{-3},10^{-4}\}$ (a)--(e); comparison of the error estimator $\E_{\mathsf{ocp};\T}$ for $\lambda\in\{1,10^{-1},10^{-2},10^{-3},10^{-4}\}$ (f).}
\label{Fig:1}
\end{figure}

\begin{figure}[!htbp]
\psfrag{example A - total-estimator}{\Huge $\E_{\mathsf{ocp};\T}$ varying $\lambda$}
\psfrag{- lambda=1.0}{\Huge $\lambda=1.0$}
\psfrag{- lambda=0.1}{\Huge $\lambda=0.1$}
\psfrag{- lambda=0.01}{\Huge $\lambda=0.01$}
\psfrag{- lambda=0.001}{\Huge $\lambda=0.001$}
\psfrag{Ndofs}{\Huge Ndof}
\psfrag{O(1)}{\huge $\textrm{Ndof}^{-1/2}$}
\psfrag{O(2)}{\huge $\textrm{Ndof}^{-1}$}
\psfrag{- total error estimator}{\Huge $\E_{\mathsf{ocp};\T}$}
\psfrag{- error estimator yh}{\Huge $\E_{\ysf}$}
\psfrag{- error estimator ph}{\Huge $\E_{\psf}$}
\psfrag{example A - total error norm}{\Huge{Example 2 - $\|(e_{\bar{\ysf}},e_{\bar{\psf}},e_{\bar{\usf}})\|_{\Omega}$}}
\psfrag{- error norm yh}{\Huge $\|\nabla(\bar{\ysf}-\bar{\ysf}_{\T})\|_{L^{2}(\rho,\Omega)}$}
\psfrag{- error norm ph}{\Huge $\|\bar{\psf}-\bar{\psf}_{\T}\|_{L^{\infty}(\Omega)}$}
\psfrag{- error norm uh}{\Huge $\|\bar{\usf}-\bar{\usf}_{\T}\|_{\mathbb{R}^{l}}$}
\psfrag{- error estimator yh}{\Huge $\E_{\ysf}$}
\psfrag{- error estimator ph}{\Huge $\E_{\psf}$}
\psfrag{- total error norm -----}{\Huge $\|(e_{\bar{\ysf}},e_{\bar{\psf}},e_{\bar{\usf}})\|_{\Omega}$}
\psfrag{- total error estimator}{\Huge $\E_{\mathsf{ocp};\T}$}
\psfrag{example A - total error estimator}{\Huge{Example 2 - $\E_{\mathsf{ocp};\T}$}}
\psfrag{example A - D=D algo - alpha=----}{}
\psfrag{- total error norm-1}{\Huge $\alpha=0.1$}
\psfrag{- total error norm-2}{\Huge $\alpha=0.5$}
\psfrag{- total error norm-3}{\Huge $\alpha=1$}
\psfrag{- total error norm-4}{\Huge $\alpha=1.5$}
\psfrag{- total error norm-5}{\Huge $\alpha=1.9$}
\psfrag{- total error estimator-1}{\Huge $\alpha=0.1$}
\psfrag{- total error estimator-2}{\Huge $\alpha=0.5$}
\psfrag{- total error estimator-3}{\Huge $\alpha=1$}
\psfrag{- total error estimator-4}{\Huge $\alpha=1.5$}
\psfrag{- total error estimator-5}{\Huge $\alpha=1.9$}
\begin{center}
\psfrag{example A - D=D algo - alpha=----}{\Huge{Example 2 - $\alpha=0.1$}}
\subfigure[]{\scalebox{.23}{\includegraphics{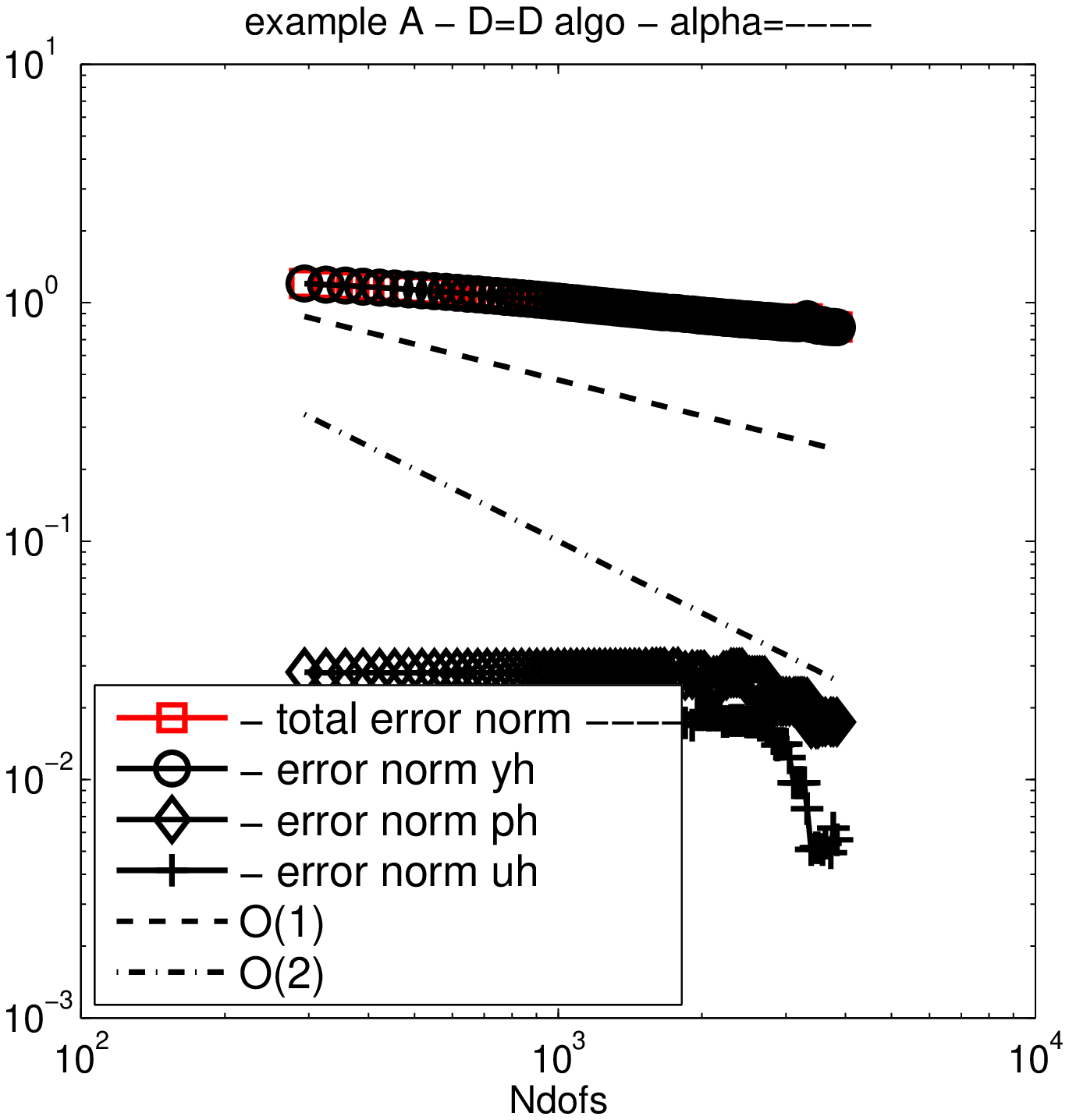}}}
\psfrag{example A - D=D algo - alpha=----}{\Huge{Example 2 - $\alpha=0.5$}}
\subfigure[]{\scalebox{.23}{\includegraphics{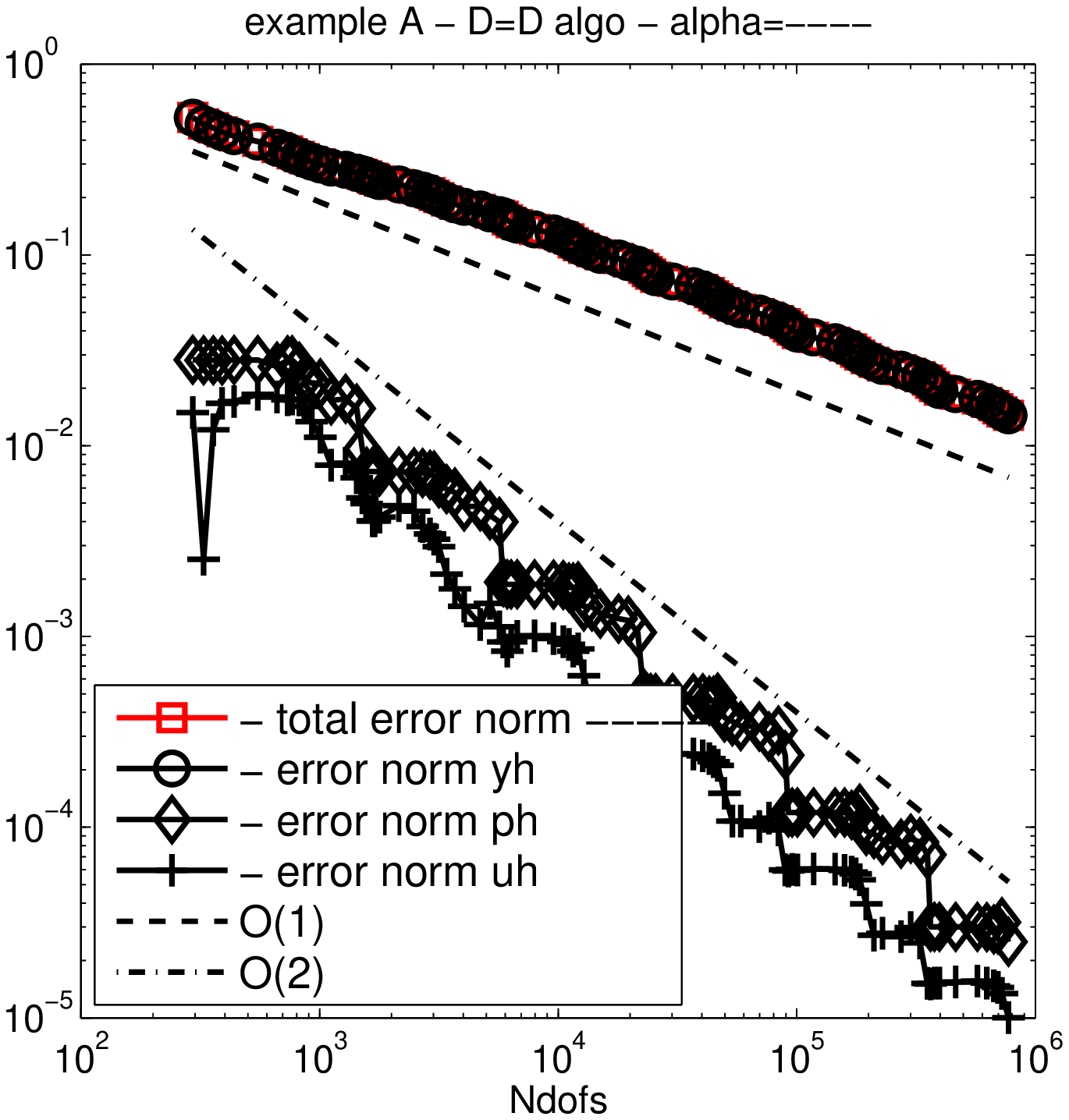}}}
\psfrag{example A - D=D algo - alpha=----}{\Huge{Example 2 - $\alpha=1$}}
\subfigure[]{\scalebox{.23}{\includegraphics{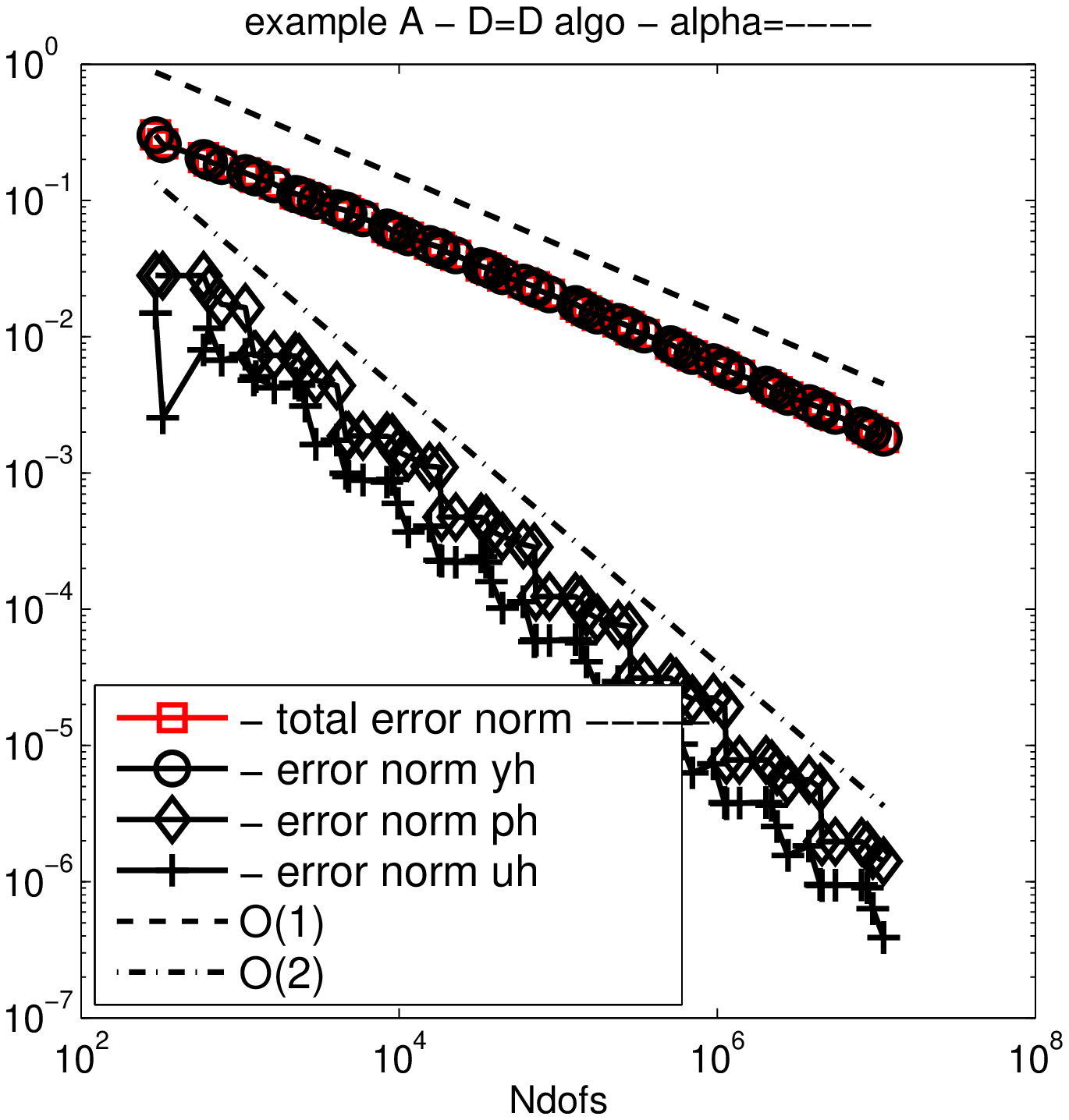}}}
\psfrag{example A - D=D algo - alpha=----}{\Huge{Example 2 - $\alpha=1.5$}}
\subfigure[]{\scalebox{.23}{\includegraphics{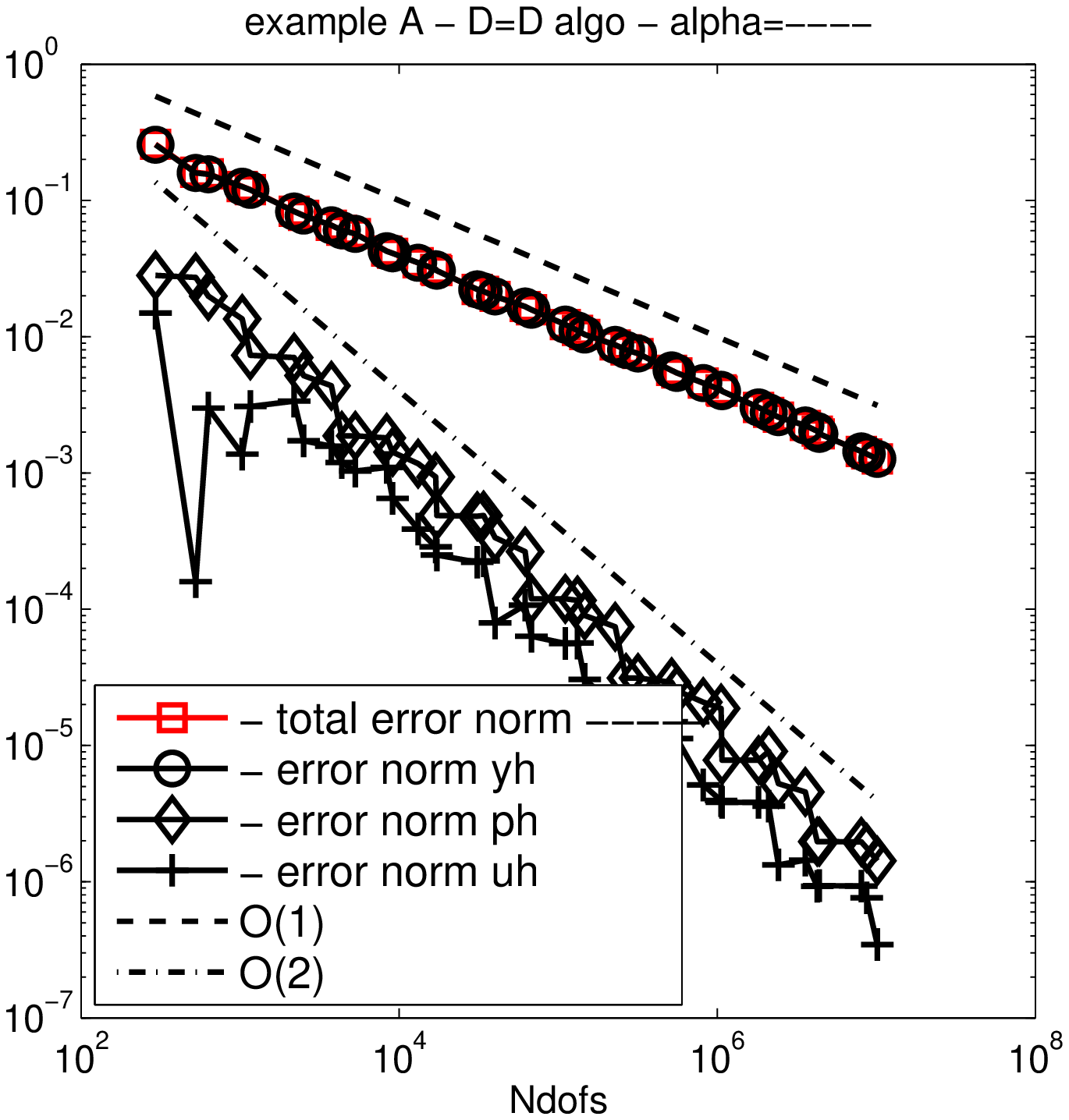}}}
\psfrag{example A - D=D algo - alpha=----}{\Huge{Example 2 - $\alpha=1.9$}}
\subfigure[]{\scalebox{.23}{\includegraphics{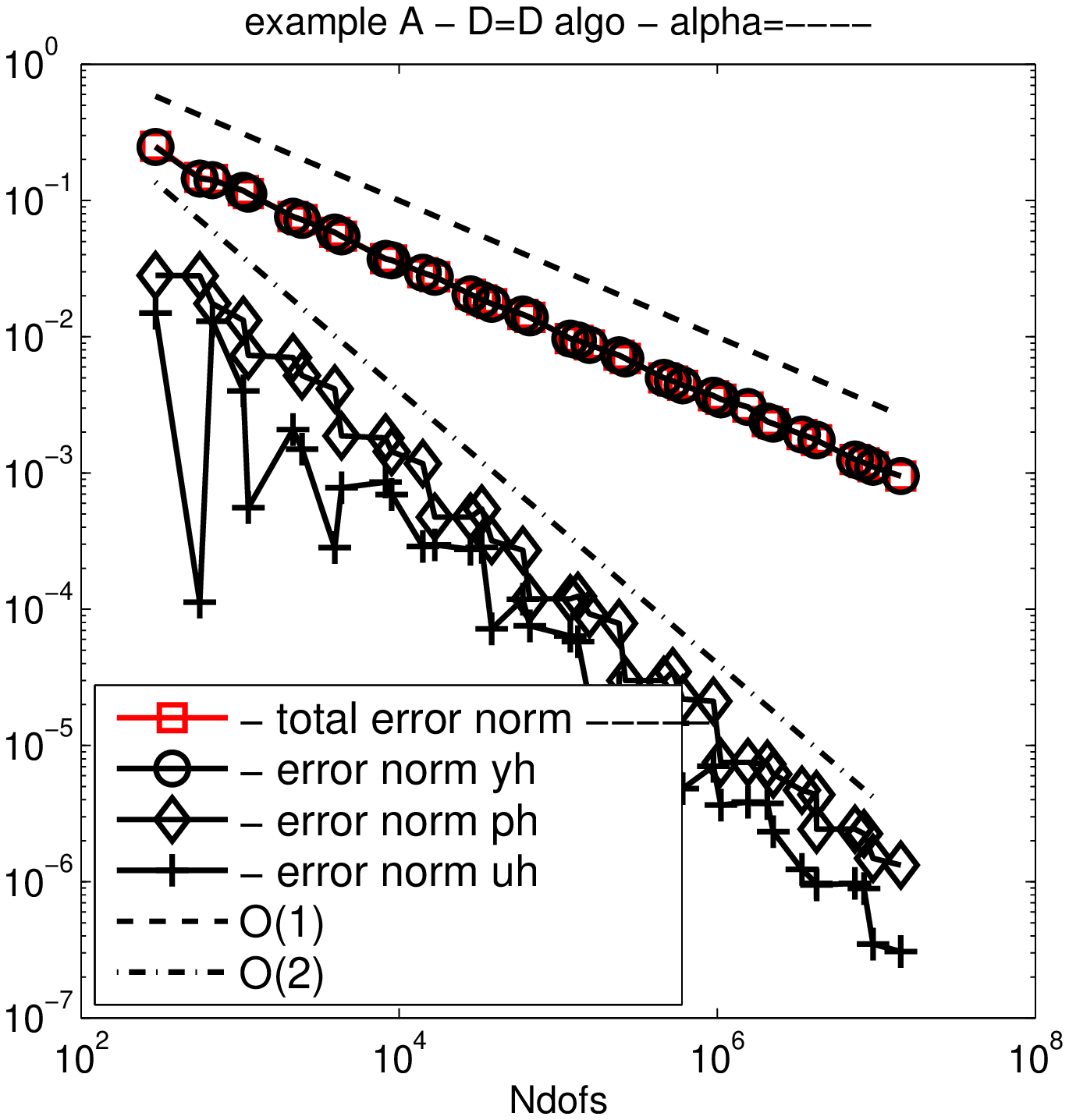}}}
\psfrag{example A - D=D algo - alpha=----}{\Huge{Example 2 - $\alpha=0.1$}}
\subfigure[]{\scalebox{.23}{\includegraphics{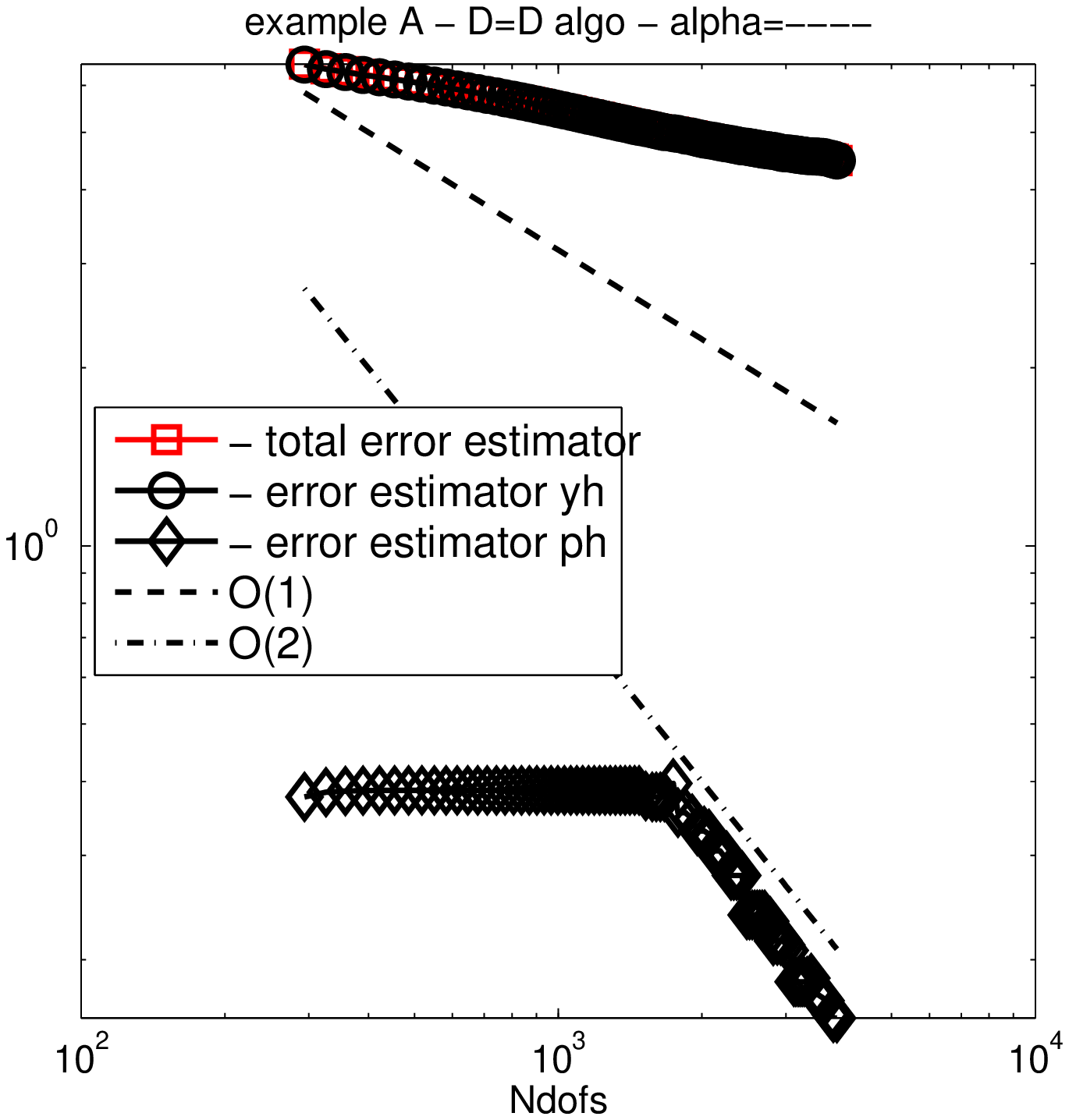}}}
\psfrag{example A - D=D algo - alpha=----}{\Huge{Example 2 - $\alpha=0.5$}}
\subfigure[]{\scalebox{.23}{\includegraphics{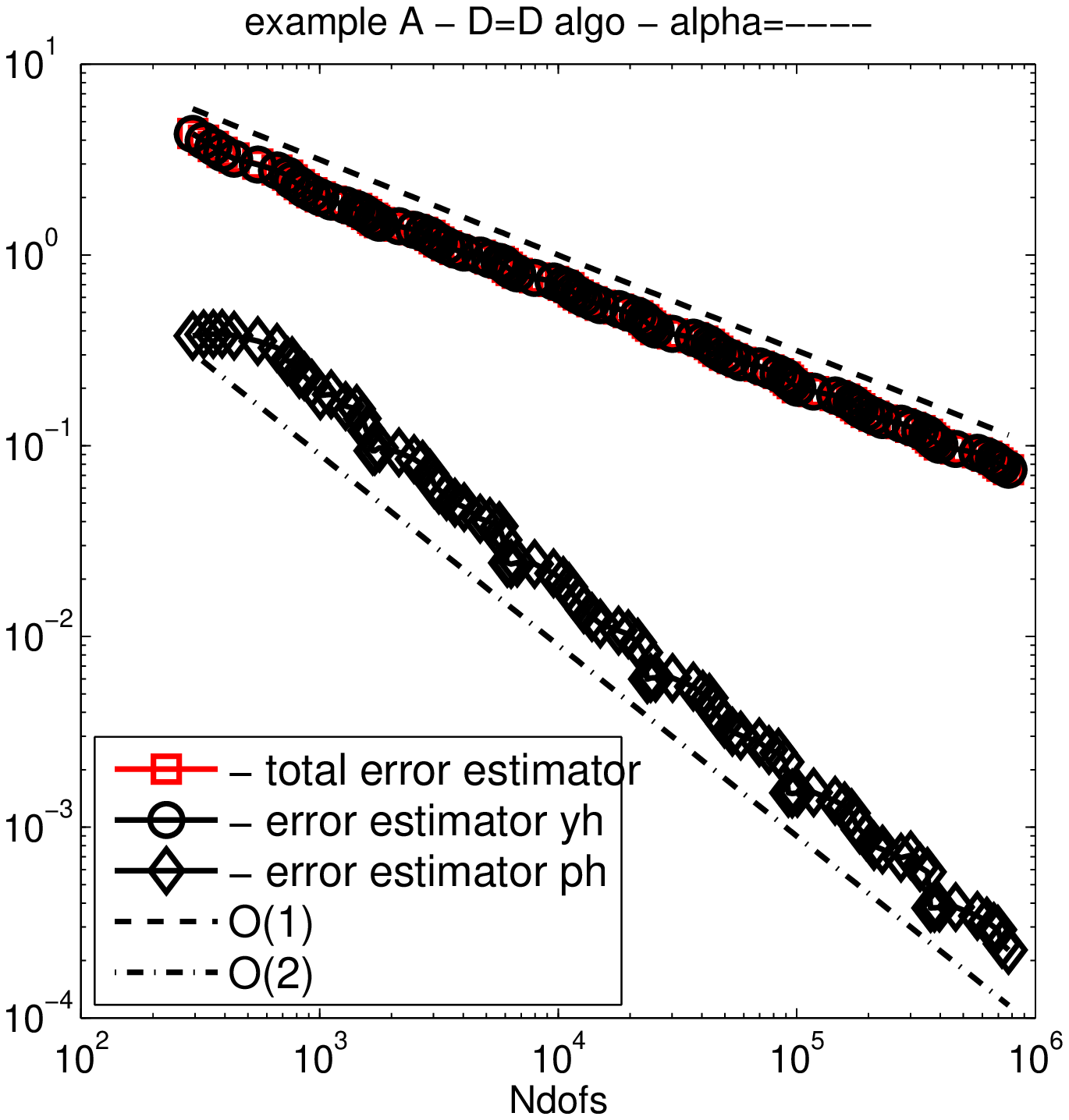}}}
\psfrag{example A - D=D algo - alpha=----}{\Huge{Example 2 - $\alpha=1$}}
\subfigure[]{\scalebox{.23}{\includegraphics{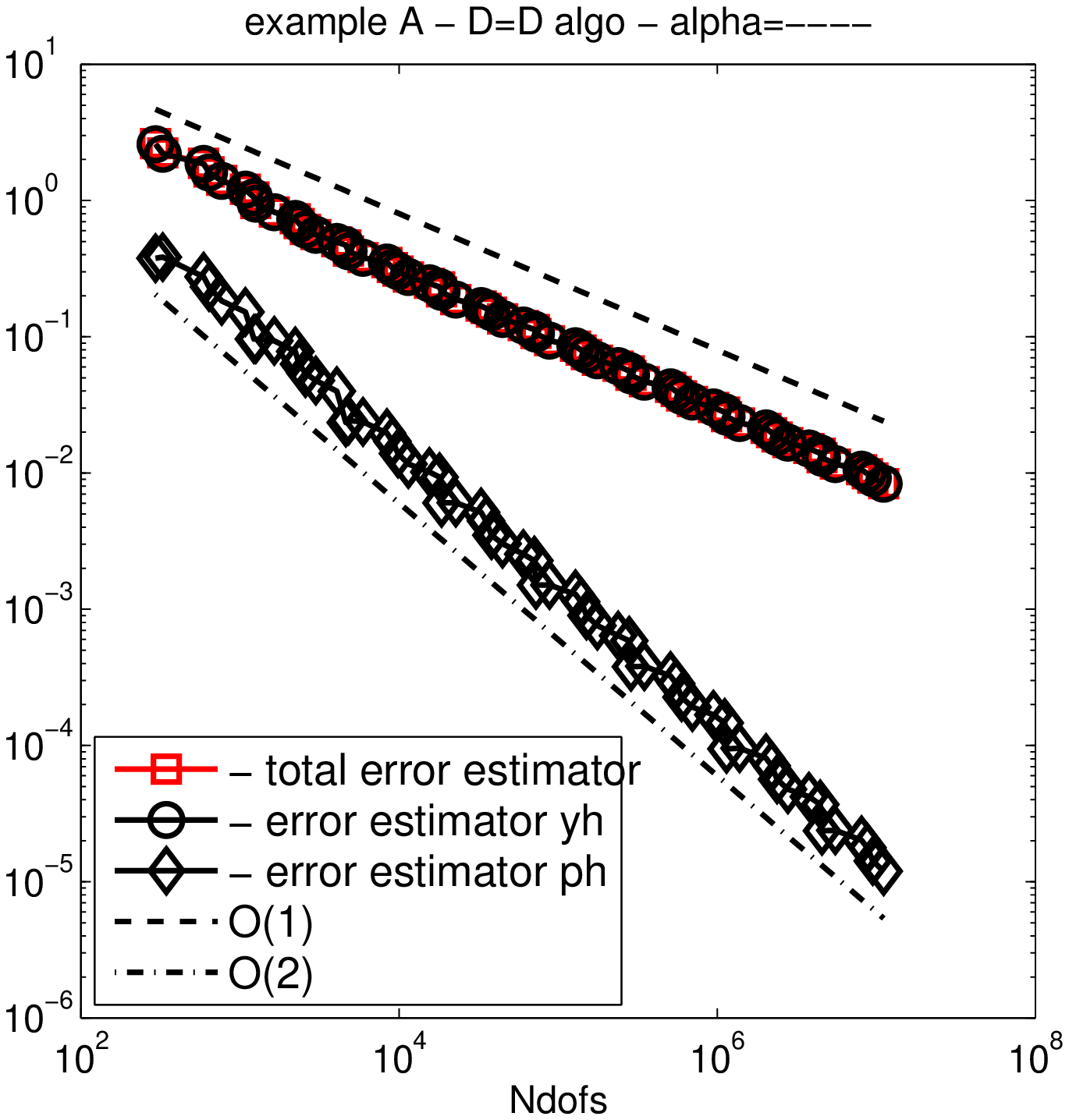}}}
\psfrag{example A - D=D algo - alpha=----}{\Huge{Example 2 - $\alpha=1.5$}}
\subfigure[]{\scalebox{.23}{\includegraphics{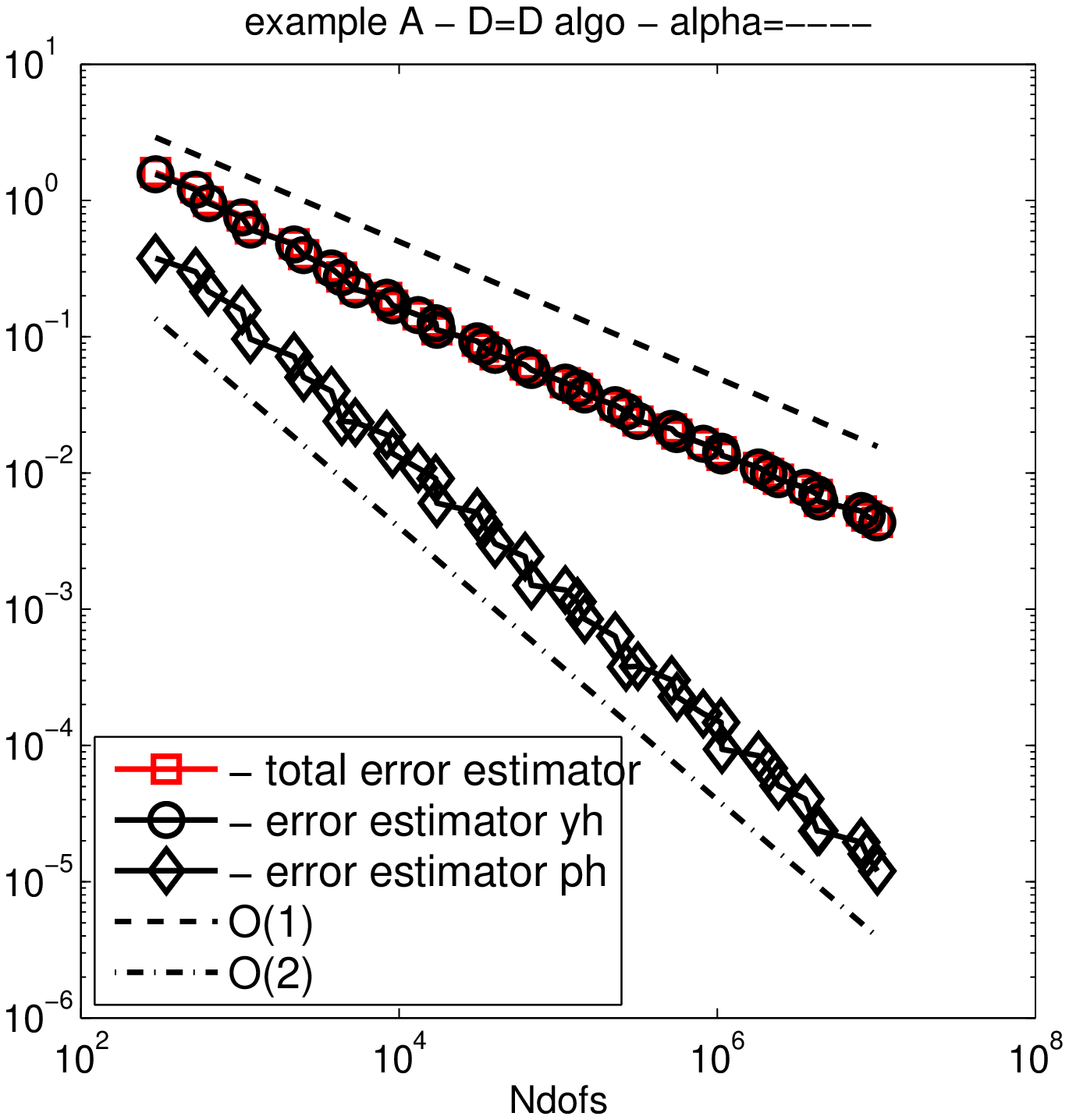}}}
\psfrag{example A - D=D algo - alpha=----}{\Huge{Example 2 - $\alpha=1.9$}}
\subfigure[]{\scalebox{.23}{\includegraphics{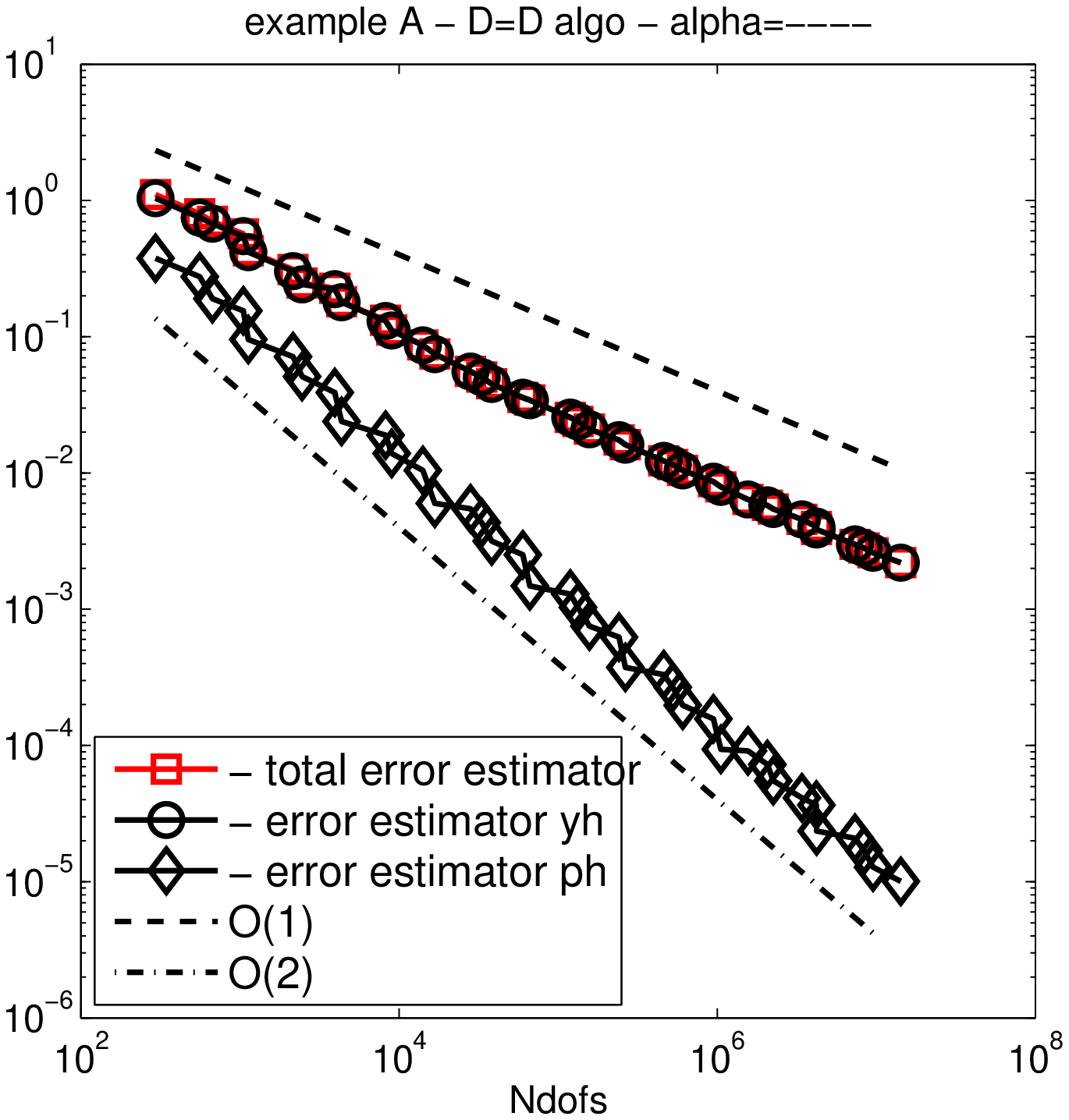}}}
\psfrag{example A - total error estimator-eficiencia}{\Huge Example 2 - $\E_{\mathsf{ocp};\T}/\|(e_{\bar{\ysf}},e_{\bar{\psf}},e_{\bar{\usf}})\|_{\Omega}$}
\psfrag{- total alpha-1}{\Huge $\alpha=0.1$}
\psfrag{- total alpha-2}{\Huge $\alpha=0.5$}
\psfrag{- total alpha-3}{\Huge $\alpha=1$}
\psfrag{- total alpha-4}{\Huge $\alpha=1.5$}
\psfrag{- total alpha-5}{\Huge $\alpha=1.9$}
\subfigure[]{\scalebox{.23}{\includegraphics{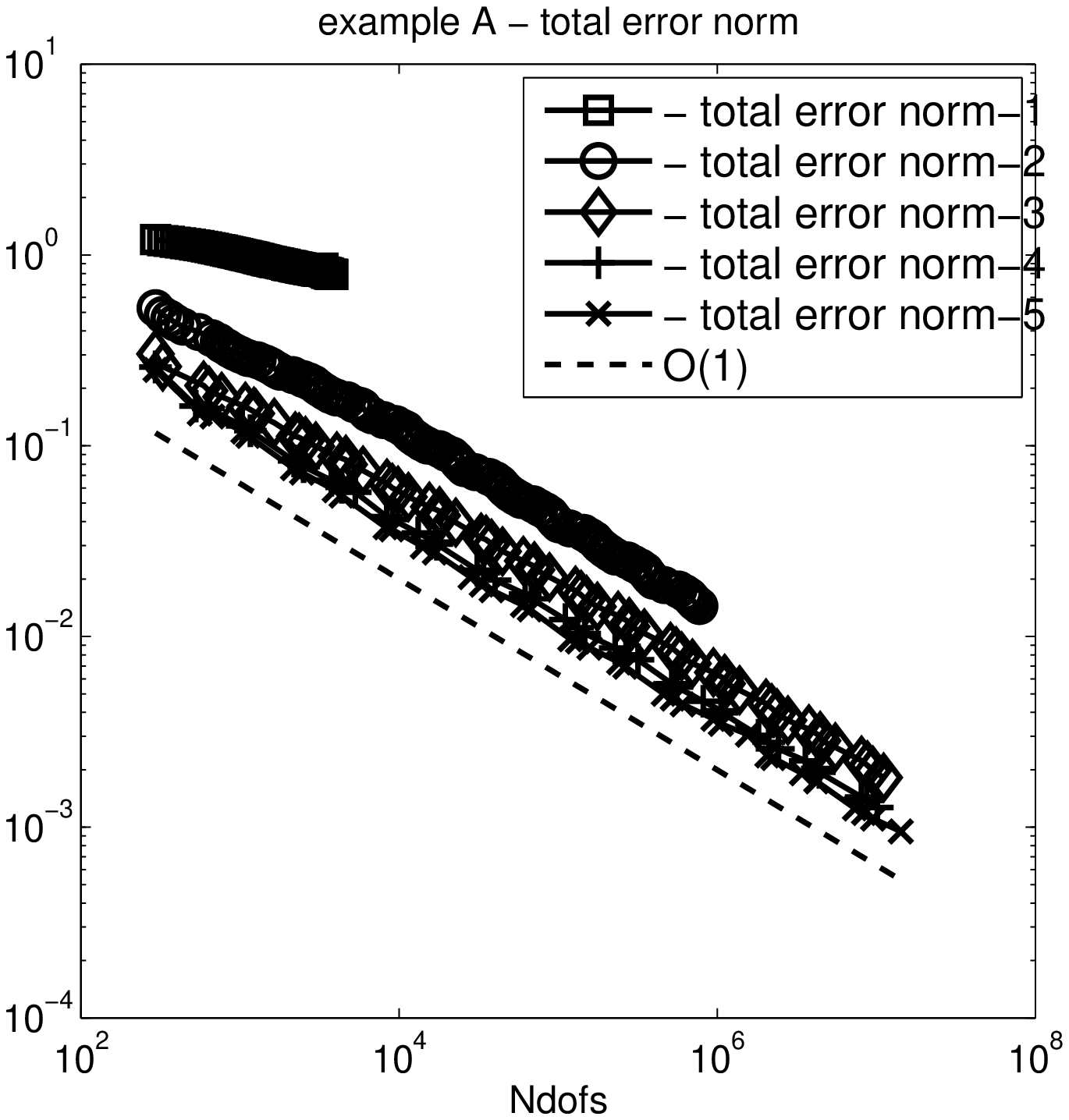}}}
\subfigure[]{\scalebox{.23}{\includegraphics{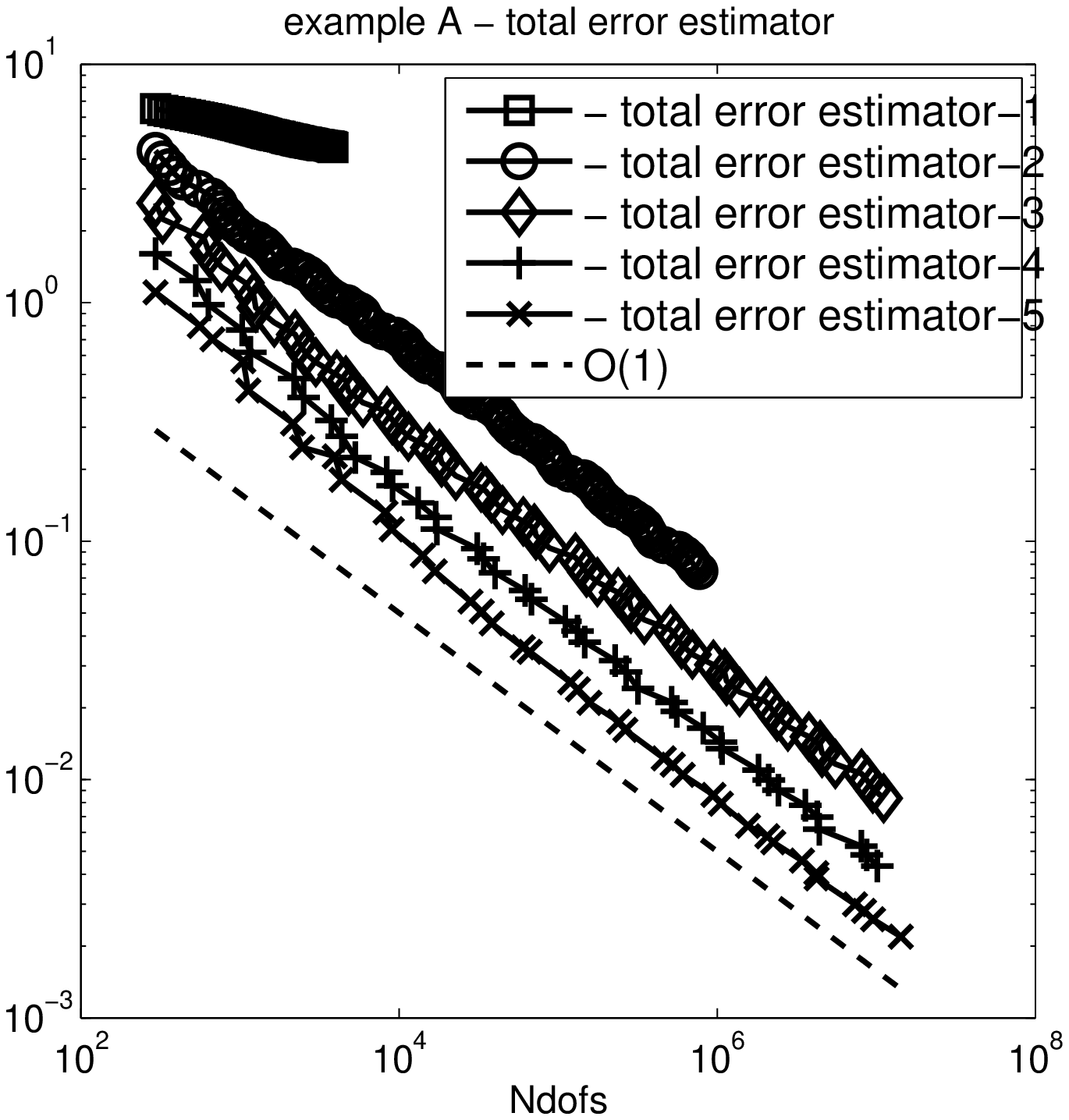}}}
\subfigure[]{\scalebox{.23}{\includegraphics{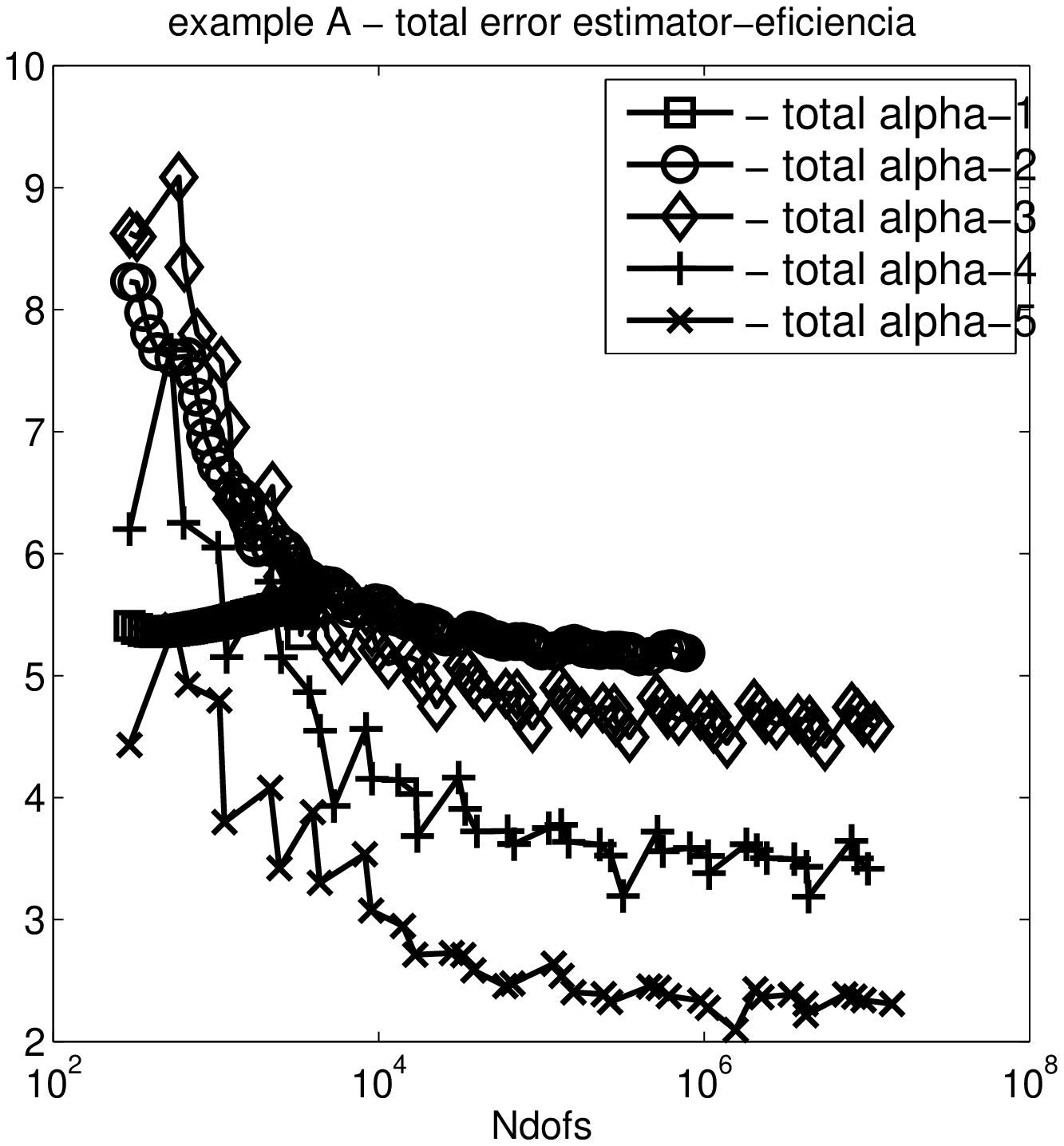}}}
\end{center}
\caption{Example 2: Experimental rates of convergence considering $\alpha\in\{0.1,0.5,1,1.5,1.9\}$, for the total error $\|(e_{\bar{\ysf}},e_{\bar{\psf}},e_{\bar{\usf}})\|_{\Omega}$ and its contributions $\|\nabla(\bar{\ysf}-\bar{\ysf}_{\T})\|_{L^{2}(\rho,\Omega)}$, $\|\bar{\psf}-\bar{\psf}_{\T}\|_{L^{\infty}(\Omega)}$ and $\|\bar{\usf}-\bar{\usf}_{\T}\|_{\mathbb{R}^{l}}$ (a)--(e); error estimator $\E_{\mathsf{ocp};\T}$ and its contributions $\E_{\ysf}$ and $\E_{\psf}$ (f)--(j); comparison of the total error $\|(e_{\bar{\ysf}},e_{\bar{\psf}},e_{\bar{\usf}})\|_{\Omega}$ for the different values of $\alpha$ (k); comparison of the error estimator $\E_{\mathsf{ocp};\T}$ for the different values of $\alpha$ (l); and effectivity indices $\E_{\mathsf{ocp};\T}/\|(e_{\bar{\ysf}},e_{\bar{\psf}},e_{\bar{\usf}})\|_{\Omega}$ for the different values of $\alpha$ (m).}
\label{Fig:2-1}
\end{figure}

\begin{figure}[!htbp]
\psfrag{example A - D=D algo - alpha=----}{}
\psfrag{- error norm yh}{\LARGE $\|\nabla(\bar{\ysf}-\bar{\ysf}_{\T})\|_{L^{2}(\rho,\Omega)}$}
\psfrag{- error norm ph}{\LARGE  $\|\bar{\psf}-\bar{\psf}_{\T}\|_{L^{\infty}(\Omega)}$}
\psfrag{- total error estimator}{\LARGE  $\E_{\mathsf{ocp};\T}$}
\psfrag{- error estimator yh}{\LARGE  $\E_{\ysf}$}
\psfrag{- error estimator ph}{\LARGE  $\E_{\psf}$}
\psfrag{- total error norm -----}{\LARGE  $\|(e_{\bar{\ysf}},e_{\bar{\psf}},e_{\bar{\usf}})\|_{\Omega}$}
\psfrag{Ndofs}{\huge Ndof}
\psfrag{O(1)}{\LARGE  $\textrm{Ndof}^{-1/2}$}
\psfrag{O(2)}{\LARGE  $\textrm{Ndof}^{-1}$}
\psfrag{yh}{}
\psfrag{ph}{}
\psfrag{- eficiencia estimador error}{\Huge $\E_{\mathsf{ocp};\T}/\|(e_{\bar{\ysf}},e_{\bar{\psf}},e_{\bar{\usf}})\|_{\Omega}$}
\begin{center}
\subfigure[]{\scalebox{.4}{\includegraphics{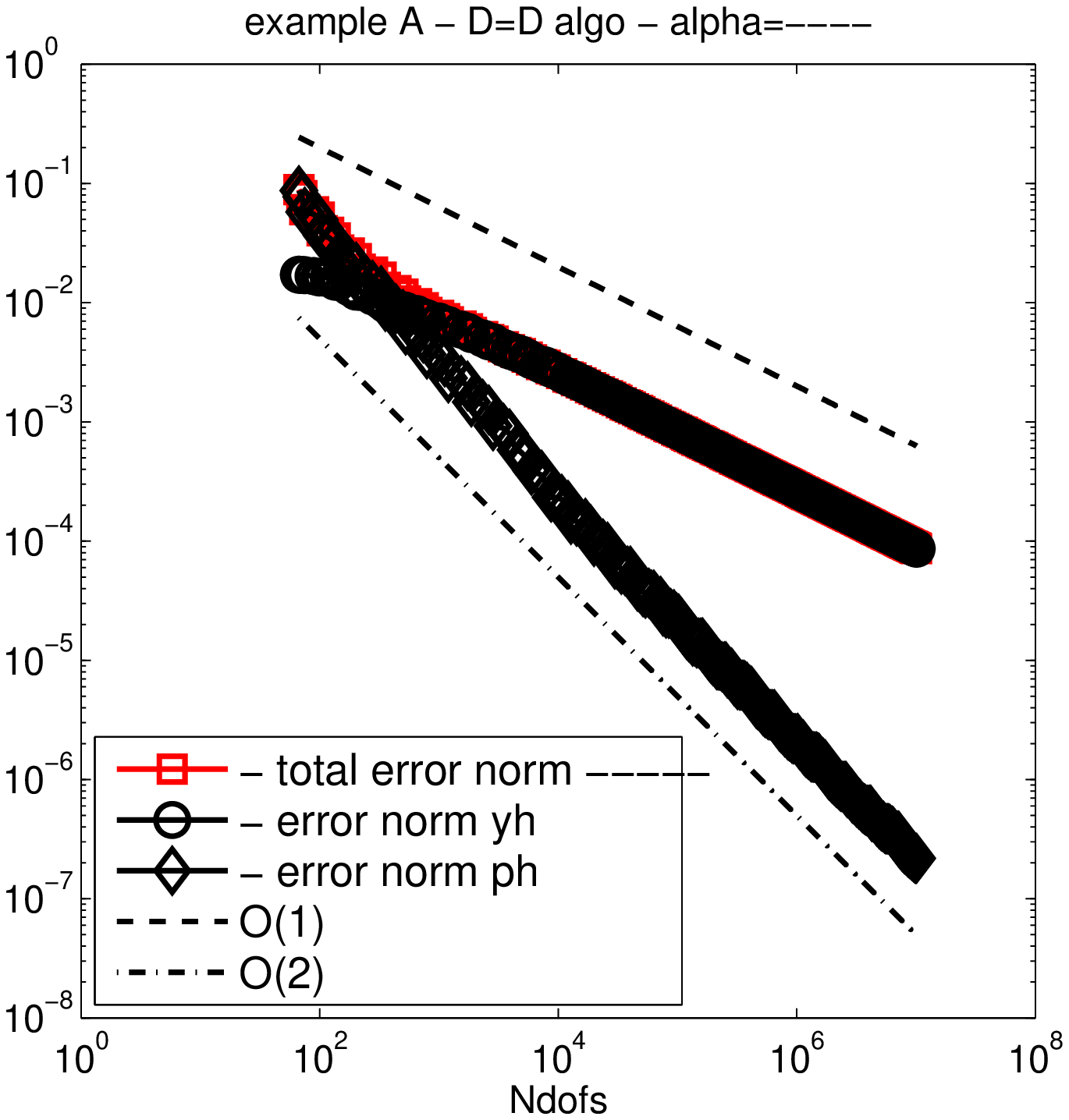}}}
\subfigure[]{\scalebox{.4}{\includegraphics{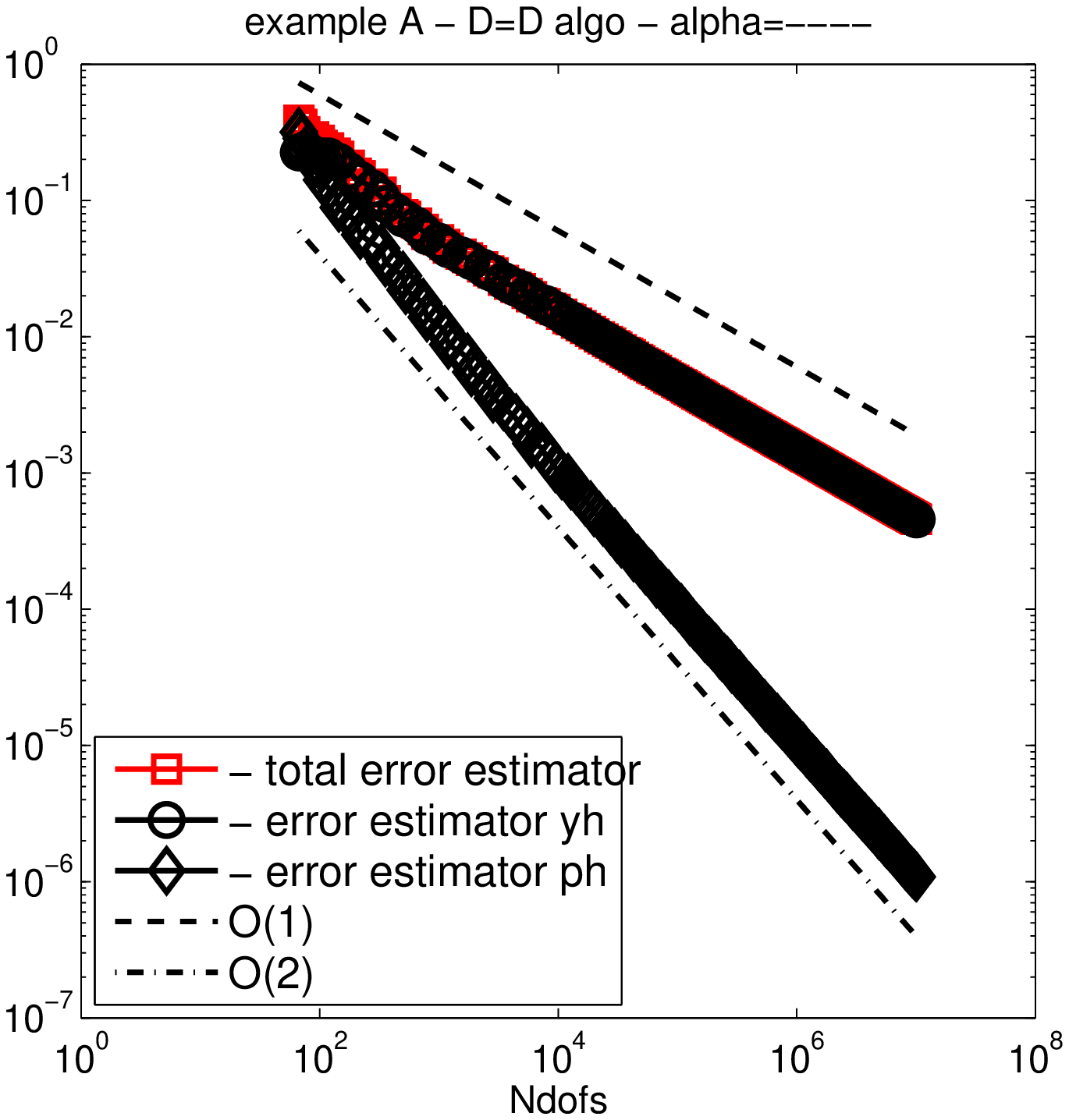}}}
\subfigure[]{\scalebox{.4}{\includegraphics{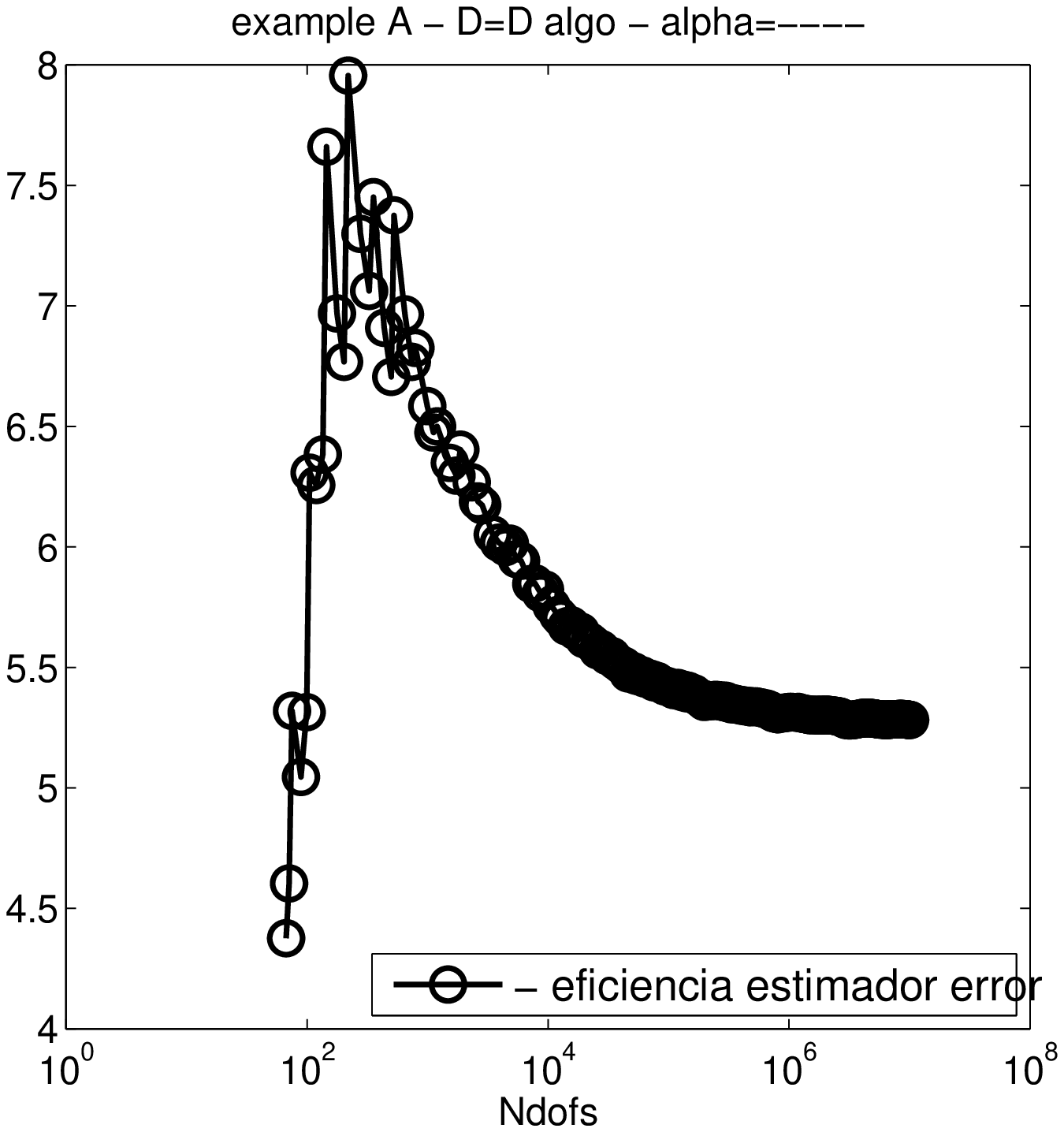}}}
\end{center}
\caption{Example 3: Experimental rates of convergence rates for the total error $\|(e_{\bar{\ysf}},e_{\bar{\psf}},e_{\bar{\usf}})\|_{\Omega}$ and its contributions $\|\nabla(\bar{\ysf}-\bar{\ysf}_{\T})\|_{L^{2}(\rho,\Omega)}$ and $\|\bar{\psf}-\bar{\psf}_{\T}\|_{L^{\infty}(\Omega)}$ (a); error estimator $\E_{\mathsf{ocp};\T}$ and its contributions $\E_{\ysf}$ and $\E_{\psf}$ (b); effectivity index $\E_{\mathsf{ocp};\T}/\|(e_{\bar{\ysf}},e_{\bar{\psf}},e_{\bar{\usf}})\|_{\Omega}$ (c).}
\label{Fig:3}
\end{figure}

We perform three examples with $n=2$ and in doing so investigate the effects of varying the parameters $\lambda$ and $\alpha$.
~\\~\\
\noindent \textbf{Example 1:} We let $\Omega=(0,1)^2$, and consider a problem with homogeneous Dirichlet boundary conditions. We set $\alpha=1.5$, $\ysf_{d}=-\sin(2\pi x)\cos(2\pi y)e^{xy}$, and
$$
D=\{(0.25,0.25),(0.75,0.25),(0.25,0.75),(0.75,0.75),(0.5,0.5)\},
$$
and consider $\asf_z=-0.5$ and $\bsf_z=1.0$ for all $z\in D$. We investigate the effect of varying the regularization parameter $\lambda$ by considering
$$
\lambda\in\{1.0,0.1,0.01,0.001,0.0001\}.
$$
The exact solutions to these problems are unknown. The results are shown in Figure \ref{Fig:1}, where we observe optimal experimental rates of convergence for the error estimator and optimal experimental decay for its contributions for all the values of the parameter $\lambda$ considered. The choice $\lambda=1$ delivers more accurate results.
~\\~\\
\noindent \textbf{Example 2:} We let $\Omega=(0,1)^2$ and consider the exact optimal state to be given by \eqref{exact_adjoint} with $\varrho_{z}=1.125$ for all $z\in D$, where
\[
D=\{(0.25,0.25),(0.75,0.25),(0.25,0.75),(0.75,0.75)\}.
\]
The exact optimal adjoint is $\bar{\psf}(x_{1},x_{2}) = -32x_1x_2(1-x_1)(1-x_2)$, $\lambda=1.0$, and $\asf_z=0.3$ and $\bsf_z=2$ for all $z\in D$. The purpose of this example is to investigate the effect of varying the exponent $\alpha$ in the Muckenhoupt weight $\rho$ defined in \eqref{eq:defofoweight}. We consider
$$
\alpha\in\{0.1,0.5,1.0,1.5,1.9\}.
$$
The results are shown in Figure \ref{Fig:2-1}. We observe that optimal experimental rates of convergence are obtained when the parameter $\alpha\in[0.5,2)$ which suggests that the meshes are being refined appropriately. However, for $\alpha=0.5$ we were unable to obtain greater accuracy than that shown in Figure \ref{Fig:2-1} by adaptively refining the mesh further. This was due to the fact that the area of some of the elements became so small that it was zero to working precision. The same situation occured when $\alpha = 0.1$.
~\\~\\
\noindent \textbf{Example 3:} We let $\asf=0.1$, $\bsf=0.9$, $\lambda=1.0$, $D=\{(0.5,0.5)\}$, $\alpha=1.0$, and $\Omega=(-1,1)^2\setminus [0,1)\times(-1,0]$ i.e., an $L$-shaped domain. The exact optimal state is given by \eqref{exact_adjoint} with $\varrho_{z}=\asf$, and the exact optimal adjoint is
$$
\bar{\psf}(x_{1},x_{2})=r^{2/3}\sin(2\theta/3),\quad\theta\in[0,3\pi/2].
$$
The results are shown in Figure \ref{Fig:3} where we observe that the total error and error estimator, together with their contributions, are decreasing at optimal rates. We note that for this example the error $\|\bar\usf-\bar\usf_\T\|_{\R^l}$ was always less than $10^{-16}$.

\subsection{Three-dimensional examples}

\begin{figure}[!htbp]
\begin{center}
\psfrag{total estimator, for alpha=1.99}{\huge The error estimator $\E_{\mathsf{ocp}}$}
\psfrag{Ndof}{\huge Ndof}
\psfrag{total estimator}{\LARGE $\E_{\mathsf{ocp}}$}
\psfrag{1d3}{\LARGE Ndof$^{-1/3}$}
\scalebox{0.4}{\includegraphics{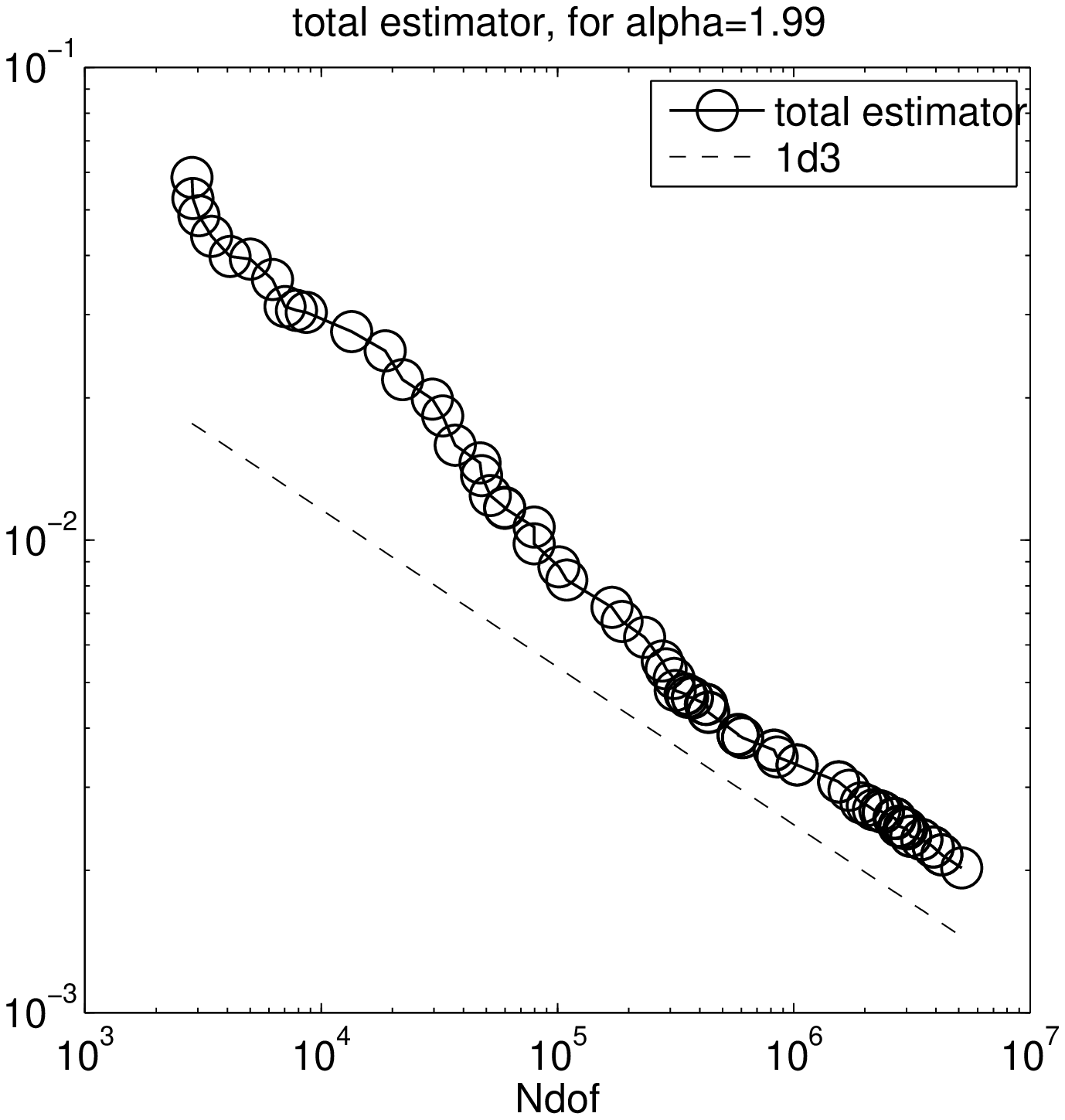}}
\psfrag{estimator contributions, for alpha=1.99}{\huge The contributions $\E_{\mathsf{y}}$ and $\E_{\mathsf{p}}$}
\psfrag{Ndof}{\huge Ndof}
\psfrag{state estimator}{\LARGE $\E_{\mathsf{y}}$}
\psfrag{adjoint estimator}{\LARGE $\E_{\mathsf{p}}$}
\psfrag{1d3}{\LARGE Ndof$^{-1/3}$}
\psfrag{2d3}{\LARGE Ndof$^{-2/3}$}
\scalebox{0.4}{\includegraphics{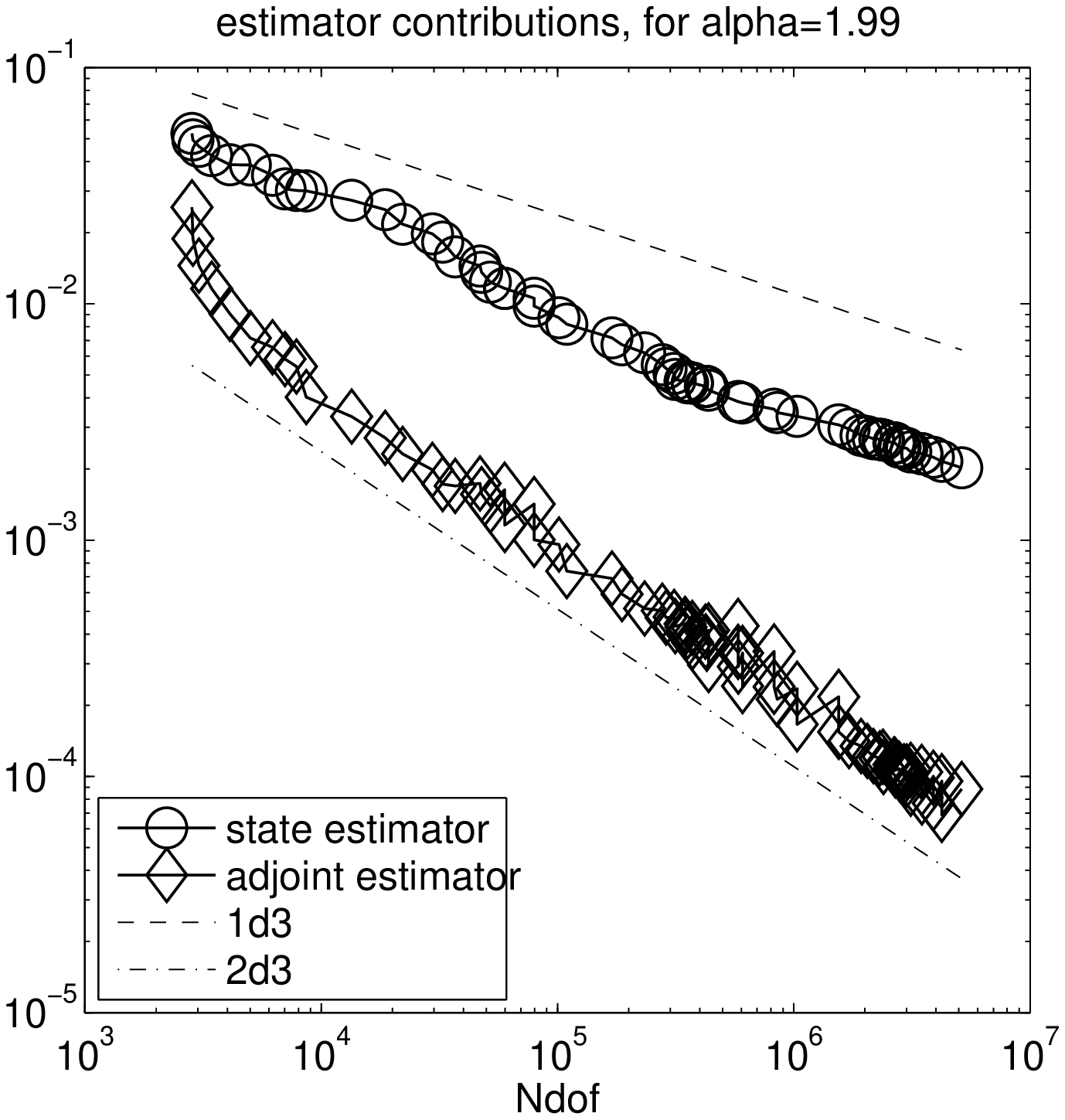}}
\end{center}
\caption{Example 4: The error estimator $\E_{\mathsf{ocp}}$ and its contributions $\E_{\mathsf{y}}$ and $\E_{\mathsf{p}}$.}
\label{Fig:Ex4}
\end{figure}

\begin{figure}[!htbp]
\begin{center}
\psfrag{the total estimator and the total error, for alpha=1.99}{\huge The error $\|(e_{\bar{\ysf}},e_{\bar{\psf}},e_{\bar{\usf}})\|_{\Omega}$ and the estimator $\E_{\mathsf{ocp}}$}
\psfrag{Ndof}{\huge Ndof}
\psfrag{total error ypu}{\LARGE $\|(e_{\bar{\ysf}},e_{\bar{\psf}},e_{\bar{\usf}})\|_{\Omega}$}
\psfrag{total estimator ypu}{\LARGE $\E_{\mathsf{ocp}}$}
\psfrag{O(1)}{\LARGE Ndof$^{-1/3}$}
\scalebox{0.4}{\includegraphics{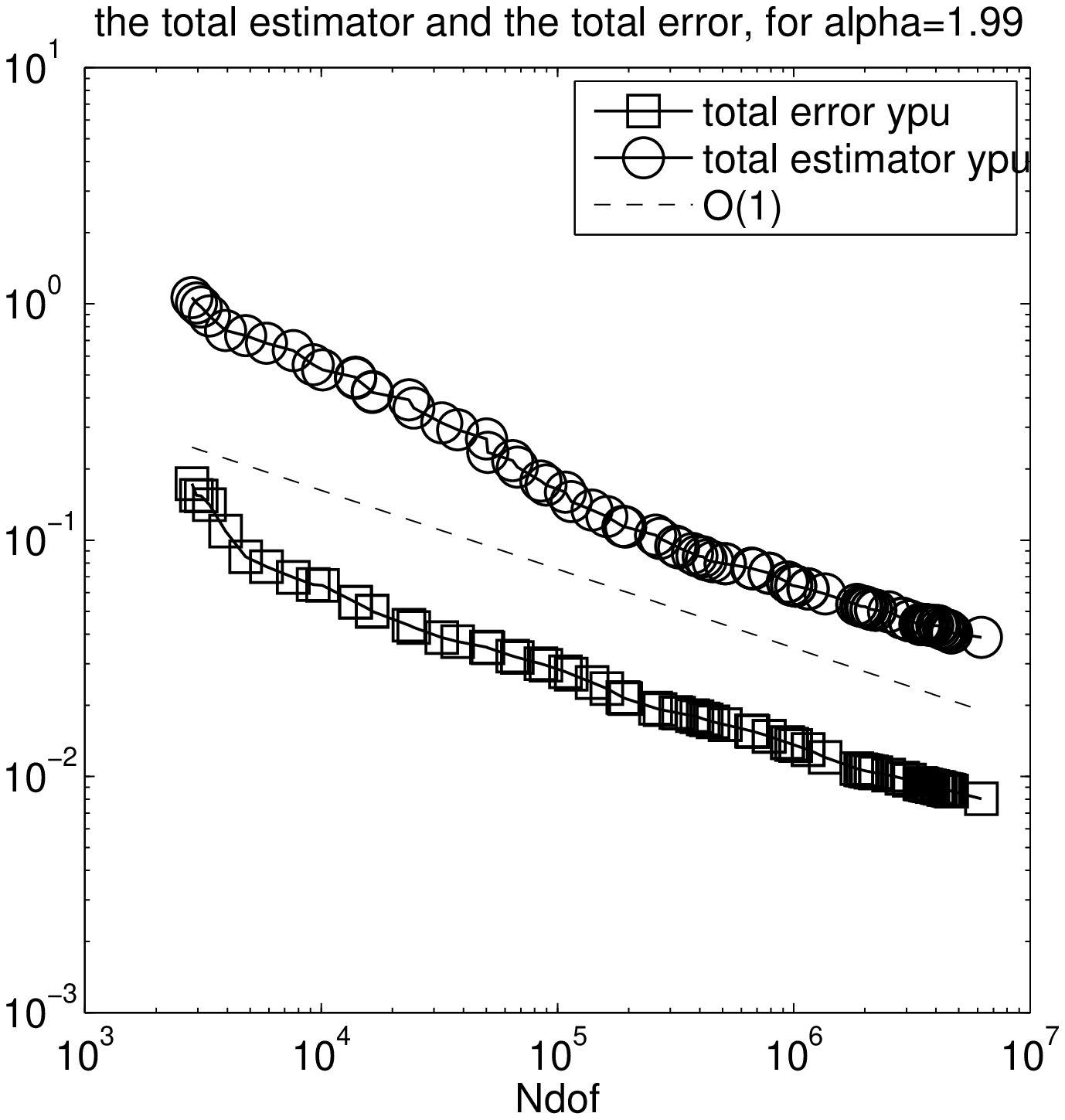}}
\psfrag{error contributions error contributions, for alpha=1.99}{\huge The contributions $\|\nabla e_{\bar{\ysf}}\|_{L^{2}(\rho,\Omega)}$, $\|e_{\bar{\psf}}\|_{L^{\infty}(\Omega)}$ and $\|e_{\bar{\usf}}\|_{\R^l}$}
\psfrag{Ndof}{\huge Ndof}
\psfrag{norm of state error}{\LARGE $\|\nabla e_{\bar{\ysf}}\|_{L^{2}(\rho,\Omega)}$}
\psfrag{norm of adjoint error}{\LARGE $\|e_{\bar{\psf}}\|_{L^{\infty}(\Omega)}$}
\psfrag{norm of control error}{\LARGE $\|e_{\bar{\usf}}\|_{\R^l}$}
\psfrag{O(1)}{\LARGE Ndof$^{-1/3}$}
\psfrag{O(2)}{\LARGE Ndof$^{-2/3}$}
\scalebox{0.4}{\includegraphics{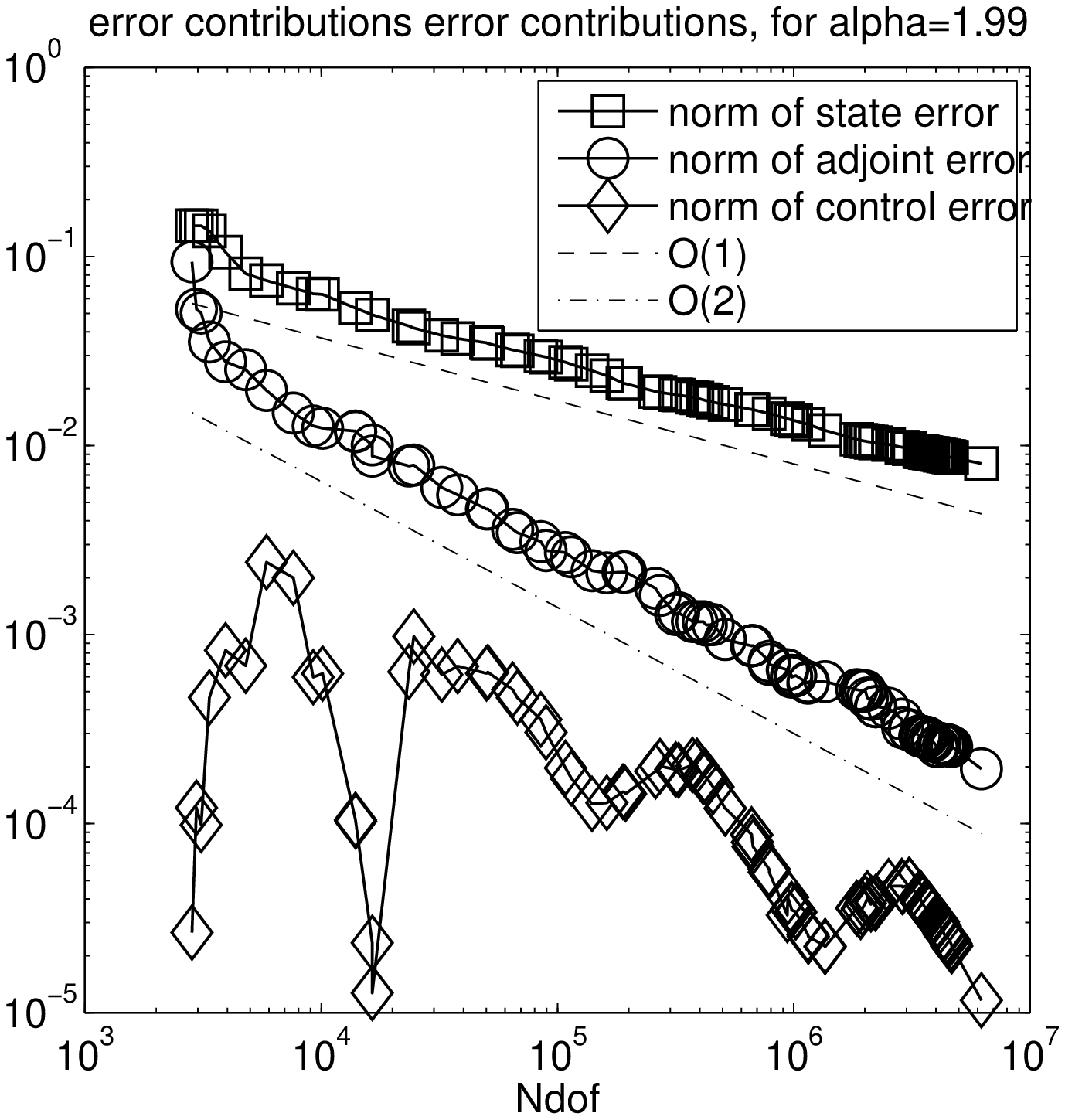}}
\psfrag{estimator contributions, for alpha=1.99}{\huge The contributions $\E_{\mathsf{y}}$ and $\E_{\mathsf{p}}$}
\psfrag{Ndof}{\huge Ndof}
\psfrag{state estimator}{\LARGE $\E_{\mathsf{y}}$}
\psfrag{adjoint estimator}{\LARGE $\E_{\mathsf{p}}$}
\psfrag{O(1)}{\LARGE Ndof$^{-1/3}$}
\psfrag{O(2)}{\LARGE Ndof$^{-2/3}$}
\scalebox{0.4}{\includegraphics{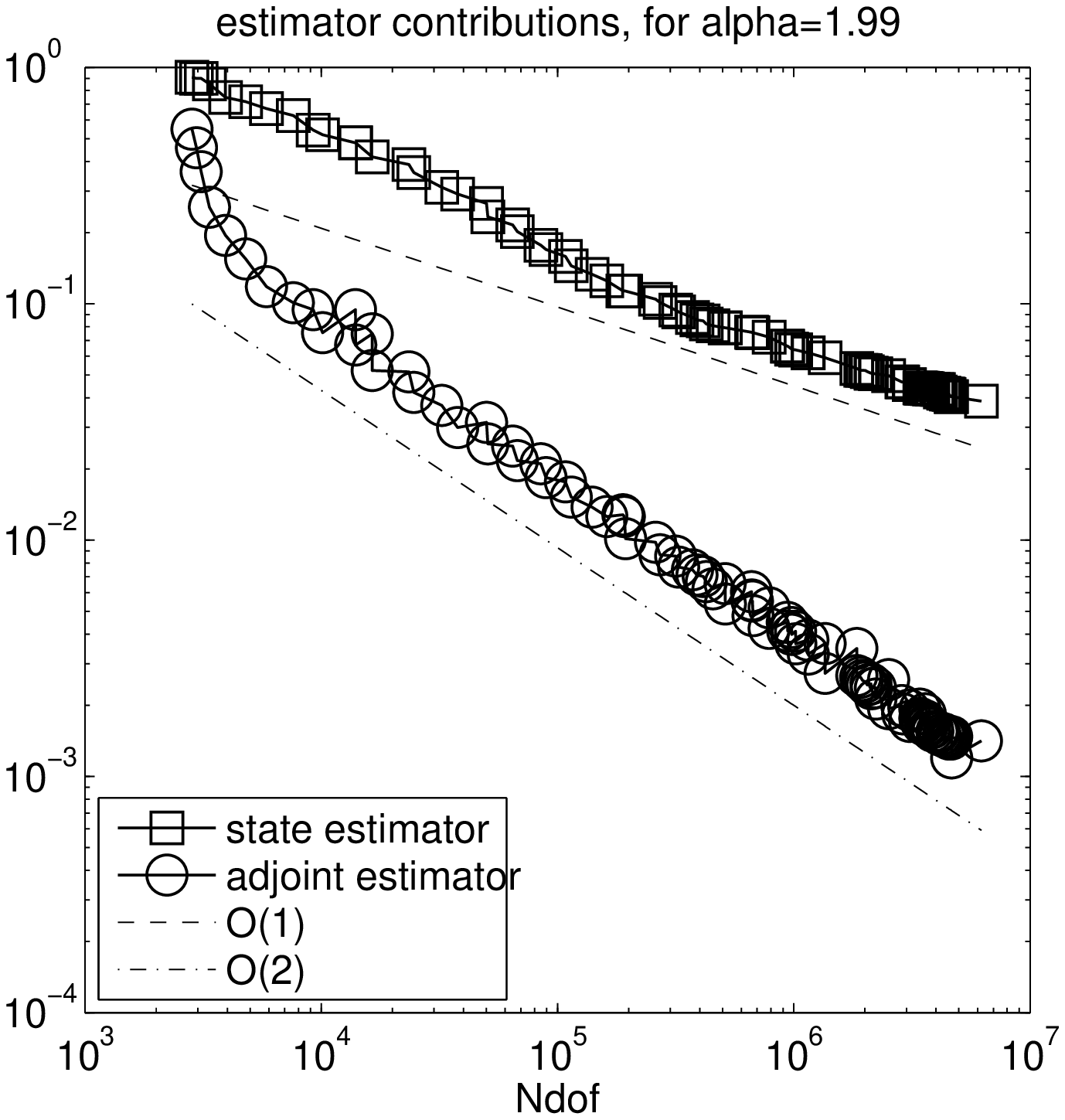}}
\psfrag{the total effectivity, for alpha=1.99}{\huge The effectivity index $\E_{\mathsf{ocp}}/\|(e_{\bar{\ysf}},e_{\bar{\psf}},e_{\bar{\usf}})\|_{\Omega}$}
\psfrag{Ndof}{\huge Ndof}
\psfrag{total effectivity index ede}{\LARGE $\E_{\mathsf{ocp}}/\|(e_{\bar{\ysf}},e_{\bar{\psf}},e_{\bar{\usf}})\|_{\Omega}$}
\psfrag{O(1)}{\LARGE Ndof$^{-1/3}$}
\scalebox{0.4}{\includegraphics{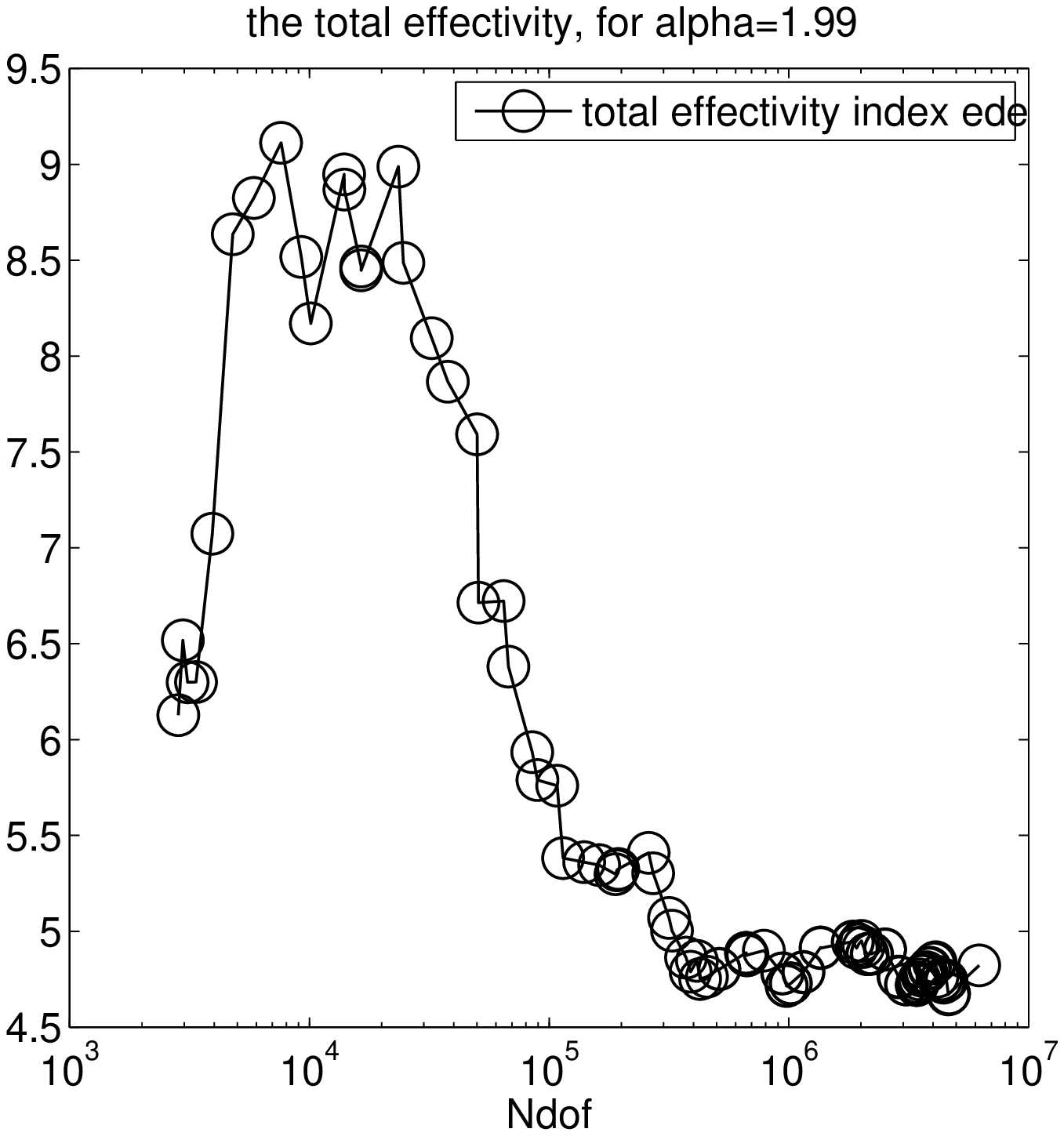}}
\end{center}
\caption{Example 5: The error $\|(e_{\bar{\ysf}},e_{\bar{\psf}},e_{\bar{\usf}})\|_{\Omega}$ and error estimator $\E_{\mathsf{ocp}}$, their contributions $\|\nabla e_{\bar{\ysf}}\|_{L^{2}(\rho,\Omega)}$, $\|e_{\bar{\psf}}\|_{L^{\infty}(\Omega)}$, $\|e_{\bar{\usf}}\|_{\R^l}$, $\E_{\mathsf{y}}$ and $\E_{\mathsf{p}}$, and the effectivity index $\E_{\mathsf{ocp}}/\|(e_{\bar{\ysf}},e_{\bar{\psf}},e_{\bar{\usf}})\|_{\Omega}$.}
\label{Fig:Ex5}
\end{figure}

\begin{figure}[!htbp]
\begin{center}
\psfrag{total estimator, for alpha=1.99}{\huge The error estimator $\E_{\mathsf{ocp}}$}
\psfrag{Ndof}{\huge Ndof}
\psfrag{total estimator}{\LARGE $\E_{\mathsf{ocp}}$}
\psfrag{1d3}{\LARGE Ndof$^{-1/3}$}
\scalebox{0.4}{\includegraphics{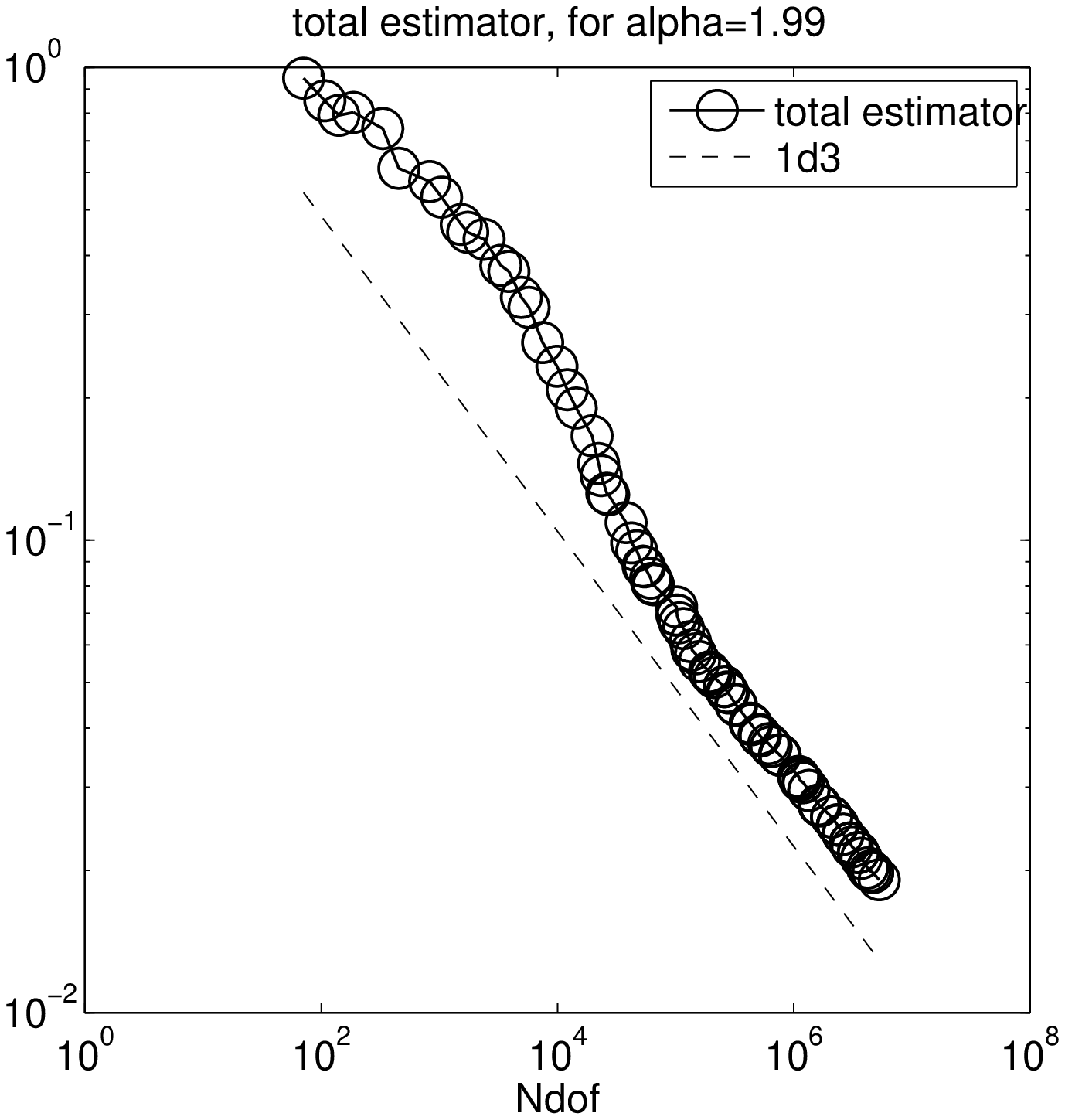}}
\psfrag{estimator contributions, for alpha=1.99}{\huge The contributions $\E_{\mathsf{y}}$ and $\E_{\mathsf{p}}$}
\psfrag{Ndof}{\huge Ndof}
\psfrag{state estimator}{\LARGE $\E_{\mathsf{y}}$}
\psfrag{adjoint estimator}{\LARGE $\E_{\mathsf{p}}$}
\psfrag{1d3}{\LARGE Ndof$^{-1/3}$}
\psfrag{2d3}{\LARGE Ndof$^{-2/3}$}
\scalebox{0.4}{\includegraphics{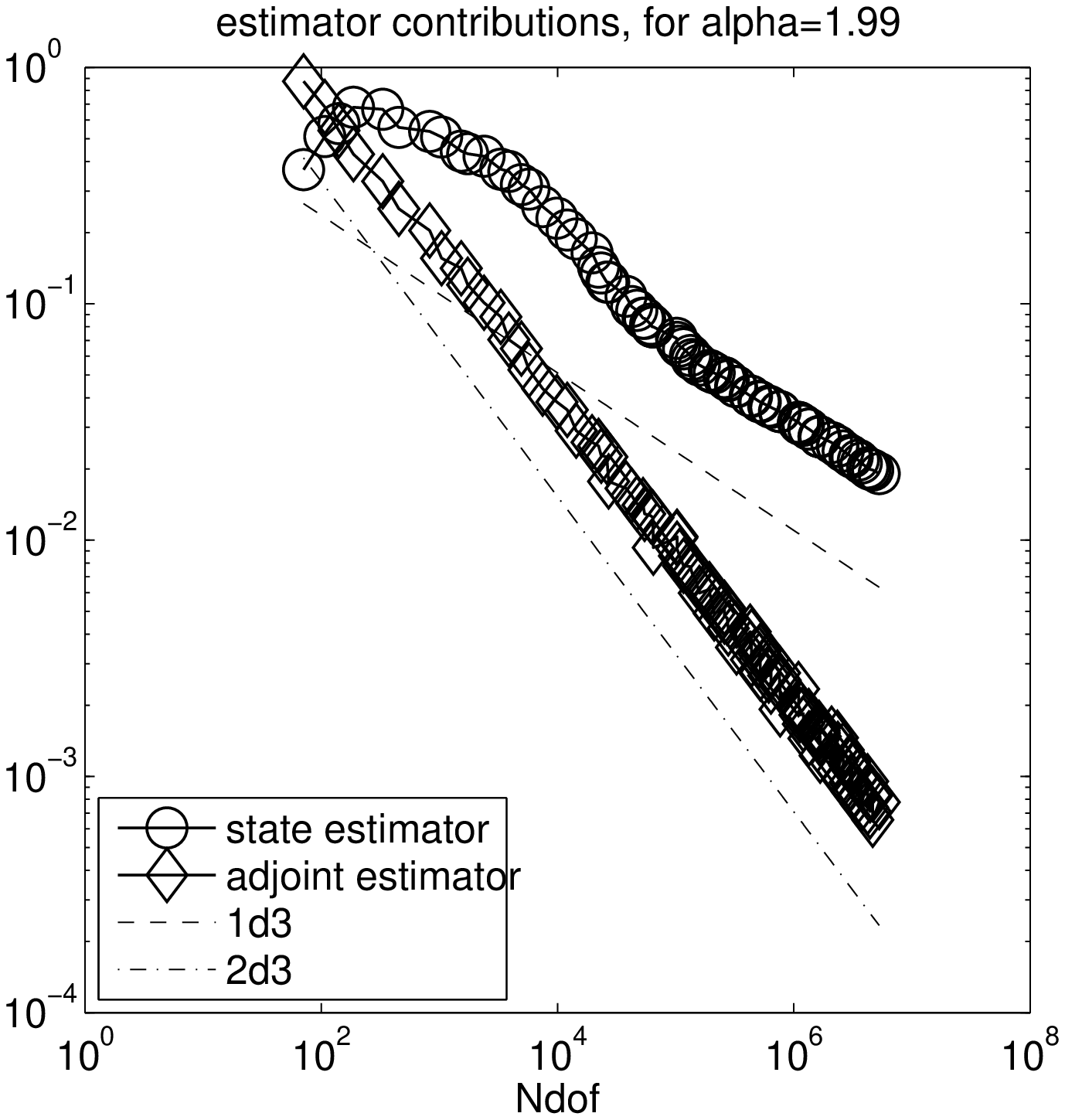}}
\end{center}
\caption{Example 6: The error estimator $\E_{\mathsf{ocp}}$ and its contributions $\E_{\mathsf{y}}$ and $\E_{\mathsf{p}}$.}
\label{Fig:Ex6}
\end{figure}

We perform three examples with $n=3$. In all of the three-dimensional examples we consider $\alpha=1.99$.
~\\~\\
\noindent \textbf{Example 4:} We let $\Omega=(0,1)^3$, and consider a problem with homogeneous Dirichlet boundary conditions, whose exact solutions are not known. We set $\lambda=1$, $\ysf_{d}=-\sin(2\pi x)\sin(2\pi y)\sin(2\pi z)e^{xyz}$, and
$$
D=\{(0.25,0.25,0.25),(0.75,0.75,0.75)\},
$$
and consider $\asf_z=-0.5$ and $\bsf_z=1$ for all $z\in D$. The results are shown in Figure \ref{Fig:Ex4}. We observe that the estimator and its contributions are decreasing at the optimal rates.
~\\~\\
\noindent \textbf{Example 5:} We let $\Omega=(0,1)^3$, $\lambda=1$, and $D=\{(0.25,0.25,0.25),(0.75,0.75,0.75)\}$, and consider $\asf_z=0$ and $\bsf_z=0.25$ for all $z\in D$. The exact optimal state is given by \eqref{exact_adjoint} with $\varrho_{(0.25,0.25,0.25)}=27/256$ and $\varrho_{(0.75,0.75,0.75)}=0.25$, and the exact optimal adjoint is
$$
\bar{\psf}(x_{1},x_{2},x_3) = -64x_1x_2x_3^{2}(1-x_1)(1-x_2)(1-x_3).
$$
The results are shown in Figure \ref{Fig:Ex5}. We observe that the error and the estimator, as well as their contributions, are decreasing at the optimal rates.
~\\~\\
\noindent \textbf{Example 6:} We consider a problem with homogeneous Dirichlet boundary conditions on the $L$-shaped domain
$$
\Omega=((-\sqrt{2},\sqrt{2})\times(-\sqrt{2},\sqrt{2})\times(0,1))\setminus([0,\sqrt{2})\times[0,\sqrt{2})\times(0,1)),
$$
whose exact solutions are not known. We set $\ysf_{d}=1$, $\asf=-1$, $\bsf=1$, $\lambda=1$, and $D=\{(0.5,0.5,0.5)\}$. The results are shown in Figure \ref{Fig:Ex6} where we can observe that the estimator and its contributions are decreasing at the optimal rates.

\bibliographystyle{plain}
\bibliography{biblio}

\end{document}